\documentclass[twoside,11pt]{article}

\usepackage{amsmath,amsthm,amssymb}

\usepackage[abbrvbib, preprint, nohyperref]{jmlr2e}

\usepackage{enumitem}
\usepackage{natbib}
\usepackage{mathrsfs}
\usepackage{graphicx,psfrag,epsf}
\usepackage{epstopdf}
\usepackage{epsfig}
\usepackage{subcaption}
\usepackage{bbm}	
\usepackage{color}
\usepackage{mathtools} 
\usepackage{tikz}
\usetikzlibrary{arrows}
\usepackage{booktabs}

\usepackage[hyperindex,breaklinks,colorlinks,citecolor=blue,urlcolor=black]{hyperref} 

\theoremstyle{plain}
\newtheorem{mylemma}{Lemma}
\newtheorem{mytheorem}{Theorem}

\newtheorem{mycorollary}{Corollary}

\theoremstyle{definition} 
\newtheorem{mydefinition}{Definition} 

\newtheorem{myremark}{Remark}

\makeatletter
\newcommand{\mylabel}[2]{#2\def\@currentlabel{#2}\label{#1}}
\makeatother

\newenvironment{myproof}{\noindent{\bfseries Proof.}}{\hfill\BlackBox \\

}
\newenvironment{myprooftitle}[1][\myproofname]{\noindent{\bfseries #1.}}{\hfill\BlackBox \\

}

\setitemize{itemsep=-3pt}
\setenumerate{itemsep=-3pt}

\usepackage{lastpage}
\jmlrheading{24}{2023}{1-\pageref{LastPage}}{4/22; Revised 6/23}{11/23}{22-0382}{Seonghyun Jeong and Veronika Ro{\v{c}}kov{\'a}}
\ShortHeadings{The art of BART}{Jeong and Ro{\v{c}}kov{\'a}}

\firstpageno{1}

\begin{document}

\title{The Art of BART: Minimax Optimality over Nonhomogeneous Smoothness in High Dimension}

\author{\name Seonghyun Jeong \email sjeong@yonsei.ac.kr \\
       \addr Department of Statistics and Data Science\\
        Department of Applied Statistics\\
       Yonsei University\\
       Seoul 03722, Republic of Korea
       \AND
       \name Veronika Ro{\v{c}}kov{\'a} \email veronika.rockova@chicagobooth.edu \\
       \addr Booth School of Business\\
       University of Chicago\\
       Chicago, IL 60637, USA}

\editor{Daniel Roy}

\maketitle

\begin{abstract}
Many asymptotically minimax procedures for function estimation often rely on somewhat arbitrary and restrictive assumptions such as isotropy or spatial homogeneity. This work enhances the theoretical understanding of Bayesian additive regression trees under substantially relaxed smoothness assumptions. We provide a comprehensive study of asymptotic optimality and posterior contraction of Bayesian forests when the regression function has anisotropic smoothness that possibly varies over the function domain. The regression function can also be possibly discontinuous. We introduce a new class of sparse {\em piecewise heterogeneous anisotropic} H\"{o}lder functions and derive their minimax lower bound of estimation in high-dimensional scenarios under the $L_2$-loss. We then find that the Bayesian tree priors, coupled with a Dirichlet subset selection prior for sparse estimation in high-dimensional scenarios, adapt to unknown heterogeneous smoothness, discontinuity, and sparsity. These results show that Bayesian forests are uniquely suited for more general estimation problems that would render other default machine learning tools, such as Gaussian processes, suboptimal. Our numerical study shows that Bayesian forests often outperform other competitors such as random forests and deep neural networks, which are believed to work well for discontinuous or complicated smooth functions. Beyond nonparametric regression, we also examined posterior contraction of Bayesian forests for density estimation and binary classification using the technique developed in this study.
\end{abstract}


\section{Introduction}
\subsection{Motivation}
Many of the existing asymptotic minimaxity results for estimating regression functions are predicated on the assumption that certain smoothness conditions hold, which can  be rarely satisfied/verified  when confronted with real data. This creates a disconnect between theory and practice, limiting the scope of many theoretical results.
For example, in nonparametric regression involving multiple predictors, the assumption of {\em isotropic smoothness}  can be unnecessarily restrictive.
A more realistic scenario is when the function exerts  different degrees of smoothness in different directions and areas, with possible discontinuities that allow further flexibility.
This study is motivated by the desire to evaluate the theoretical performance of Bayesian forests, one of the workhorses of Bayesian machine learning, in such broad scenarios.

Bayesian trees and their ensembles have achieved notable empirical success in statistics and machine learning \citep{chipman1998bayesian,denison1998bayesian,chipman2010bart}. 
Relative to other Bayesian machine learning alternatives, tree-based methods require comparatively less tuning and can be scaled to higher dimensions \citep{lakshminarayanan2013top,bleich2014variable,he2019xbart}.
The popularity of Bayesian forests, such as Bayesian additive regression trees (BART), \citep{chipman2010bart} is growing rapidly  in many areas including causal inference \citep{hill2011bayesian,hahn2020bayesian}, mean-variance function estimation \citep{pratola2020heteroscedastic}, smooth function estimation  \citep{linero2018bayesian-jrssb}, variable selection \citep{bleich2014variable,linero2018bayesian-jasa}, interaction detection \citep{du2019interaction}, survival analysis \citep{sparapani2016nonparametric}, time series  \citep{taddy2011dynamic}, count and categorical data analysis \citep{murray2021log}, and density regression \citep{orlandi2021density,li2022adaptive}. For comprehensive overviews and surveys, refer to \citet{linero2017review}, \citet{tan2019bayesian}, and \citet{hill2020bayesian}.

Despite remarkable success in empirical studies, the theoretical properties of Bayesian forests remained unavailable until the emergence of recent literature \citep{rockova2020posterior,linero2018bayesian-jrssb,rockova2019theory,castillo2021uncertainty}. 
Although these pioneering findings divulge why tree-based methods perform well,
they are limited to isotropic regression function surfaces, which exhibit the same level of smoothness in every direction.
Isotropy is an archetypal assumption in theoretical studies, but it can be restrictive in real-world applications.
This assumption is particularly unattractive in higher dimensions wherein the function can behave very poorly in certain directions.

\begin{figure}[t!]
	\centering
	\begin{subfigure}[b]{5.5in}
		\includegraphics[width=5.5in]{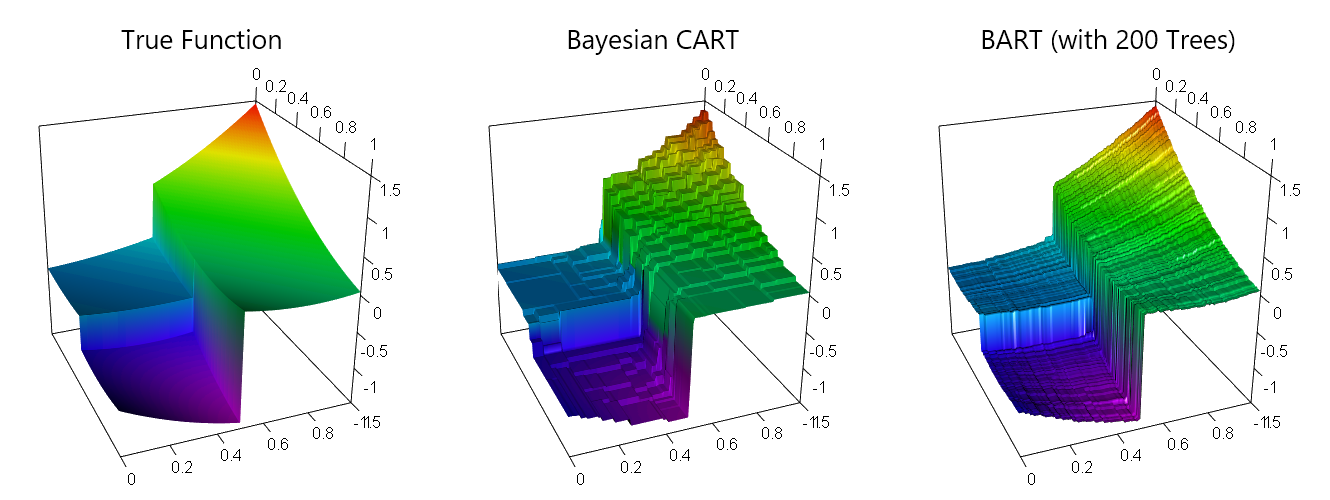}
		\caption{Piecewise smooth function estimation }
	\end{subfigure}
	\begin{subfigure}[b]{5.5in}
		\includegraphics[width=5.5in]{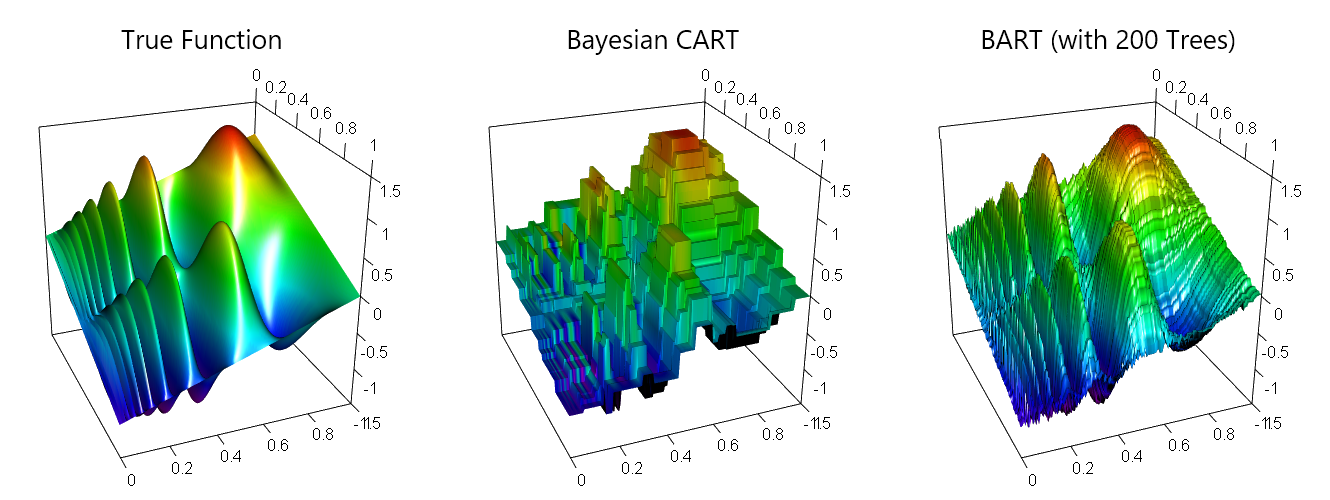}
		\caption{Doppler-type function estimation }
	\end{subfigure}
	\caption{Function estimation in nonparametric regression with complicated smoothness using Bayesian CART and BART.}
	\label{fig:bart}
\end{figure}

However, empirical evidence suggests that Bayesian forests are expected to adapt to more intricate smoothness situations. For example, Figure~\ref{fig:bart} shows that BART successfully adapts to a piecewise smooth function or a Doppler-type function.
The successful performance beyond isotropy is attributable to at least three reasons: (i)  tree methods are based on top-down recursive partitioning, wherein splits occur more often in areas where the function is locally uneven or bumpy, making the procedure spatially adaptive; (ii) the choice of coordinates for the split is data-driven, dividing the domain more often in directions in which the function is less smooth; and (iii)  tree-based learners are piecewise constant and, as such, are expected to adapt to discontinuous functions by detecting smoothness boundaries and jumps. 
These considerations naturally  create an expectation that Bayesian forests achieve optimal estimation properties in more complex function classes  {\em without any prior modification.}

\subsection{Our Contribution}

The main goal of this study is to examine optimality and posterior contraction of Bayesian forests under relaxed smoothness assumptions.
We introduce a class of functions the domain of which has been cleaved into hyper-rectangles, where each rectangular piece has its own anisotropic smoothness (with the same harmonic mean). We allow for possible discontinuities at the boundaries of the pieces.
We call this new class of functions {\em piecewise heterogeneous anisotropic functions} (see Definitions~\ref{def:hol}--\ref{def:genhol} in Section~\ref{sec:function}). We then establish an approximation theory for this general class, which blends anisotropy with spatial inhomogeneity and which, to the best of our knowledge, has not yet been pursued in the literature. Our results complement the body of existing work on  piecewise isotropic smoothness classes \citep[e.g.,][]{candes2000curvelets,candes2004new,le2005sparse,petersen2018optimal,imaizumi2019deep}.
Our function class subsumes the usual (homogeneous) anisotropic space for which adaptive procedures exist with optimal convergence rate guarantees, including the dyadic classification and regression trees (CART) of \citet{donoho1997cart}. We refer to \citet{barron1999risk}, \citet{neumann1997wavelet}, \citet{hoffman2002random}, \citet{lepski2015adaptive}, and references therein for a more complete list.
There are also adaptive Bayesian procedures for anisotropic function estimation with desired asymptotic properties \citep[e.g.,][]{bhattacharya2014anisotropic,shen2015adaptive}.
There appear to be {\em no}  theoretical properties for adaptation in the more general case of piecewise heterogeneous anisotropic smoothness.
Indeed, existing theoretical studies for discontinuous piecewise smooth classes impose the isotropy assumption \citep[e.g.,][]{candes2000curvelets,candes2004new,le2005sparse,petersen2018optimal,imaizumi2019deep} and
the convergence rates in spatially adaptive estimation depend on global smoothness parameters \citep[e.g.,][]{pintore2006spatially,liu2010data,wang2013smoothing,tibshirani2014adaptive}.
In this respect, our study appears to be the first theoretical investigation of piecewise anisotropic function classes.

The majority of  frequentist/Bayesian methods for anisotropic function estimation rely on multiple scaling (bandwidth) parameters, one for each direction.
As noted by \citet{bhattacharya2014anisotropic}, selecting optimal scaling parameters in a frequentist way can be computationally difficult, as adaptation in anisotropic spaces presents several challenges \citep{lepski1999adaptive}.
The Bayesian paradigm provides an effective remedy  by assigning priors over these unknown  parameters.
One such example is the generalized Gaussian process priors or spline basis representations
\citep{bhattacharya2014anisotropic,shen2015adaptive}. 
Although these priors enjoy elegant theoretical guarantees in typical anisotropic spaces, whether they can adapt to piecewise heterogeneous anisotropic spaces without substantial modification remains unclear.
Contrariwise, Bayesian forests are expected to work in these more complex scenarios  without any additional scaling parameters. The approximability is controlled merely by the depth of a tree and the orientation of its branches, where no prior modifications should be required to  achieve optimal performance.
Moreover, computation with Gaussian processes can be quite costly \citep{banerjee2013efficient,liu2020gaussian}, while Bayesian forests are more scalable and faster than their competitors.

In the context of regression or classification, Bayesian forests often rely on observed covariate values for  splits in recursive partitioning \citep{chipman1998bayesian,denison1998bayesian,chipman2010bart}. This facilitates  theoretical investigation under the fixed regression design.
In the context of nonparametric Gaussian regression, \citet{rockova2020posterior} and \citet{rockova2019theory} investigated posterior contraction for BART based on this conventional manner of partitioning, whereas the dyadic CART \citep{donoho1997cart} splits at dyadic midpoints of the domain and can achieve optimal performance as well \citep{castillo2021uncertainty}. We generalize the dyadic CART by introducing the notion of split-nets, which form a collection of candidate split-points that are not necessarily observed covariate values and/or dyadic midpoints. Our findings show that optimality can be achieved with split-nets that are sufficiently evenly distributed. By allowing the split-points to occur beyond observed values, we show that Bayesian forests enjoy the general recipe of the posterior contraction theory \citep{ghosal2000convergence,ghosal2007convergence}, which applies to other statistical setups such as density estimation or regression/classification with random design.

Asymptotic minimaxity is often used to evaluate the optimality of statistical procedures.
\citet{yang2015minimax} derived the minimax rates of sparse function estimation in high dimensions, but their results are restricted to the isotropic cases.
In fixed (low) dimensions, minimax rates over anisotropic function spaces have been extensively studied in the literature \citep{ibragimov2013statistical,nussbaum1985spline,birge1986estimating}. If the true function only depends on a subset of coordinates, the minimax rate is improved and determined by the smoothness parameters of active coordinates \citep{hoffman2002random}.
However, to the best of our knowledge, there are no available studies on minimax rates over piecewise anisotropic function spaces like ours.  While there exist results on piecewise isotropic classes \citep[e.g.,][]{imaizumi2019deep}, even the simpler fixed-dimensional setup without sparsity has {\em not} been studied for piecewise anisotropic classes.
Focusing on Gaussian nonparametric regression, we derive the minimax lower bound for our piecewise heterogeneous anisotropic spaces under the high-dimensional scenario. 
This result verifies the finding that our obtained contraction rates for Bayesian forests are indeed minimax-optimal up to a logarithmic factor.

We summarize the contribution of this study as follows.
\begin{itemize}
	\item {\bf Approximation theory}: The true function should be approximable by tree-based learners to establish the optimal rate of posterior contraction. Approximation theory for piecewise heterogeneous anisotropic classes is much more intricate when there are discontinuities and heterogeneity. We establish such approximation theory here under suitable regularity conditions (with smoothness up to $1$ owing to the limitation of piecewise constant learners).
	\item {\bf Posterior contraction}: For piecewise heterogeneous anisotropic functions, posterior contraction of Bayesian forests is established under the high-dimensional setup with a Dirichlet sparse prior. The derived rates consist of the risk of variable selection uncertainty and the risk of function estimation, similar to isotropic cases \citep{yang2015minimax,rockova2020posterior}.
	\item {\bf Minimax optimality}: Minimax rates in high-dimensional spaces  have  been unavailable even for simple anisotropic classes.
	For Gaussian nonparametric regression with high-dimensional inputs, we formally derive the minimax lower bound over piecewise heterogeneous anisotropic spaces.  This certifies that our obtained contraction rate for Bayesian forests is optimal up to a logarithmic factor.
	\item  {\bf Applications beyond regression}:
	Unlike the asymptotic studies of the traditional tree priors \citep{rockova2020posterior,rockova2019theory}, our findings show that splits for recursive partitioning do not necessarily have to be at observed covariate values.
	This implies that our technique of proofs extends beyond fixed-design regression to other estimation problems such as density estimation or regression/classification with random design.
\end{itemize}

\subsection{Preview and Outline of the Paper}
The main results of this study begin to appear in Section~\ref{sec:app} after excessive preliminary steps. Before going into the preparatory phase, here we provide a preview of our main results. Let us focus on a fixed design regression setup,
\begin{align}
	Y_i=f_0(x_i)+\varepsilon_i, \quad\varepsilon_i\sim \text{N}(0,\sigma_0^2), \quad i=1,\dots,n,
	\label{eqn:fixreg}
\end{align}
with a response $Y_i\in\mathbb R$ and a covariate $x_i\in[0,1]^p$, where $f_0:[0,1]^p\rightarrow \mathbb R$ and $\sigma_0^2<\infty$. Assume that $f_0$ depends only on $d$ variables among $p$ coordinates.
Assume further that $f_0$ is a piecewise heterogeneous anisotropic function with a global smoothness harmonic mean $\bar\alpha\in(0,1]$ (see Definitions~\ref{def:hol}--\ref{def:sparse} for a more precise definition).
Assigning the BART prior on $f_0$, the posterior contraction rate is obtained as
$\sqrt{(d\log p)/n}+(\log n)^c n^{-\bar\alpha/(2\bar\alpha+d)}$ for some $c>0$ (Theorem~\ref{thm:nonreg}). This rate is minimax-optimal up to a log factor (Theorem~\ref{thm:minimax}). The same contraction rates are also achieved in other statistical setups (Theorems~\ref{thm:nonregrd}--\ref{thm:binary}). For the additive true function, the rate has an additive form (Theorems~\ref{thm:addreg}).

The rest of this paper is organized as follows.
In Section~\ref{sec:background}, we describe the background of function spaces and Bayesian forests. In high-dimensional scenarios, the tree priors on functions are specified in Section~\ref{sec:prior}. In Section~\ref{sec:appgen}, we illuminate the approximation theory for our function spaces. In Section~\ref{sec:nonreg}, we study posterior contraction of Bayesian forests and their minimax optimality in nonparametric regression with a fixed design. The section also includes a numerical study that shows the outstanding performance of BART over other methods such as random forests and deep neural networks, which are believed to work well for discontinuous or complicated smooth functions.
Posterior contraction properties in other statistical models such as density estimation and binary classification are investigated in Section~\ref{sec:furapp}. An example of additive regression is also considered in Section~\ref{sec:furapp} to emphasize a theoretical advantage of Bayesian forests over single tree models.
Section~\ref{sec:disc} concludes.
All technical proofs are presented in Appendix.

\section{Preliminaries}\label{sec:background}

\subsection{Notation and Terminology}
Although the main focus of this study is BART for regression in \eqref{eqn:fixreg}, we work with a general statistical experiment $P_f$ indexed by a measurable function $f:[0,1]^p\rightarrow\mathbb R$ for some $p>0$, which will be modeled by Bayesian forests.
This allows us to incorporate other statistical setups, such as density estimation, into our theoretical framework.
Each statistical model we are dealing with will be specified for our examples in Sections~\ref{sec:nonreg}--\ref{sec:furapp}. 
We observe $n$ observations with the true function denoted by $f_0$ and assume that $p$ is possibly increasing with the sample size $n$.
The notation $\mathbb E_0$ denotes the expectation operator under the true model with $f_0$.

For sequences $a_n$ and $b_n$, we write $a_n\lesssim b_n$ (or $b_n\gtrsim a_n$ equivalently) if $a_n\le C b_n$ for some constant $C>0$, and $a_n\asymp b_n$ implies $a_n\lesssim b_n\lesssim a_n$.  We also write $a_n\ll b_n$ (or $b_n\gg a_n$ equivalently) if $a_n/b_n\rightarrow 0$ as $n\rightarrow\infty$.
For a subspace $E$ of the Euclidean space, $\mathcal C(E)$ denotes a class of continuous functions $f:E\rightarrow \mathbb R$.
For a given measure $\mu$ and a measurable function $f$, we denote by $\lVert f \rVert_{v,\mu}=(\int |f|^v d\mu)^{1/v}$ the $L_v(\mu)$-norm, $1\le v<\infty$.
We denote by $\mathcal L_2(\mu)$ the linear space of real valued functions equipped with inner product $\langle f,g \rangle_\mu=\int f g d\mu$ and norm $\lVert f \rVert_{2,\mu}=\langle f,f \rangle_\mu^{1/2}$. For the sake of brevity, with the Lebesgue measure on a unit hypercube, $\mathcal L_2$ denotes the $L_2$ space and $\lVert f \rVert_v$ denotes the $L_v$-norm.
In particular,  $\lVert f \rVert_\infty$ denotes the $L_\infty$-norm of a function $f$ defined by the essential supremum, i.e., $\lVert f \rVert_\infty = \inf\{C\ge 0 : |f(x)|\le C \text{ for almost every $x$}\}$.\footnote{We use the $L_\infty$-norm to measure the difference of discontinuous functions while ignoring possible disagreement at jump surfaces. The $L_\infty$-norm is reduced to the supremum-norm for continuous functions if the domain is not a null set.}
 The support of a measure $\mu$ is denoted by ${\rm supp}(\mu)$.
For a given vector $u$, the notations $\lVert u \rVert_v$ and $\lVert u \rVert_\infty$ represent the $\ell_v$-norms, $1\le v<\infty$, and the maximum-norm, respectively. 
For a semimetric space $(\mathcal F,\rho)$ endowed with a semimetric $\rho$, the expressions $D(\epsilon,\mathcal F, \rho)$ and $N(\epsilon,\mathcal F, \rho)$ are $\epsilon$-packing and $\epsilon$-covering numbers of $\mathcal F$, respectively.
For a subset $S\subseteq\{1,\dots, p\}$ and $x=(x_1,\dots, x_p)^\top\in\mathbb R^p$, let $x_S=(x_j,j\in S)\in\mathbb R^{|S|}$ be the indices chosen by $S$.

A $q$-dimensional hyper-rectangle $\Psi\subseteq[0,1]^q$ with any $q>0$ is simply called a {\em box}. Precisely, a box is defined as the Cartesian product of open, closed, or semi-closed intervals; therefore, a box can be open, closed, or neither (e.g., $[a_1,b_1)\times(a_2,b_2]$) depending on the context.
A partition $\mathfrak Y=\{\Psi_1,\dots,\Psi_J\}$ of $[0,1]^q$, consisting of $J$ disjoint boxes $\Psi_r\subseteq[0,1]^q$, $r=1,\dots,J$, is called a {\em box partition}.
For the Cartesian product of $q$ subsets of $\mathbb R$, i.e., $E\subseteq \mathbb R^q$, we denote the $j$th projection mapping of $E$ by $[E]_j=\{x_j\in\mathbb R: (x_1,\dots,x_q)^\top\in E\}$.  The length and interior of an interval $I\in\mathbb R$ is denoted by $\mathsf{len}(I)$ and $\mathsf{int}(I)$, respectively.

\subsection{Heterogeneous Anisotropic Function Spaces with Sparsity}
\label{sec:function}
In this subsection, we introduce our  function spaces with  heterogeneous smoothness and sparsity in high dimensions.
The first assumption is that the true regression function $f_0:[0,1]^p\rightarrow \mathbb R$ is $d$-sparse, i.e., it depends on a small subset of $d$ variables.
This means that there exist a function $h_0:[0,1]^d\rightarrow\mathbb R$ and a subset $S_0\subseteq\{1,\dots, p\}$ with $|S_0|=d$, such that
$f_0(x)=h_0(x_{S_0})$ for any $x\in[0,1]^p$. For example, suppose the true function is defined as $f_0(x_1,x_2)=\sin(x_1)$ on $[0,1]^2$ with $p=2$. This function can be completely expressed by the one-dimensional function $h_0(x_1)=\sin(x_1)$ on $[0,1]$, and hence is 1-sparse by definition.

For now, we focus on the function $h_0$ on the low-dimensional domain $[0,1]^d$. The complete characterization of $f_0$ will soon be discussed.
We assume that $[0,1]^d$ partitioned into many boxes and  $h_0$ is H\"older  continuous with possibly different smoothness in each box.
The smoothness inside each box is anisotropic, i.e.,  different for each coordinate.
Focusing on a single box, we first define an {\em anisotropic H\"older space} in the usual sense.

\begin{mydefinition}[Anisotropic H\"older space]
	\label{def:hol}
	For smoothness $\alpha=(\alpha_1,\dots,\alpha_d)^\top\in(0,1]^d$, a box $\Psi\subseteq[0,1]^d$, and a H\"older coefficient $\lambda<\infty$, we denote by $\mathcal H_\lambda^{\alpha,d}(\Psi)$ an anisotropic $\alpha$-H\"older space on $\Psi$, i.e.,
	\begin{align*}
		\mathcal H_\lambda^{\alpha,d}(\Psi) =\left\{ h:\Psi \rightarrow \mathbb R ; ~ |h(x)-h(y)|\le \lambda\sum_{j=1}^d |x_j-y_j|^{\alpha_j} ,~ x,y\in\Psi\right\}.
	\end{align*}
\end{mydefinition}

Note that the definition above imposes a restriction $\alpha\in(0,1]^d$.
Although one can generalize this definition to smoother classes \citep[e.g.][]{bhattacharya2014anisotropic},
we do not consider such extensions here, as step function estimators cannot be optimal in classes smoother than Lipschitz.

As discussed above, our targeted function class is not necessarily globally anisotropic over the entire domain $[0,1]^d$.
Instead, we assume that $h_0$ has different anisotropic smoothness on $R\ge1$ disjoint boxes of the domain with the same harmonic mean (the same harmonic mean is an important assumption for obtaining the minimax lower bound in Section~\ref{sec:minimax}). To be more precise, we define a set of $R$-tuples for smoothness parameters,
\begin{align*}
	\mathcal A_{\bar\alpha}^{R,d} = \left\{(\alpha_1,\dots,\alpha_R): \alpha_r=(\alpha_{r1},\dots,\alpha_{rd})^\top\in(0,1]^d, ~ \bar\alpha^{-1}=d^{-1}\sum_{j=1}^d \alpha_{rj}^{-1},~r=1,\dots, R \right\}.
\end{align*}
We assume that the anisotropic smoothness of $h_0$, the nonsparse proxy of $f_0$, is specified on an unknown underlying box partition $\mathfrak X_0=\{\Xi_1,\dots,\Xi_R\}$ of $[0,1]^d$ with $R\ge1$ boxes. If $R=1$, we write $\mathfrak X_0=\{[0,1]^d\}$ with $\Xi_1=[0,1]^d$. Note that each $\Xi_r$ can be open, closed, or neither.
The function space is formed by agglomerating anisotropic H\"older spaces for all boxes. 
We emphasize that the resulting function space is not necessarily continuous, which provides a lot more flexibility relative to the conventional H\"{o}lderian class.
Considering that smoothness parameters can vary across boxes and functions can be discontinuous at their boundaries, we call this new class a  {\em piecewise heterogeneous anisotropic H\"older space}. We define these functions formally below.

\begin{mydefinition}[Piecewise heterogeneous anisotropic H\"older space]
	\label{def:genhol}
	Consider a smoothness parameter $A_{\bar\alpha}=(\alpha_r)_{r=1}^R\in\mathcal A_{\bar\alpha}^{R,d}$ for some $\bar\alpha\in(0,1]$ and a box partition $\mathfrak Y=\{\Psi_1,\dots,\Psi_R\}$ of $[0,1]^d$ with boxes $\Psi_r\subseteq[0,1]^d$.\footnote{
		For any $q>1$, we write $\mathfrak Y=\{\Psi_r\}_r$ to denote an arbitrary box partition of $[0,1]^q$ with boxes $\Psi_r\subseteq[0,1]^q$, $r=1,2,\dots$, and write $\Psi\subseteq[0,1]^q$ to denote an arbitrary $q$-dimensional box.
	} We define a piecewise  heterogeneous anisotropic H\"older space as
	\begin{align*}
		\mathcal H_\lambda^{A_{\bar\alpha},d}(\mathfrak Y) =\left\{ h:[0,1]^d \rightarrow \mathbb R ; ~ h|_{\Psi_r} \in \mathcal H_\lambda^{\alpha_r,d}(\Psi_r) ,~r=1,\dots,R\right\}.
	\end{align*}
\end{mydefinition}

\begin{figure}[t!]
	\centering
	\resizebox{2.25in}{!}{
		\begin{tikzpicture}
			\draw  (0,0) rectangle (7,7);
			\draw  (0,2.5) -- (7,2.5);
			\draw  (2,2.5) -- (2,7);
			\draw (4.5,2.5) -- (4.5,0);
			\draw (2,5) -- (7,5);
			\node at (1,5) {$\Xi_1$};
			\node at (4.5,6) {$\Xi_2$};
			\node at (4.5,3.75) {$\Xi_3$};
			\node at (2.25,1.25) {$\Xi_4$};
			\node at (5.75,1.25) {$\Xi_5$};
			\node [rotate=270,blue] at (0.25,5) {\footnotesize$\alpha_{12}$};
			\draw [->,blue] (0.25,5.35) to (0.25,5.75);
			\draw [->,blue] (0.25,4.65) to (0.25,4.25);
			\node [rotate=270,blue] at (4.75,1.25) {\footnotesize$\alpha_{52}$};
			\draw [->,blue] (4.75,1.6) to (4.75,2);
			\draw [->,blue] (4.75,0.9) to (4.75,0.5);
			\node [rotate=270,blue] at (0.25,1.25) {\footnotesize$\alpha_{42}$};
			\draw [->,blue] (0.25,1.6) to (0.25,2);
			\draw [->,blue] (0.25,0.9) to (0.25,0.5);
			\node [rotate=270,blue] at (2.25,6) {\footnotesize$\alpha_{22}$};
			\draw [->,blue] (2.25,6.35) to (2.25,6.75);
			\draw [->,blue] (2.25,5.65) to (2.25,5.25);
			\node [rotate=270,blue] at (2.25,3.75) {\footnotesize$\alpha_{32}$};
			\draw [->,blue] (2.25,4.1) to (2.25,4.5);
			\draw [->,blue] (2.25,3.4) to (2.25,3);
			\node [blue] at (4.5,2.75) {\footnotesize$\alpha_{31}$};
			\draw [->,blue] (4.85,2.75) to (5.25,2.75);
			\draw [->,blue] (4.15,2.75) to (3.75,2.75);
			\node [blue] at (4.5,5.25) {\footnotesize$\alpha_{21}$};
			\draw [->,blue] (4.85,5.25) to (5.25,5.25);
			\draw [->,blue] (4.15,5.25) to (3.75,5.25);
			\node [blue] at (2.25,0.25) {\footnotesize$\alpha_{41}$};
			\draw [->,blue] (2.6,0.25) to (3,0.25);
			\draw [->,blue] (1.9,0.25) to (1.5,0.25);
			\node [blue] at (5.75,0.25) {\footnotesize$\alpha_{51}$};
			\draw [->,blue] (6.1,0.25) to (6.5,0.25);
			\draw [->,blue] (5.4,0.25) to (5,0.25);
			\node [blue] at (1,2.75) {\footnotesize$\alpha_{11}$};
			\draw [->,blue] (1.35,2.75) to (1.75,2.75);
			\draw [->,blue] (0.65,2.75) to (0.25,2.75);
		\end{tikzpicture}
	}
	\caption{A graphical illustration of a piecewise heterogeneous anisotropic H\"older space with five boxes. Each piece has its own smoothness parameter, but the harmonic mean is assumed to be the same.}
	\label{fig:anisohol}
\end{figure}

A graphical illustration of the piecewise heterogeneous anisotropic H\"older spaces is given in Figure~\ref{fig:anisohol}.
Clearly, Definition~\ref{def:genhol} subsumes the anisotropic H\"older space in Definition~\ref{def:hol} with $R=1$.
According to  Definition~\ref{def:genhol}, any $h\in\mathcal H_\lambda^{A_{\bar\alpha},d}(\mathfrak Y)$ is anisotropic on each $\Psi_r$ with a smoothness parameter $\alpha_r\in(0,1]^d$ and the same harmonic mean $\bar\alpha$ for all $\Psi_r$. We again emphasize that discontinuities are allowed at the boundaries of boxes $\Psi_r$, $r=1,\dots,R$.

Definition~\ref{def:genhol} does not impose a specific structure on the partition $\mathfrak Y$ other than a box partition.
However, we will later see that, depending on the approximation metric, our approximation theory will require $\mathfrak X_0$ to be a tree-based recursive structure defined in the next section (see Figure~\ref{fig:partition} below). Nonetheless, as every box partition can be extended to the required form by adding more splits, this discrepancy can be addressed, but it may harm our posterior contraction rate. We refer the reader to Section~\ref{sec:dennet} for more discussion.

\begin{myremark}
	We compare Definition~\ref{def:genhol} with piecewise smooth function spaces widely investigated in the literature. Approximation rates for piecewise smooth functions with smooth jump curves/surfaces have been extensively studied in two dimensions \citep[e.g.,][]{candes2000curvelets,candes2004new,guo2007optimally} as well as in higher dimensions \citep{chandrasekaran2008representation,petersen2018optimal,imaizumi2019deep}. 
	All these studies deal with smooth functions with smooth jump curves/surfaces under the isotropy assumption. 
	contrariwise, our definition deals with different anisotropic smoothness parameters for the boxes in a box partition, and hence seems to offer some flexibility. Our jump surfaces, however, are restricted to hyper-planes parallel to the coordinates. 
\end{myremark}

\begin{myremark}
We believe that our function class is not a subset of a popular one, but is originally defined in our work. For example, anisotropic and mixed smooth Besov spaces are highly flexible classes that render discontinuity and spatially varying smoothness \citep{suzuki2019adaptivity,suzuki2021deep}, but they do not account for our piecewise heterogeneous anisotropic smoothness in Definition~\ref{def:genhol}. In our construction, the axis-aligned box partition appears to be an important assumption in obtaining the optimal posterior contraction rate using our theory. Later we will see that our contraction rate depends on $R$, which is translated as the number of binary splits required to approximate the true $\mathfrak X_0^\ast$ (see Section~\ref{sec:dennet}). If the partition is not axis-aligned, infinitely many splits are needed, which will deteriorate our rate. Whether this is a fundamental limitation of BART is still unclear.
\end{myremark}

Note that Definition~\ref{def:genhol} can be used for the mapping $h_0$ from the lower dimensional domain $[0,1]^d$ while the true function $f_0$ maps the entire $[0,1]^p$ to $\mathbb R$.
We now characterize a {\em sparse} elaboration of Definition~\ref{def:genhol} for the mapping $f_0:[0,1]^p\rightarrow \mathbb R$.
For any $S\subseteq\{1,\dots,p\}$, we denote with
$W_S^p:\mathcal C(\mathbb R^{|S|})\rightarrow \mathcal C(\mathbb R^p)$ the map that transmits 
$h\in\mathcal C(\mathbb R^{|S|})$ onto 
$W_S^p h : x\mapsto h(x_S)$.
Similar to \citet{yang2015minimax} for the isotropic cases, we now formalize $d$-sparse function spaces as follows.

\begin{mydefinition}[Sparse function space]
	For the space $\mathcal H_\lambda^{A_{\bar\alpha},d}(\mathfrak Y)$ in Definition~\ref{def:genhol}, we define a $d$-sparse piecewise heterogeneous anisotropic H\"older space as 
	\begin{align*}
		\Gamma^{A_{\bar\alpha},d,p}_{\lambda}(\mathfrak Y)=\bigcup_{S\subseteq\{1,\dots,p\}: |S|=d} W_S^p\big(\mathcal H_\lambda^{A_{\bar\alpha},d}(\mathfrak Y) \big).
	\end{align*}
	\label{def:sparse}
\end{mydefinition}

That is, $\Gamma^{A_{\bar\alpha},d,p}_{\lambda}(\mathfrak Y)$ is read as the collection of $p$-dimensional $d$-sparse functions over $\mathfrak Y$ with piecewise anisotropic $\bar\alpha$ smoothness and a Lipschitz constant $\lambda$.
For an unknown smoothness parameter $A_{\bar\alpha}=(\alpha_r)_{r=1}^R\in\mathcal A_{\bar\alpha}^{R,d}$ (with possibly decreasing $\bar\alpha$) and
model components $R$, $d$, $p$, and $\lambda$ (which are possibly increasing with $n$), the true function $f_0$ is assumed to belong to the class $\Gamma^{A_{\bar\alpha},d,p}_{\lambda}(\mathfrak X_0)$ which allows for discontinuities, or to its continuous variant $\Gamma^{A_{\bar\alpha},d,p}_{\lambda}(\mathfrak X_0)\cap \mathcal C([0,1]^p)$.
This means that there exists a function $h_0:[0,1]^d\rightarrow\mathbb R$ and a subset $S_0\subseteq\{1,\dots,p\}$ with $|S_0|=d$ such that $f_0=W_{S_0}^p h_0$.
The continuous variant $\Gamma^{A_{\bar\alpha},d,p}_{\lambda}(\mathfrak X_0)\cap \mathcal C([0,1]^p)$ achieves approximability under more relaxed assumptions (see Theorem~\ref{thm:approx2} in Section~\ref{sec:app}). The two spaces are identical if $R=1$.

Note that the true underlying $\mathfrak X_0$ is the box partition of the $d$-dimensional cube $[0,1]^d$.
Considering the domain $[0,1]^p$ of $f_0$, it will be convenient to extend $\mathfrak X_0$ to the corresponding box partition of the $p$-dimensional cube $[0,1]^p$. To this end, we extend each $\Xi_r$ to the $p$-dimensional box $\Xi_r^\ast=\{x\in[0,1]^p: x_{S_0}\in \Xi_r ,x_{S_0^c}\in[0,1]^{p-d} \} \subseteq[0,1]^p$ using the true sparsity index $S_0$; that is, $\Xi_r$ is the projection of $\Xi_r^\ast$ onto the coordinates in $S_0$.
The boxes $\Xi_r^\ast$ then constitute the box partition $\mathfrak X_0^\ast=\{\Xi_1^\ast,\dots,\Xi_R^\ast\}$ of $[0,1]^p$.\footnote{
	The notations $\mathfrak X_0=\{\Xi_1,\dots,\Xi_R\}$ and $\mathfrak X_0^\ast=\{\Xi_1^\ast,\dots,\Xi_R^\ast\}$ are used only to denote the true underlying box partition for the anisotropic smoothness of $h_0$ and its extension to the $p$-dimensional space for $f_0$, respectively.
}
We emphasize that $\mathfrak X_0^\ast$ is determined by the unknown sparsity index $S_0$ of the true function $f_0$. Observe also that our definition gives rise to $\mathfrak X_0^\ast=\{[0,1]^p\}$ with  $\Xi_1^\ast=[0,1]^p$ if $R=1$.

Apart from the notion of sparsity for functions, we also introduce sparsity of box partitions as follows.

\begin{figure}[t!]
	\centering
	\begin{subfigure}[b]{2in}
		\resizebox{2in}{!}{
			\begin{tikzpicture}
				\draw [-stealth, dashed] (0,0,5) -- (0,0,6) node [left] {\Large$1$};
				\draw [-stealth, dashed] (5,0,0) -- (6,0,0) node [above right] {\Large$2$};
				\draw [-stealth, dashed] (0,5,0) -- (0,6,0) node [above left] {\Large$3$};
				\draw (0,0,0) -- (5,0,0) -- (5,5,0) -- (0,5,0)  -- (0,0,0) ;
				\draw (5,5,5)-- (0,5,5)-- (0,0,5)-- (5,0,5)-- (5,5,5);
				\draw (0,0,0)-- (0,0,5);
				\draw (0,5,0)-- (0,5,5);
				\draw (5,0,0)-- (5,0,5);
				\draw (5,5,0)-- (5,5,5);
				\draw [dashed,blue] (2.5,0,0)-- (2.5,0,5)-- (2.5,5,5)-- (2.5,5,0)-- (2.5,0,0);
				\draw [dashed,blue] (2.5,0,2)-- (2.5,5,2)-- (5,5,2)-- (5,0,2)-- (2.5,0,2);
				\draw [dashed,blue] (2.5,0,3)-- (2.5,5,3)-- (0,5,3)-- (0,0,3)-- (2.5,0,3);
			\end{tikzpicture}
		}
		\caption{$\{1,2\}$-chopped partition}
		\label{fig:sppar1}
	\end{subfigure}
	\begin{subfigure}[b]{2in}
		\resizebox{2in}{!}{
			\begin{tikzpicture}
				\draw [-stealth, dashed] (0,0,5) -- (0,0,6) node [left] {\Large$1$};
				\draw [-stealth, dashed] (5,0,0) -- (6,0,0) node [above right] {\Large$2$};
				\draw [-stealth, dashed] (0,5,0) -- (0,6,0) node [above left] {\Large$3$};
				\draw (0,0,0) -- (5,0,0) -- (5,5,0) -- (0,5,0)  -- (0,0,0) ;
				\draw (5,5,5)-- (0,5,5)-- (0,0,5)-- (5,0,5)-- (5,5,5);
				\draw (0,0,0)-- (0,0,5);
				\draw (0,5,0)-- (0,5,5);
				\draw (5,0,0)-- (5,0,5);
				\draw (5,5,0)-- (5,5,5);
				\draw [dashed,blue] (2.5,0,0)-- (2.5,0,5)-- (2.5,5,5)-- (2.5,5,0)-- (2.5,0,0);
				\draw [dashed,blue] (1.2,0,0)-- (1.2,0,5)-- (1.2,5,5)-- (1.2,5,0)-- (1.2,0,0);
				\draw [dashed,blue] (4,0,0)-- (4,0,5)-- (4,5,5)-- (4,5,0)-- (4,0,0);
			\end{tikzpicture}
		}
		\caption{$\{2\}$-chopped partition}
		\label{fig:sppar2}
	\end{subfigure}
	\caption{Examples of sparse partitions in three dimensions.}
	\label{fig:sppar}
\end{figure}

\begin{mydefinition}[Sparse partition]
	Consider a box partition $\mathfrak Y=\{\Psi_1,\dots,\Psi_J\}$ of $[0,1]^p$ with boxes $\Psi_r\subseteq [0,1]^p$, $r=1,\dots,J$.
	For a subset $S\subseteq\{1,\dots,p\}$, the partition $\mathfrak Y$ is called $S$-chopped if $\max_{j \in S}\mathsf{len}([\Psi_r]_j)<1$ and $\min_{j\notin S}\mathsf{len}([\Psi_r]_j)=1$ for every $r=1,\dots,J$.
	\label{def:sppartition}
\end{mydefinition}

A graphical illustration of sparse partitions is provided in Figure~\ref{fig:sppar}.
According to Definition~\ref{def:sppartition}, the extended box partition $\mathfrak X_0^\ast$ is $S$-chopped for some $S\subseteq S_0$. Observe that $\mathfrak X_0^\ast$ is not always $S_0$-chopped, since $\mathfrak X_0$ may not have been cleaved in some coordinates. 
For example, if $f_0(x_1,x_2,x_3)=h_0(x_1,x_3)=\sin(x_1)\cos(x_3)\mathbbm 1(0\le x_1\le 0.5)\mathbbm 1(0\le x_2\le 1)$ with $p=3$ and $d=2$, then $S_0=\{1,3\}$, but $\mathfrak X_0^\ast=\{[0,0.5]\times[0,1]^2,(0.5,1]\times[0,1]^2\}$ is $\{1\}$-chopped.
In particular, $\mathfrak X_0^\ast$ is $\varnothing$-chopped if $R=1$ irrespective of what $S_0$ is. 
It is then clear that sparsity of $\mathfrak X_0^\ast$ is not the same as sparsity of $f_0$.
In what follows, we write $S_0^\ast\subseteq S_0$ to denote sparsity of $\mathfrak X_0^\ast$;  that is, $\mathfrak X_0^\ast$ is $S_0^\ast$-chopped.

\begin{myremark}
	Throughout the study, the model parameters $\bar\alpha$, $R$, $d$, $p$, and $\lambda$ are treated as positive sequences of $n$, which can vary at appropriate rates so that our target posterior contraction rate in \eqref{eqn:rate} changes. Accordingly, the model objects related to these sequences, e.g., $\mathfrak X_0$, $\mathfrak X_0^\ast$, and $A_{\alpha}\in\mathcal A_\alpha^{R,d}$, can also vary with $n$.
	The only exception is the minimax study in Section~\ref{sec:minimax}, where a fixed $d$ provides a correct interpretation of the obtained minimax lower bound (see the lower bound in Theorem~\ref{thm:minimax}). With a slight abuse of notation, we usually suppress the dependency on $n$ for the sake of notational simplicity.
\end{myremark}

\subsection{Tree-Based Partitions}
\label{sec:partition}
In this work, for estimators of the true function $f_0$,
we focus on piecewise constant learners, i.e., step functions that are constant on each piece of a box partition of $[0,1]^p$.
A precise description of piecewise constant learners requires an underlying partitioning rule that produces a partition for these step functions. In tree-structured models, the idea is based on recursively applying binary splitting rules to split the domain $[0,1]^p$.
Here we shed light on this mechanism to construct tree-based partitions, while deferring a complete description of the induced step functions to Section~\ref{sec:bayesforest}.

For a given box $\Psi\subseteq[0,1]^p$, 
choose a {\em splitting coordinate} $j\in\{1,\dots, p\}$ and a {\em split-point} $\tau_j\in{\mathsf{int}}([\Psi]_j)$. The pair $(j,\tau_j)$ then dichotomizes $\Psi$ along the $j$th coordinate into two boxes: $\{ x\in\Psi: x_j \le \tau_j \}$ and $\{ x\in\Psi:x_j > \tau_j \}$, where $x_j$ is $j$th entry of $x$.
Starting from the root node $[0,1]^p$, the procedure is iterated $K-1$ times in a top-down manner by picking one box for a split each time. This generates $K$ disjoint boxes $\Psi_1,\dots,\Psi_K$, called {\em terminal nodes}, which constitute a tree-shaped partition of $[0,1]^p$, called a {\em tree partition}.
We call this iterative procedure the {\em binary tree partitioning}. We will further refer to the resulting tree partitions as {\em flexible tree partitions} to emphasize that splits can occur everywhere in the domain $[0,1]^p$ (not necessarily at dyadic midpoints or observed covariate values). 
According to Definition~\ref{def:sppartition}, we say that a flexible tree partition is $S$-chopped if splitting coordinates $j$ are restricted to a subset $S\subseteq\{1,\dots,p\}$.
Note that while flexible tree partitions are always box partitions, the reverse is not generally true; see Figure~\ref{fig:partition}.

\begin{figure}[t!]
	\centering
	\begin{subfigure}[b]{1.8in}
		\resizebox{1.8in}{!}{\quad
			\begin{tikzpicture}
				\draw  (-33.5,0) rectangle (-26.5,7);
				\draw (-31.5,7) -- (-31.5,5) -- (-28.5,5) -- (-28.5,0);
				\draw (-33.5,2) -- (-28.5,2);
				\draw (-31.5,2) -- (-31.5,5);
				\draw (-26.5,5) -- (-28.5,5);
			\end{tikzpicture}\quad
		}
		\caption{A non-tree box partition}
		\label{fig:partition1}
	\end{subfigure}
	\begin{subfigure}[b]{1.8in}
		\resizebox{1.8in}{!}{\quad
			\begin{tikzpicture}
				\draw  (-33.5,0) rectangle (-26.5,7);
				\draw (-30.5,7) -- (-30.5,5) -- (-28,5) -- (-28,0);
				\draw (-33.5,2.5) -- (-30.5,2.5);
				\draw (-30.5,0) -- (-30.5,5);
				\draw (-26.5,5) -- (-28,5);
			\end{tikzpicture}\quad
		}
		\caption{A tree partition.}
		\label{fig:partition2}
	\end{subfigure}
	\caption{Examples of non-tree box partitions and tree partitions
	}
	\label{fig:partition}
\end{figure}

Although the binary tree partitioning allows splits to occur anywhere in the domain, Bayesian tree models usually take advantage of priors that choose split-points from a predetermined discrete set.
For example, in regression with continuous covariates, observed covariate values are typically used for split-points \citep{chipman1998bayesian,denison1998bayesian,chipman2010bart}.
Following this manner,
\citet{rockova2020posterior} and \citet{rockova2019theory} investigated posterior contraction of BART in Gaussian nonparametric regression with fixed covariates. Here, we relax this restriction while keeping split-points chosen from a discrete set.
To this end, we define a discrete collection of locations where splits can occur, which we call a split-net.

\begin{mydefinition}[split-net]
	For an integer sequence $b_n$, a split-net $\mathcal Z=\{z_i\in[0,1]^p,~i=1,\dots,b_n\}$ is a set of $b_n$ points $z_i=(z_{i1},\dots,z_{ip})^\top\in[0,1]^p$ at which possible splits occur along coordinates. 
	\label{def:split-net}
\end{mydefinition}

\begin{figure}[t!]
	\centering
	\begin{subfigure}[b]{2in}
		\includegraphics[width=\linewidth]{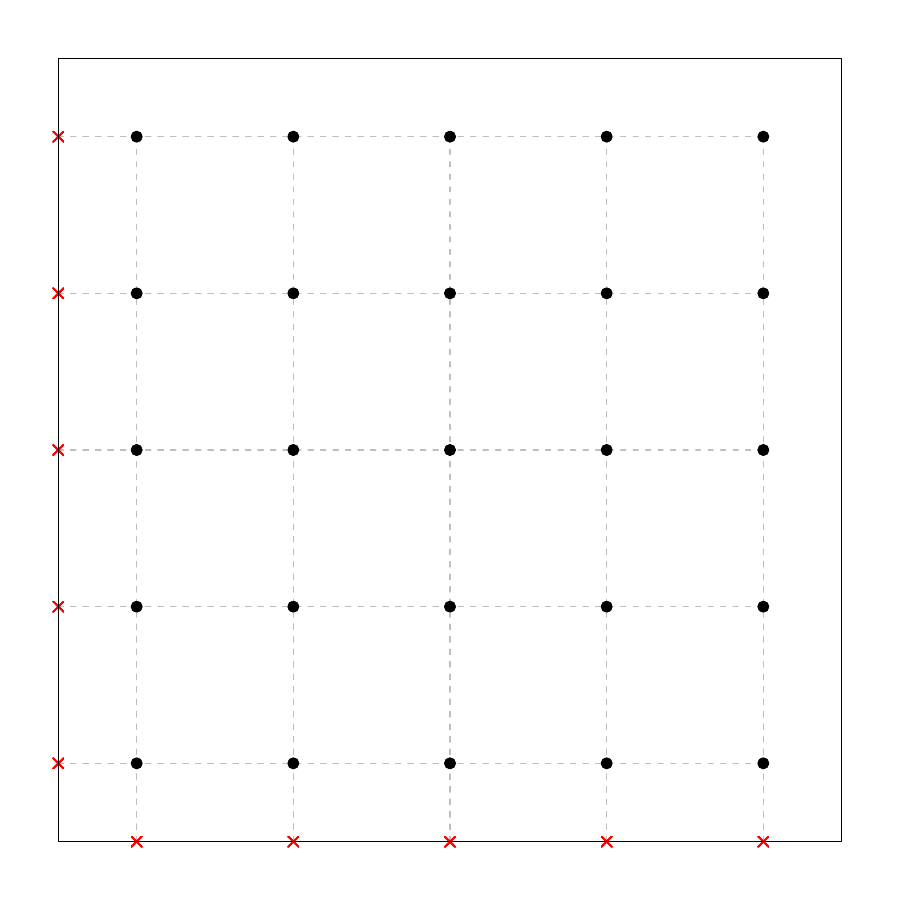}\vspace*{-0.1in}		
		\caption{Regular grid system}
		\label{fig:spn1}
	\end{subfigure}
	\begin{subfigure}[b]{2in}
		\includegraphics[width=\linewidth]{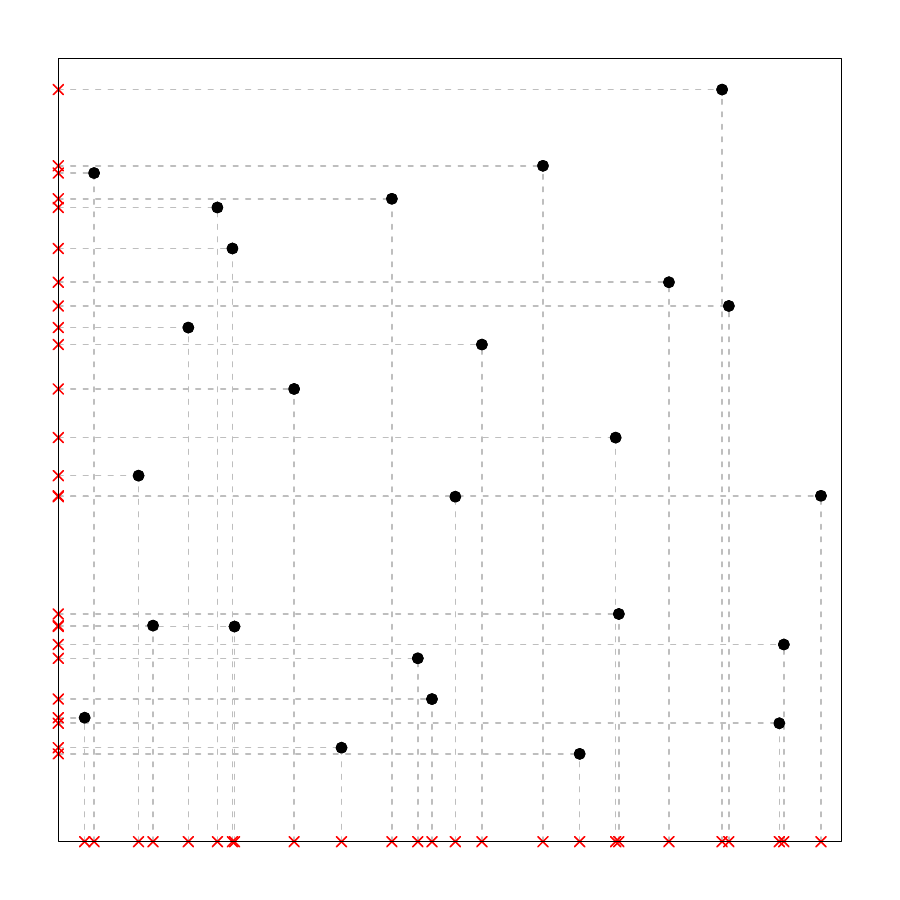}\vspace*{-0.1in}
		\caption{Split-net without duplication}
		\label{fig:spn2}
	\end{subfigure}
	\caption{Examples of the split-net with $b_n=25$ in two dimensions.
		For the regular grid in (a), one can easily see that $b_j(\mathcal Z)=5$, $j=1,2$; hence, initial splits eliminate the possibility of other splits.
		The split-candidates of the split-net in (b) are unique in every coordinate, so $b_j(\mathcal Z)=b_n$, $j=1,2$.
	}
	\label{fig:spn}
\end{figure}

For a given split-net $\mathcal Z$, we call each point $z_i=(z_{i1},\dots,z_{ip})^\top$ a {\em split-candidate}. For a given splitting coordinate $j$ and a split-net $\mathcal Z$, a split-point will be chosen from $[\mathcal Z]_j\cap\mathsf{int}([\Psi]_j)$ to dichotomize a box $\Psi$.
Note that $[\mathcal Z]_j=\{z_{ij}\in[0,1],i=1,\dots,b_n\}$ may have fewer elements than $\mathcal Z$ owing to duplication.
We denote by $b_j(\mathcal Z)$ the cardinality of $[\mathcal Z]_j$, i.e., the number of unique values in the $b_n$-tuple $(z_{1j},\dots,z_{b_n j})$. We then obtain $\max_{1\le j\le p} b_j(\mathcal Z)\le b_n$ by definition. For example, consider a regular (equidistant) grid system illustrated in Figure~\ref{fig:spn1}, wherein $b_j(\mathcal Z)=b_n^{1/p}<b_n$, $j=1,\dots,p$. This simplest split-net will be further discussed in Section~\ref{sec:grid}.
It is also possible to construct a split-net such that $b_j(\mathcal Z)=b_n$, $j=1,\dots, p$, as shown in Figure~\ref{fig:spn2}.
As noted above, another typical example of $\mathcal Z$ is the observed covariate values in fixed-design nonparametric regression with $b_n=n$ (supposing that all $x_i$ are different). This specific example will be discussed in Section~\ref{sec:fixeddesign}.
Our definition of split-nets yields additional flexibility in situations when no deterministic covariate values are available, such as density estimation or in the analysis of nonparametric regression with random covariates. 
A subset of the observed covariate values can also be used in a fixed-design regression setup.

In assigning a prior over tree partitions, we will assume that splits in the binary partitioning rule occur only at the points in $\mathcal Z$; that is, for every splitting box $\Psi\subseteq[0,1]^p$ with a splitting coordinate $j$, a split-point $\tau_j$ is chosen such that $\tau_j\in[\mathcal Z]_j\cap\mathsf{int}([\Psi]_j)$. As a split is restricted to the interior of a given interval, some split-candidates may have already been eliminated in the previous steps of the splitting procedure (see Figure~\ref{fig:spn1}).
Clearly, a tree partition constructed by $\mathcal Z$ is an instance of flexible tree partitions, but the reverse is not the case.
To distinguish between the two more clearly, we make the following definition.

\begin{mydefinition}[$\mathcal Z$-tree partition]
	For a given split-net $\mathcal Z$, a flexible tree partition $\mathcal T=\{\Omega_1,\dots,\Omega_K\}$ of $[0,1]^p$ with boxes $\Omega_k\subseteq[0,1]^p$, $k=1,\dots,K$, is called a $\mathcal Z$-tree partition if every split occurs at points $z_i\in\mathcal Z$.\footnote{
		The notation $\mathcal T=\{\Omega_k\}_k$ is used only for the $\mathcal Z$-tree partitions with a split-net $\mathcal Z$, with some suitable superscript and/or superscript if required. We denote flexible tree partitions by $\mathfrak Y=\{\Psi_k\}_k$ as general box partitions.
	}
\end{mydefinition}

In summary, we obtain the following relationship among the three types of partitions: $\{\text{$\mathcal Z$-tree partitions}\}\subseteq\{\text{Flexible tree partitions}\}\subseteq\{\text{Box partitions}\}$.
Similar to flexible tree partitions, $\mathcal Z$-tree partitions can be $S$-chopped for a subset $S\subseteq\{1,\dots,p\}$ irrespective of what $\mathcal Z$ is employed.
As we aim to do sparse estimation in high-dimensional setups, we are primarily interested in $S$-chopped $\mathcal Z$-tree partitions for some low-dimensional $S$.
In what follows, we denote by $\mathscr T_{S,K,\mathcal Z}$ the set of all $S$-chopped $\mathcal Z$-tree partitions with $K$ boxes.

\begin{myremark}
	The definition of a $\mathcal Z$-tree partition is introduced to restrict possible splits to a discrete set. This means that we assign a discrete prior on the tree topologies (see Section~\ref{sec:prior}). One may instead assign a prior on the topology of flexible tree partitions, in which case a split-net $\mathcal Z$ is not needed. 
	For regression problems, most of the recent BART procedures deploy a discrete set of split-candidates in their prior constructions using the observed covariate values. 
	We aim to generalize this conventional idea while incorporating it into our framework. 
	A discrete prior has an advantage in that it is invariant to a transformation of predictor variables \citep{chipman1998bayesian}.
We only consider placing a discrete tree prior using a given split-net $\mathcal Z$, and a continuous prior on flexible tree partitions is not considered.	
\end{myremark}

\subsection{Bayesian Trees and Forests}
\label{sec:bayesforest}

We now describe our piecewise constant learners using $\mathcal Z$-tree partitions.
While single tree learners have received some attention \citep{chipman1998bayesian, denison1998bayesian}, it is widely accepted that additive aggregations of small trees are much more effective for prediction \citep{chipman2010bart}. Noting that single trees are a special case of tree ensembles (forests), we will focus on forests throughout the rest of the paper.

We consider a fixed number $T$ of trees.  For a given split-net $\mathcal Z$ and for each $ t\leq T$,  we denote with $\mathcal T^t=\{\Omega_1^t,\dots,\Omega_{K^t}^t\}$ a $\mathcal Z$-tree partition of size $K^t$ and  with $\beta^t=(\beta_1^t,\dots,\beta_{K^t}^t)^\top\in \mathbb R^{K^t}$ the heights of the step function, called the {\em step-heights}.
An additive tree-based learner is then fully described by a tree ensemble $\mathcal E =\{\mathcal T^1,\dots, \mathcal T^T\}$ and terminal node parameters
$B=({\beta^{1\top}},\dots,{\beta^{T\top}})^\top\in \mathbb R^{\sum_{t=1}^T K^t}$  through
\begin{align}
	f_{\mathcal E, B}(x)=\sum_{t=1}^T \sum_{k=1}^{K^t} \beta_k^t \mathbbm 1(x\in\Omega_k^t).
	\label{eqn:treelearner}
\end{align}
That is, $f_{\mathcal E, B}$ is constant on the boxes constructed by overlapping $\mathcal Z$-tree partitions $\mathcal T^1,\dots,\mathcal T^T$.
\citet{chipman2010bart} recommends the choice $T=200$, which was seen to provide good empirical results.
For a given ensemble $\mathcal E$, we henceforth define $\mathcal F_{\mathcal E}=\{f_{\mathcal E, B}: B\in\mathbb R^{\sum_{t=1}^T K^t}\}$ the set of functions in \eqref{eqn:treelearner}.
If $\mathcal E$ consists of a single tree $\mathcal T$, we instead write $\mathcal F_{\mathcal T}$ to denote $\mathcal F_{\mathcal E}$.

Our objective is to characterize the posterior asymptotic properties of the tree learners in \eqref{eqn:treelearner} in estimating the true function $f_0$ belonging to $\Gamma^{A_{\bar\alpha},d,p}_{\lambda}(\mathfrak X_0)$ or $\Gamma^{A_{\bar\alpha},d,p}_{\lambda}(\mathfrak X_0)\cap \mathcal C([0,1]^p)$.
This goal requires two nice attributes of the procedure.
First, appropriate prior distributions should be assigned to the tree learners $f_{\mathcal E, B}$ in \eqref{eqn:treelearner} so that the induced posterior can achieve the desired asymptotic properties.
Second, there should exist a piecewise tree learner approximating $f_0$ with a suitable approximation error matched to our target rate. In the following two sections, we elucidate these in detail.

\section{Tree and Forest Priors in High Dimensions}
\label{sec:prior}

\subsection{Priors over Tree Topologies with Sparsity}
\label{sec:treeprior}

Conventional tree priors \citep{chipman1998bayesian,denison1998bayesian} are not designed for high-dimensional data with a sparse underlying structure. Prior modifications are thus required for trees to meet demands of high-dimensional applications \citep{linero2018bayesian-jasa,linero2018bayesian-jrssb,rockova2020posterior}. 
\citet{rockova2020posterior} adopted a spike-and-slab prior for BART to achieve adaptability to unknown sparsity levels, but the computation of the posterior distribution is much more challenging than the original BART algorithm owing to the nature of a point mass prior. \citet{linero2018bayesian-jasa} and \citet{linero2018bayesian-jrssb} considered a sparse Dirichlet prior on splitting coordinates for a computationally feasible algorithm, while achieving the theoretical optimality in the high-dimensional scenario. We deploy the sparse Dirichlet prior developed by \citet{linero2018bayesian-jasa} for ease of computation for the posterior distribution.

Unlike the original tree priors, the BART model with the sparse Dirichlet prior chooses a splitting coordinate $j$ is from a proportion vector $\eta=(\eta_1,\dots,\eta_p)^\top$ belonging to the $p$-dimensional simplex $\mathbb S^p=\{(x_1,\dots,x_p)^\top\in\mathbb R^p: \sum_{j=1}^p x_j=1, x_j\ge 0,j=1,\dots,p\}$. A proportion vector $\eta$ has a Dirichlet prior with $\zeta>0$ and $\xi>1$,
\begin{align}
	\eta=(\eta_1,\dots,\eta_p)^\top\sim \text{Dir}(\zeta/p^{\xi},\dots,\zeta/p^{\xi}).
	\label{eqn:dir}
\end{align}
The requirement $\xi>1$ is needed for technical reasons. The prior imposes a sparsity into splitting variables (we refer the reader to Figure~2 of \citet{linero2018bayesian-jasa}).
Given a proportion vector $\eta$, the BART prior is assigned, as in \citet{chipman2010bart}, with a minor modification. Assuming an independent product prior for $\mathcal E$, i.e., $\Pi(\mathcal E)=\prod_{t=1}^T\Pi(\mathcal T^t)$, a Bayesian CART prior \citep{chipman1998bayesian} is assigned to each $\mathcal T^t$. The procedure begins with the root node $[0,1]^p$ of depth $\ell=0$, where the depth of a node means the number of nodes along the path from the root node down to that node. For each $\ell=0,1,2,\dots$, each node at depth $\ell$ is split with prior probability $\nu^{\ell+1}$ for $\nu\in(0,1/2)$. If a node corresponding to a box $\Omega$ is split, a splitting coordinate $j$ is drawn from the proportion vector $\eta$ and a split-point $\tau_j$ will be chosen randomly from $[\mathcal Z]_j\cap\mathsf{int}([\Omega]_j)$ for a given $\mathcal Z$. The procedure repeats until all nodes are terminal. 

The original CART prior proposed by \citet{chipman1998bayesian} uses a splitting probability that decays polynomially. \citet{rockova2019theory} showed that this decay may not be fast enough, and suggested using an exponentially decaying probability as ours.
This modification gives rise to the desirable exponential tail property of tree sizes. \citet{linero2018bayesian-jrssb} handled this issue by assigning a prior on the number $T$ of trees. As we want to fix $T$ as in the practical implementation of BART, we use the exponentially decaying prior probability for splits.

\subsection{Prior on Step-Heights}
\label{sec:betaprior}
To complete the prior on the sparse function space,
what remains to be specified is the prior on step-heights $B$ in \eqref{eqn:treelearner}. Given $K^1,\dots,K^T$ induced by $\mathcal E$, \citet{chipman2010bart} suggests using a Gaussian prior on $B$ (after shifting and rescaling the responses):
\begin{align*}
	d\Pi(B| K^1,\dots, K^T) = \prod_{t=1}^T \prod_{k=1}^{K^t} \phi(\beta_k^t ; 0,c_\beta/T),
\end{align*}
where $c_\beta>0$ is a constant and $\phi(\, \cdot\, ;\mu,\tau^2)$ is the Gaussian density with mean $\mu$ and variance $\tau^2$.
The variance $c_\beta/T$ shrinks step-heights toward zero, limiting the effect of individual components by keeping them small enough for large $T$.
This choice is preferred in view of the practical performance, but any zero-mean multivariate Gaussian prior on $B$ gives rise to the same optimal properties as soon as the eigenvalues of the covariance matrix are bounded below and above.
Throughout the paper, we place a Gaussian prior on the step-heights $B$ in most cases.
From the computational point of view, this choice is certainly appealing in Gaussian nonparametric regression owing to its semi-conjugacy.
For theoretical purposes, a prior with exponentially decaying thicker tails, such as a Laplace distribution, can easily replace a Gaussian prior for the same optimality under relaxed conditions.
Although such a prior may loosen a restriction on $\lVert f_0\rVert_\infty$ \citep{rockova2019semi,jeong2021posterior}, we primarily consider normal priors throughout the paper, even for non-Gaussian models for the sake of simplicity.
We consider non-Gaussian priors only when required for theoretical purposes; see, for example, a truncated prior for regression with random design in Section~\ref{sec:furapp}.

\section{Approximating the True Function}
\label{sec:appgen}
Recall that tree learners $f_{\mathcal E, B}$ in \eqref{eqn:treelearner} are piecewise constant, whereas the true function $f_0$ does not have to be. This will not be an issue as long as there exists a tree learner that can approximate  $f_0$ sufficiently well.
In this section, we establish the approximation theory for tree ensembles in the context of our targeted function spaces.

For isotropic classes, balanced $k$-d trees \citep{bentley1979multidimensional} are known to give rise to rate-optimal approximations under mild regularity conditions \citep{rockova2020posterior}.
This is not necessarily the case for our general setup where smoothness may vary over the domain and where cycling repeatedly through the coordinates (as is done in the $k$-d tree) may not be enough to capture localized features of $f_0$. We thus generalize the notion of $k$-d trees and show that there exists a good partitioning scheme for piecewise heterogeneous anisotropic classes. Although our primary interest lies in additive tree aggregations in \eqref{eqn:treelearner}, we show that a single deep tree can approximate well. We thereby consider only single trees $\mathcal T$ and suppress the superscript $t$ throughout this section.

\subsection{Split-Nets for Approximation}
\label{sec:splitnet}
Approximation properties of tree-based estimators are driven by the granularity and fineness of a chosen split-net. Roughly speaking, a good approximation requires that a split-net have two properties: (i) it should be dense enough so that the boundaries of the box partition $\mathfrak X_0^\ast=\{\Xi_1^\ast,\dots,\Xi_R^\ast\}$, extended from $\mathfrak X_0=\{\Xi_1,\dots,\Xi_R\}$, can be detected by a $\mathcal Z$-tree partition with a minimal error; and (ii) it should be regular enough so that there exists a $\mathcal Z$-tree partition that captures local/global features of $f_0$ on each $\Xi_r^\ast$.
We elucidate these two properties.

\subsubsection{Dense Split-Nets: Global Approximability}
\label{sec:dennet}

Recall that the underlying partition $\mathfrak X_0^\ast=\{\Xi_1^\ast,\dots, \Xi_R^\ast\}$ for the true function is unknown.
From the sheer flexibility of binary tree partitioning, we expect that the boundaries can be detected well enough by a $\mathcal Z$-tree partition if $\mathfrak X_0^\ast$ is a flexible tree partition.
If the prior rewards partitions that are sufficiently close to $\mathfrak X_0^\ast$, Bayesian CART (BART) is expected to adapt to unknown $\mathfrak X_0^\ast$ without much loss of efficiency. We examine when this adaptivity can be achieved in more detail below.

\begin{figure}[t!]
	\centering
	\resizebox{5in}{!}{
		\begin{tikzpicture}
			\draw  (-15.25,0) rectangle (-8.25,7);
			\draw  (-13.25,7) -- (-13.25,0);
			\draw  (-13.25,3) -- (-8.25,3);
			\node at (-14.25,3.5) {\Large$\Psi_1^1$};
			\node at (-10.75,5) {\Large$\Psi_2^1$};
			\node at (-10.75,1.5) {\Large$\Psi_3^1$};
			\draw  (-7.5,0) rectangle (-0.5,7);
			\draw [dashed]  (-4.75,7) -- (-4.75,0);
			\draw [dashed]  (-4.75,4) -- (-0.5,4);
			\node at (-6,3.5) {\Large$\Psi_1^2$};
			\node at (-2.75,5.25) {\Large$\Psi_2^2$};
			\node at (-2.75,2.25) {\Large$\Psi_3^2$};
			\draw  (1,0) rectangle (8,7);
			\draw  (3,7) -- (3,0);
			\draw  (3,3) -- (8,3);
			\draw  (1,0) rectangle (8,7);
			\draw [dashed]  (3.75,7) -- (3.75,0);
			\draw [dashed]  (3.75,4) -- (8,4);
			\draw [stealth-stealth,blue] (7.75,4) -- (7.75,3);
			\node [blue] at (6.5,3.5) {\Large$\Upsilon(\mathfrak Y^1,\mathfrak Y^2)$};
			\node at (0.25,3.5) {\LARGE$\Rightarrow$};
			\node at (-11.75,7.5) {\Large$\mathfrak Y^1$};
			\node at (-4,7.5) {\Large$\mathfrak Y^2$};
		\end{tikzpicture}
	}
	\caption{A two-dimensional example of the Hausdorff-type divergence in Definition~\ref{def:haus}. The divergence is the maximum dependency of the boxes in the partitions.}
	\label{fig:diver}
\end{figure}

The ability to detect $\mathfrak X_0^\ast$ is thus closely tied to the density of the split-net $\mathcal Z$; it should be dense enough so that a $\mathcal Z$-tree partition can be constructed that is sufficiently close to $\mathfrak X_0^\ast$.  Therefore, we need a gadget to measure the closeness between two partitions. To this end, we introduce a Hausdorff-type divergence; see Figure~\ref{fig:diver} for an illustration.

\begin{mydefinition}[Hausdorff-type divergence]
	For any two box partitions $\mathfrak Y^1=\{\Psi_1^1,\dots,\Psi_J^1\}$ and $\mathfrak Y^2=\{\Psi_1^2,\dots,\Psi_J^2\}$ with the same number $J$ of boxes, we define a divergence between $\mathfrak Y^1$ and $\mathfrak Y^2$ as
	\begin{align*}
		\Upsilon(\mathfrak Y^1,\mathfrak Y^2)
		=\min_{(\pi(1)\dots\pi(J))\in P_\pi[J]}\, \max_{1\le r\le J} \, {\rm Haus}(\Psi_r^1,\Psi_{\pi(r)}^2) ,
	\end{align*}
	where $P_\pi[J]$ denotes the set of all permutations $(\pi(1)\dots\pi(J))$ of $\{1,\dots,J\}$ and ${\rm Haus}(\cdot,\cdot)$ is the Hausdorff distance.
	\label{def:haus}
\end{mydefinition}

The permutation in Definition \ref{def:haus} makes the specification immune to the ordering of boxes.
We want the split-net $\mathcal Z$ to produce a  $\mathcal Z$-tree partition $\mathcal T$ such that $\Upsilon(\mathfrak X_0^\ast,\mathcal T)$ is smaller than some threshold.
Section~\ref{sec:app} establishes how small these thresholds should be  so that the tree learner is close to $f_0$  (for various approximation metrics). 
The following definition  will  be useful in characterizing the details.

\begin{mydefinition}[Dense split-net]
	For a given subset $S\subseteq\{1,\dots,p\}$ and an integer $J\ge1$, consider an $S$-chopped partition $\mathfrak Y=\{\Psi_1,\dots,\Psi_J\}$ of $[0,1]^p$ with boxes $\Psi_r\subseteq [0,1]^p$, $r=1,\dots,J$. For any given $c_n\ge0$, a split-net $\mathcal Z=\{z_i\in[0,1]^p,i=1,\dots, b_n\}$ is said to be $(\mathfrak Y, c_n)$-dense if  there exists an $S$-chopped $\mathcal Z$-tree partition $\mathcal T=\{\Omega_1,\dots,\Omega_J\}$ of $[0,1]^p$ such that $\Upsilon(\mathfrak Y,\mathcal T)\le c_n$.
	\label{def:den-net}
\end{mydefinition}

In Section~\ref{sec:app}, the approximation theory will require that $\mathcal Z$ be $(\mathfrak X_0^\ast, c_n)$-dense for some suitable $c_n\ge0$. 
Note that the ideal case $c_n=0$ can be achieved only when $\mathfrak X_0^\ast$ is a $\mathcal Z$-tree partition.
This condition, while obviously satisfied in the case $R=1$, is very restrictive in the most situations. This is because, if $J=1$, i.e., $\mathfrak Y=\{[0,1]^p\}$, we obtain $\Upsilon(\mathfrak Y,\mathcal T)=0$ for $\mathcal T=\{[0,1]^p\}$. Hence, every split-net $\mathcal Z$ is $(([0,1]^p),0)$-dense.	
However, we will see  in Theorem~\ref{thm:approx2} that, in many cases, it is sufficient that $c_n$ tends to zero at a suitable rate.
This means that $\mathfrak X_0^\ast$ should be at least a flexible tree partition, but not necessarily a $\mathcal Z$-tree partition. If $\mathfrak X_0^\ast$ is a box partition but not a flexible tree partition, we can redefine $\mathfrak X_0^\ast$  by adding more splits to make it a flexible tree partition. 
For example, the non-tree box partition in Figure~\ref{fig:partition} can be extended to a tree partition with a single extra split. 
However, this approach increases $R$ and hence may deteriorate the result (observe that our rate in \eqref{eqn:rate} is dependent on $R$). In particular, if $\mathfrak X_0^\ast$ is not a box partition (e.g., jumps are not axis-parallel), the redefined $R$ increases to infinity. For our theory to be valid, $\mathfrak X_0^\ast$ must be at least a box partition.
In Section~\ref{sec:exmspnet}, we present some examples of dense split-nets.

Dense split-nets have nested properties. That is, a $(\mathfrak Y, c_n)$-dense split-net is also $(\mathfrak Y, \tilde c_n)$-dense for every $\tilde c_n\ge c_n$. We are interested in the smallest possible $c_n$. In particular, every split-net $\mathcal Z$ is $(\mathfrak Y, 1)$-dense for any box partition $\mathfrak Y$.

\subsubsection{Regular Split-Nets: Local Approximability}
\label{sec:regular}
Beyond closely tracking smoothness boundaries,  good  tree partitions should be able to capture local/global  smoothness features of $f_0$.
In other words, there should exist a $\mathcal Z$-tree partition that achieves an optimal approximation error determined by our target rate.
In Section~\ref{sec:dennet}, we focused on {\em global} approximability of underlying partitions, which requires split-nets to be suitably dense. Now, we focus on {\em local} approximability.

Assume that $\mathfrak X_0^\ast$ can be approximated well (as discussed in the previous section) by an $S_0^\ast$-chopped $\mathcal Z$-tree partition $\mathcal T^\ast=\{\Omega_1^\ast,\dots,\Omega_R^\ast\}$,\footnote{
	The notation $\mathcal T^\ast=\{\Omega_1^\ast,\dots,\Omega_R^\ast\}$ with an asterisk is only used to denote an $S_0^\ast$-chopped $\mathcal Z$-tree partition approximating $\mathfrak X_0^\ast=\{\Xi_1^\ast,\dots,\Xi_R^\ast\}$.
}
which is formally written as
\begin{align}
	\mathcal T^\ast=\underset{\mathcal T \in\mathscr T_{S_0^\ast,R,\mathcal Z}}{\text{arg\,min}}  \Upsilon(\mathfrak X_0^\ast,\mathcal T).
	\label{eqn:treestar}
\end{align}
We now focus on local approximability inside each box $\Omega_r^\ast$. Ideally, one would want to construct a sub-tree partition of this local box that balances out approximation errors in all coordinates.
Therefore, we first need to devise a splitting scheme to achieve this balancing condition. The regularity of split-nets can then be spelled out based on such a law.

We  now zoom onto a single box $\Omega_r^\ast$. Recall that the true function $f_0$ has anisotropic smoothness on each of $\Xi_r^\ast$.
Intuitively, denser subdivisions are required for less smooth coordinates to capture the local features. 
Allowing  splits to occur more often in certain directions, we define the {\em anisotropic $k$-d tree}, which achieves the desired approximation error for anisotropic smoothness. The definition requires the notion of {\em midpoint-splits} defined as follows. For a given box $\Psi$ and a splitting coordinate $j$, a midpoint-split picks up the $\lceil \tilde b_j(\mathcal Z,\Psi)/2 \rceil$th split-candidate in $[\mathcal Z]_j\cap {\mathsf{int}}([\Psi]_j)$ as a split-point $\tau_j$, where  $\tilde b_j(\mathcal Z,\Psi)$ is the cardinality of $[\mathcal Z]_j\cap {\mathsf{int}}([\Psi]_j)$. 

\begin{mydefinition}[Anisotropic $k$-d tree]
	Consider a smoothness vector $\alpha=(\alpha_1,\dots,\alpha_d)^\top\in(0,1]^d$, a box $\Psi\subseteq[0,1]^p$, a split-net  $\mathcal Z=\{z_i\in[0,1]^p,~i=1,\dots,b_n\}$, an integer $L >0$, and an index set $S=\{s_1,\dots, s_d\}\subseteq\{1,\dots,p\}$ with $|S|=d$.
	We define the {\em anisotropic $k$-d tree} $\mathsf{Akd}(\Psi;\mathcal Z,\alpha,L,S)$ as the iterative splitting procedure that partitions $\Psi$ into disjoint boxes as follows. 
	\begin{enumerate}
		\item Start from the root node by setting $\Omega_1^\circ=\Psi$ and set $l_j=0$, $j=1,\dots, d$.
		\item \label{sp:2} For splits at iteration $1+\sum_{j=1}^d l_j$, choose $j$ corresponding to the smallest $l_j\alpha_j$. If the smallest $l_j\alpha_j$ is duplicated with multiple $j$s, choose the smallest $j$ among such $j$'s.
		\item For all boxes $\Omega_k^\circ$, $k=1,\dots, 2^{\sum_{j=1}^d l_j}$, at the current iteration, do the midpoint-splits with the given $\mathcal Z$ and the splitting coordinate $s_j$ chosen by $j$. Relabel the generated new boxes as $\Omega_k^\circ$, $k=1,\dots, 2^{1+\sum_{j=1}^d l_j}$, and then increase $l_j$ by one for chosen $j$.
		\item Repeat 2--3 until either $\sum_{j=1}^d l_j=L$ or the midpoint-split is no longer available. Return $(l_1,\dots,l_d)^\top$ and $\mathcal T^\circ=\{\Omega_1^\circ,\dots, \Omega_{2^{L^\circ}}^\circ\}$, where $L^\circ = \sum_{j=1}^d l_j$.
	\end{enumerate}
	\label{def:sp}
\end{mydefinition}

\begin{figure}[t!]
	\centering
	\includegraphics[width=2.5in]{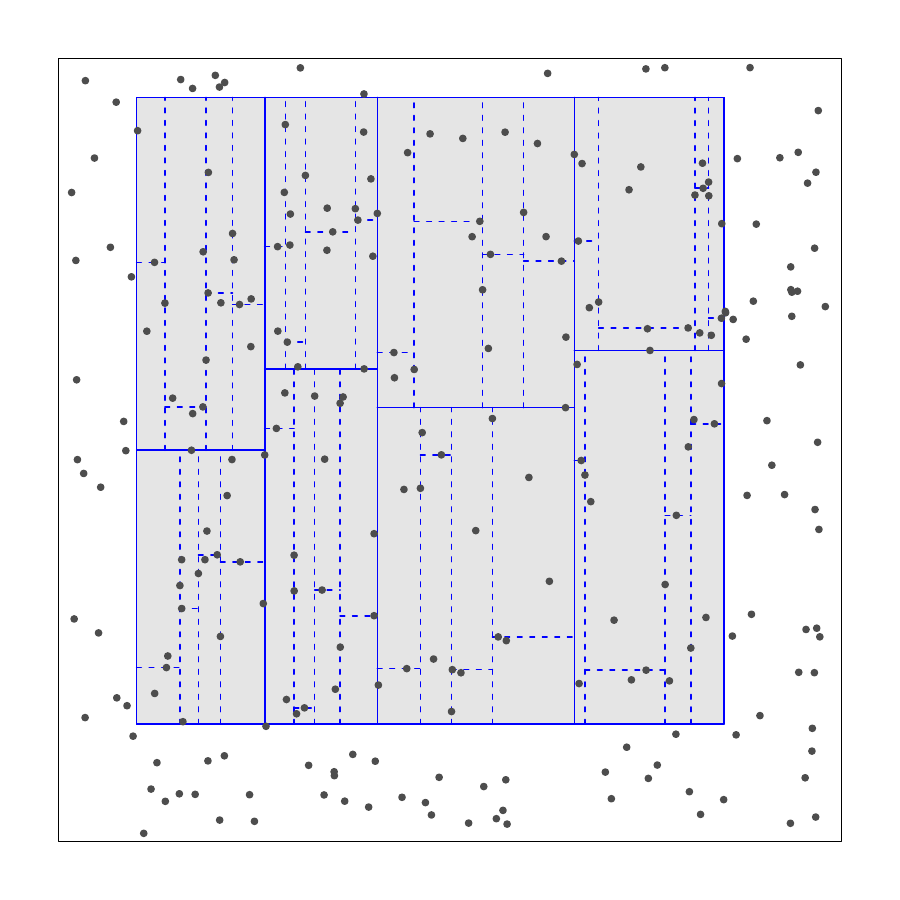}\vspace*{-0.1in}
	\caption{A realization of the anisotropic $k$-d tree with smoothness parameters $\alpha_1=0.25$ (for the horizontal axis) and $\alpha_2=0.5$ (for the vertical axis), and a box $\Psi$ (the shaded box) that is a subspace of $[0,1]^2$ (the outer square). Because $2\alpha_1=\alpha_2$, the subset $\Psi$ splits twice as often in the vertical direction than in the horizontal direction. }
	\label{fig:akd}
\end{figure}
Note that the anisotropic $k$-d tree construction depends on the smoothness that is unknown. Rather than a practical estimator, we use this to show that there exists a good tree approximator in the technical proof.
One possible realization of the anisotropic $k$-d tree generating process is given in Figure~\ref{fig:akd}. 
Observe that $\mathsf{Akd}(\Psi;\mathcal Z,\alpha,L,S)$ returns a tree partition $\mathcal T^\circ=\{\Omega_1^\circ,\dots,\Omega_{2^{L^\circ}}^\circ\}$ of $\Psi$ and a vector $(l_1,\dots,l_d)^\top$ such that $L^\circ=\sum_{j=1}^d l_j\le L$.\footnote{
	The notation $\mathcal T^\circ = \{\Omega_k^\circ\}_k$ with a circle is used only for tree partitions of some box $\Psi\subseteq[0,1]^p$, returned by the anisotropic $k$-d trees, with some suitable subscript if required.
} 
Although these returned items clearly depend on the inputs of the anisotropic $k$-d tree procedure (i.e., $\Psi$, $\mathcal Z$, $\alpha$, $L$, and $S$), we suppress them throughout the paper.
Each $l_j$ is a counter of how many times the $j$th coordinate has been used. The procedure is designed so that every $l_j$ is approximately proportional to $\alpha_j^{-1}$ after enough iterations. The total number of splits for the $j$th coordinate is thus close to $2^{C/\alpha_j}$ for every $j$ with some $C>0$. In the proof of Theorem~\ref{thm:approx2}, this matching is indeed clearly optimal and minimizes the induced bias.

To play a role as a `sieve' for approximation, $\Psi$ needs to be sufficiently finely subdivided to capture the global/local behavior of a function.
The threshold $L$ determines the resolution of the returned tree partition $\mathcal T^\circ=\{\Omega_1^\circ,\dots,\Omega_{2^{L^\circ}}^\circ\}$. 
For a good approximation, we are particularly interested in the situation when ${L^\circ}=L$, i.e., the resulting tree has the desired depth.
If ${L^\circ}<L$ owing to insufficient split-candidates, the resolution may not be good enough.

Now, we can define the regularity of a split-net on $\Psi\subseteq[0,1]^p$ using $\mathcal T^\circ$.
The desirable situation is  when all the splits occur nearly at the center of boxes such that, for any given $j\in S$, all $\mathsf{len}([\Omega_k^\circ]_j)$, $k=1,\dots,2^L$, are balanced well. 
The evenness of the returned partition is solely determined by the regularity of a split-net $\mathcal Z$.
Intuitively, the split-net should be regularly distributed to give rise to an appropriate partition, in which we say a split-net is regular. We make the definition technically precise below, which will be used as a basis for approximating the function classes. See \citet{verma2009spatial} for a related regularity condition.

\begin{mydefinition}[Regular split-net]
	For a given box $\Psi\subseteq[0,1]^p$, an integer $L>0$, and an index set $S=\{s_1,\dots,s_d\}\subseteq\{1,\dots,p\}$, we say that
	a split-net $\mathcal Z$ is  $(\Psi,\alpha,L, S)$-regular if  $\mathcal T^\circ = \{\Omega_1^\circ,\dots,\Omega_{2^{L^\circ}}^\circ\}$ and $ (l_1,\dots,l_d)^\top$, returned by $\mathsf{Akd}(\Psi;\mathcal Z,\alpha,L,S)$, satisfy $L^\circ=L$ and
	$\max_k\mathsf{len}([\Omega_k^\circ]_{s_j})\lesssim  \mathsf{len}([\Psi]_{s_j}) 2^{-l_j}$ for every $j=1,\dots, d$.
	\label{def:reg-net}
\end{mydefinition}

The condition $\max_k\mathsf{len}([\Omega_k^\circ]_{s_j})\lesssim \mathsf{len}([\Psi]_{s_j}) 2^{-l_j}$ is the key to obtaining optimal approximation results. In the ideal case that all the splits occur exactly at the center, this condition is trivially satisfied as $\max_k\mathsf{len}([\Omega_k^\circ]_{s_j})= \mathsf{len}([\Psi]_{s_j}) 2^{-l_j}$.
The inequality provides a lot more flexibility where the condition can be satisfied in most cases except for very extreme situations.
See Section~\ref{sec:exmspnet} for examples of regular split-nets.

Similar to dense split-nets, regular split-nets also have nested properties. If a split-net $\mathcal Z$ is $(\Psi,\alpha,L, S)$-regular for some $\Psi$, $\alpha$, $L$, and $S$, then it is also $(\Psi, \alpha,\tilde L, S)$-regular for any $\tilde L\le L$. This can be easily shown by noting that the latter is determined only by a pruned tree of the full-blown tree for the former. We are particularly interested in the largest possible $L$.

\begin{myremark}
	\label{rmk:reg}
	As regular split-nets require the desired depth, i.e., $L^\circ=L$, it is of interest to see which $L$ achieves this precondition.
	Consider a box $\Psi\subseteq[0,1]^p$ and a split-net $\mathcal Z=\{z_i\in[0,1]^p,~i=1,\dots,b_n\}$.
	If there are no ties in $\mathcal Z$ for any coordinate, i.e., $b_j(\mathcal Z)=b_n$, $j=1,\dots,p$, it can be easily checked that any integer $L\le\lfloor\log_2(\tilde b_j(\mathcal Z;\Psi)+1)\rfloor$ gives rise to $L^\circ=L$ with the anisotropic $k$-d tree. (Observe that all $\tilde b_j(\mathcal Z;\Psi)$ are identical in this case.) If there are ties, $L$ may need to be much smaller to achieve $L^\circ=L$, but a tight upper bound may not be obtained for the general case.
\end{myremark}

\subsection{Approximation Theory}
\label{sec:app}
Our goal is to establish the contraction rate of the posterior distribution. The construction requires that tree learners be able to approximate functions in the spaces $\Gamma^{A_{\bar\alpha},d,p}_{\lambda}(\mathfrak X_0)$ and $\Gamma^{A_{\bar\alpha},d,p}_{\lambda}(\mathfrak X_0)\cap \mathcal C([0,1]^p)$ appropriately.
Here, we investigate the approximation properties for these sparse function spaces.

\begin{figure}[t!]
	\centering
	\resizebox{1.5in}{!}{
		\begin{tikzpicture}
			\draw[opacity=0]  (0,0) rectangle (6,6);
			\node at (2,4.45) {\large$\Xi_1^\ast$};
			\node at (4.8,4.45) {\large$\Xi_2^\ast$};
			\node at (1.6,1.55) {\large$\Xi_3^\ast$};
			\node at (4.4,1.55) {\large$\Xi_3^\ast$};
		\end{tikzpicture}
	}
	\hspace{-1.5in}\includegraphics[width=1.5in]{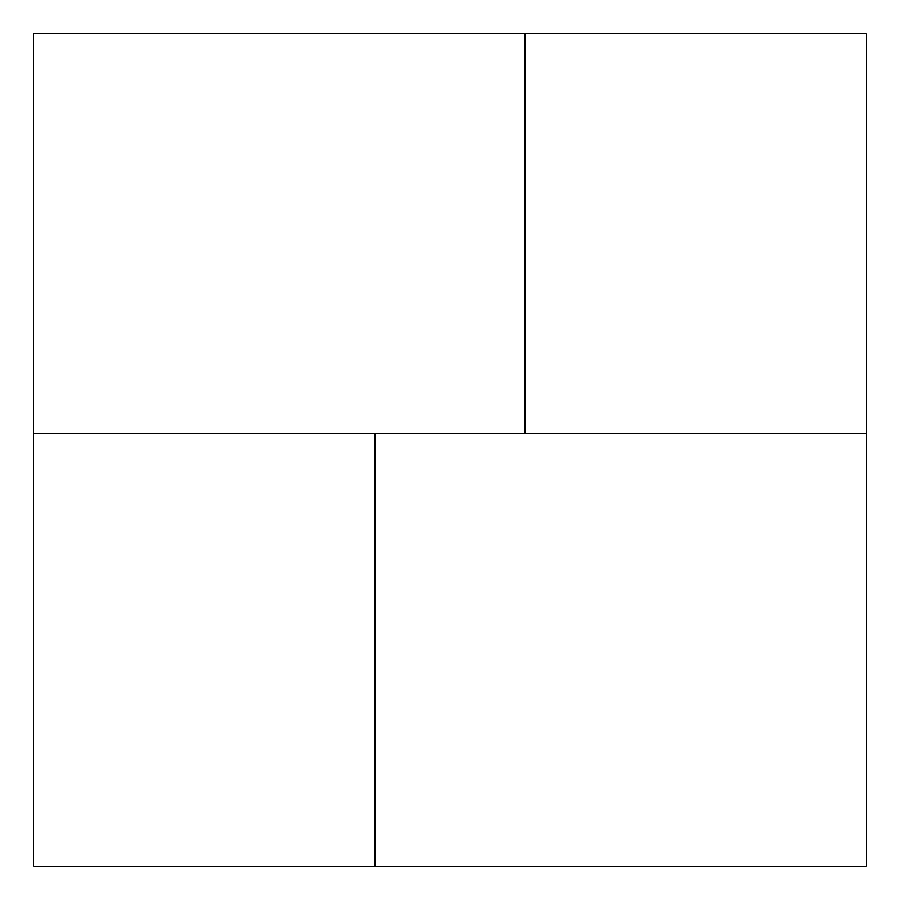}
	\resizebox{1.5in}{!}{
		\begin{tikzpicture}
			\draw[opacity=0]  (0,0) rectangle (6,6);
			\node at (2,4.45) {\large$\Omega_1^\ast$};
			\node at (4.8,4.45) {\large$\Omega_2^\ast$};
			\node at (1.6,1.55) {\large$\Omega_3^\ast$};
			\node at (4.4,1.55) {\large$\Omega_3^\ast$};
		\end{tikzpicture}
	}
	\hspace{-1.5in}\includegraphics[width=1.5in]{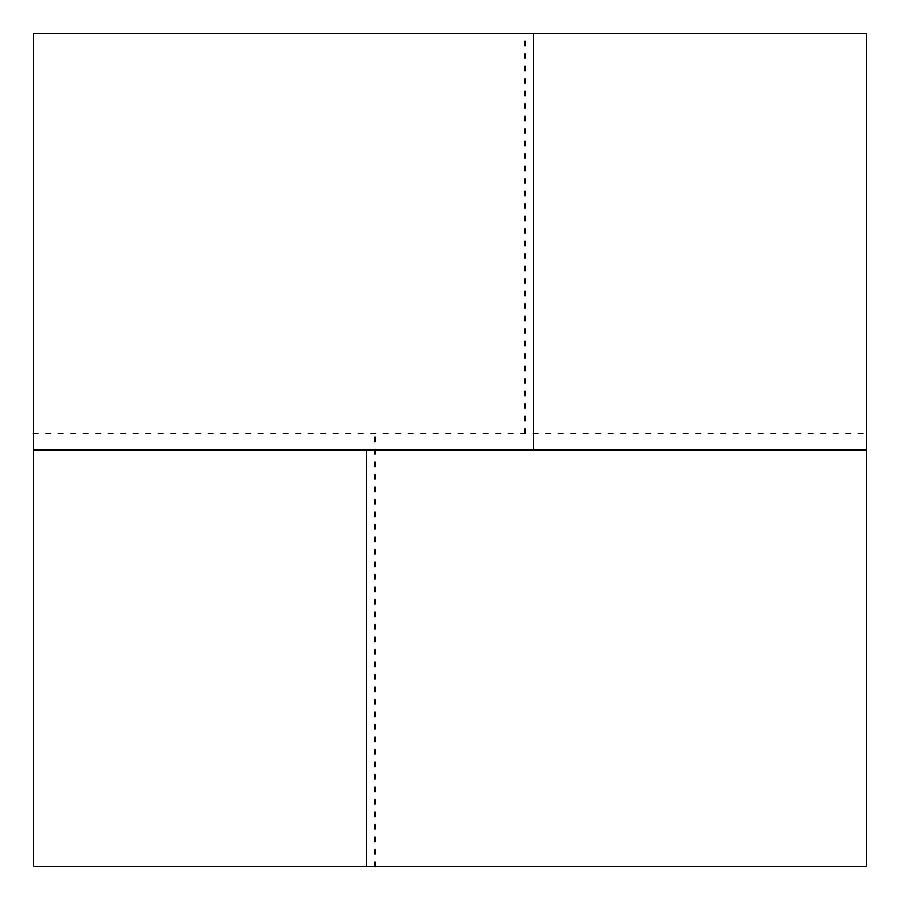}
	\resizebox{1.5in}{!}{
		\begin{tikzpicture}
			\draw[opacity=0]  (0,0) rectangle (6,6);
		\end{tikzpicture}
	}
	\hspace{-1.5in}\includegraphics[width=1.5in]{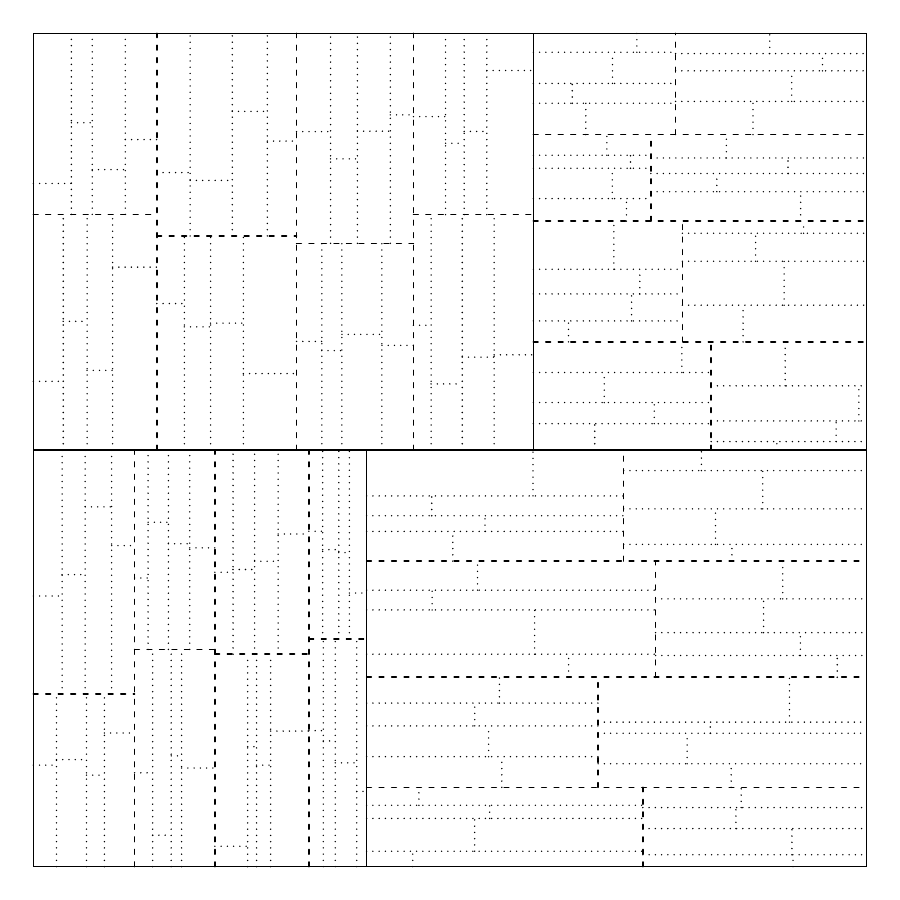}
	\caption{An example of constructing $\widehat{\mathcal T}$. First, $\mathfrak X_0^\ast=\{\Xi_1^\ast,\dots,\Xi_4^\ast\}$ is approximated by $\mathcal T^\ast=\{\Omega_1^\ast,\dots,\Omega_4^\ast\}$.
		Then, each $\Omega_r^\ast$ is subdivided by the anisotropic $k$-d tree, producing $\mathcal T_r^\circ$ a constituent of $\widehat{\mathcal T}$ displayed on the rightmost panel.}
	\label{fig:apptree}
\end{figure}

Recall that a split-net $\mathcal Z$ is required to be suitably dense and regular.
First, a split-net $\mathcal Z$ should be $(\mathfrak X_0^\ast, c_n)$-dense for some appropriate $c_n$, so that the boundaries of $\mathfrak X_0^\ast=\{\Xi_1^\ast,\dots, \Xi_R^\ast\}$ can be detected well by the binary tree partitioning rule. As $\mathfrak X_0^\ast$ is approximated by a $\mathcal Z$-tree partition with a given $\mathcal Z$, the underlying partition $\mathfrak X_0^\ast$ should be at least a flexible tree partition, but a stronger result is obtained if $\mathfrak X_0^\ast$ is a $\mathcal Z$-tree partition (see Theorem~\ref{thm:approx2} below).
Denoting by $\mathcal T^\ast=\{\Omega_1^\ast,\dots,\Omega_R^\ast\}$ the $S_0^\ast$-chopped $\mathcal Z$-tree partition in \eqref{eqn:treestar}, each box $\Omega_r^\ast$ should be appropriately subdivided to capture the local/global nature of the true function on $\Xi_r^\ast$. (If $R=1$, we write $\mathcal T^\ast=\mathfrak X_0^\ast=\{[0,1]^p\}$ with $\Omega_1^\ast=[0,1]^p$.)
Hence, for a smoothness parameter $A_{\bar\alpha}\in\mathcal A_{\bar\alpha}^{R,d}$ and some suitably chosen $L>0$, $\mathcal Z$ should also be $(\Xi_r^\ast,\alpha_r,L,S_0)$-regular, $r=1,\dots,R$. The integer sequence $L$ will eventually be chosen such that the approximation error is balanced with our target rate (see $L_0$ in Theorem~\ref{thm:approx2}).
Let $\mathcal T_r^\circ=\{\Omega_{r1}^\circ,\dots, \Omega_{r2^{L}}^\circ\}$ be the tree partition of $\Omega_r^\ast$ returned by $\mathsf{Akd}(\Omega_r^\ast;\mathcal Z,\alpha_r,L,S_0)$, $r=1,\dots,R$.
Then, the approximating partition $\widehat{\mathcal T}$ is formed by agglomerating all sub-tree partitions $\mathcal T_r^\circ$, leading to an $S_0$-chopped $\mathcal Z$-tree partition 
\begin{align}
\widehat{\mathcal T}=\left\{\Omega_{11}^\circ,\dots, \Omega_{12^{L}}^\circ,\dots, \Omega_{R1}^\circ,\dots, \Omega_{R2^{L}}^\circ\right\}.
\label{eqn:treeapp}
\end{align}
(Note that each $\mathcal T_r^\circ$ is $S_0$-chopped, not $S_0^\ast$-chopped.) A graphical illustration of constructing $\widehat{\mathcal T}$ is given in Figure~\ref{fig:apptree}.

The strongest approximation results relative to the $L_\infty$-norm for $\Gamma^{A_{\bar\alpha},d,p}_{\lambda}(\mathfrak X_0)$ are of particular interest.
Owing to the possible discontinuity or heterogeneity at the unknown boundaries of $\mathfrak X_0^\ast$, however, such results are not practically obtained except for the case $R=1$. 
As the following theorem shows, the conditions can be relaxed if we opt for weaker metrics, which often suffice in many statistical setups.
For example, in our examples of Gaussian nonparametric regression in Section~\ref{sec:nonregcon}, we only need an approximation rate in $L_2$- or empirical $L_2$-sense. The approximation results for the continuous variant $\Gamma^{A_{\bar\alpha},d,p}_{\lambda}(\mathfrak X_0)\cap \mathcal C([0,1]^p)$ require even milder conditions.

\begin{mytheorem}[Approximation theory]
For $c_n\ge0$ specified below, assume that a split-net $\mathcal Z$ is $(\mathfrak X_0^\ast,c_n)$-dense. For a smoothness parameter $A_{\bar\alpha}\in \mathcal A_{\bar\alpha}^{R,d}$ and an integer $L>0$, assume that $\mathcal Z$ is $(\Xi_r^\ast,\alpha_r,L, S_0)$-regular 
	for every $r=1,\dots, R$. Let $\tilde\epsilon_n$ be a sequence satisfying $\tilde\epsilon_n\gtrsim \lambda d 2^{-\bar\alpha L/d} $ and construct  $\widehat{\mathcal T}$ as in \eqref{eqn:treeapp} (through $\mathcal T^\ast$ in \eqref{eqn:treestar}). 
		Then, for any $f_0\in \Gamma^{A_{\bar\alpha},d,p}_{\lambda}(\mathfrak X_0)$, there exists  $\hat f_0\in\mathcal F_{\widehat{\mathcal T}}$ such that 
	\begin{enumerate}[label=\rm(\roman*)]
		\item \label{lmm0:sup1} $\lVert f_0 - \hat f_0\rVert_\infty\lesssim \tilde\epsilon_n$ if $c_n=0$;
		\item \label{lmm0:Leb1} $\lVert f_0 - \hat f_0\rVert_{v}\lesssim \tilde\epsilon_n$ if $c_n  \lesssim (\tilde\epsilon_n/\lVert f_0\rVert_\infty)^v  \min_{r,j}\mathsf{len}([\Xi_r^\ast]_j)/|S_0^\ast|$ for any $v\ge1$;
		\item \label{lmm0:emp1} $\lVert f_0 - \hat f_0\rVert_{v,P_{\mathcal Z}}\lesssim \tilde \epsilon_n$ for any $v\ge1$, where $P_{\mathcal Z}(\cdot)=b_n^{-1} \sum_{i=1}^{b_n} \delta_{z_i}(\cdot)$.
	\end{enumerate}
	Further, for any $f_0\in\Gamma^{A_{\bar\alpha},d,p}_{\lambda}(\mathfrak X_0)\cap \mathcal C([0,1]^p)$, there exists  $\hat f_0\in\mathcal F_{\widehat{\mathcal T}}$ such that
	\begin{enumerate}[label=\rm(\roman*),resume]
		\item \label{lmm0:sup2} $\lVert f_0 - \hat f_0\rVert_\infty\lesssim \tilde\epsilon_n$ if $c_n^{\min_{r,j}\alpha_{rj}}\lesssim \tilde\epsilon_n/(\lambda |S_0^\ast|)$;
		\item \label{lmm0:Leb2} $\lVert f_0 - \hat f_0\rVert_{v}\lesssim \tilde\epsilon_n$ if $c_n^{1+v\min_{r,j}\alpha_{rj}} \lesssim (\tilde\epsilon_n/\lambda )^v  \min_{r,j}\mathsf{len}([\Xi_r^\ast]_j)/|S_0^\ast|^{v+1}$ for any $v\ge1$.
	\end{enumerate}
In particular, if we choose $L=L_0$ such that $2^{L_0}\asymp (n(\lambda d)^2/(R\log n))^{d/(2\bar\alpha+d)}$, then the above assertions hold for $\tilde\epsilon_n=\bar\epsilon_n := (\lambda d)^{d/(2\bar\alpha+d)}(({R\log n})/{n})^{\bar\alpha/(2\bar\alpha+d)}$.	
	\label{thm:approx2}
\end{mytheorem}

\begin{myproof}
	See Section~\ref{sec:prooflmm1} in Appendix.
\end{myproof}

Although Theorem~\ref{thm:approx2} holds for any $\tilde\epsilon_n\gtrsim \lambda d 2^{-\bar\alpha L/d} $, the results are particularly useful for our purposes when combined with $L_0$ and $\bar\epsilon_n$, motivated by our target rate $\epsilon_n$ in \eqref{eqn:rate}.
The assertion in \ref{lmm0:sup1} gives the strongest result with the $L_\infty$-norm.  However, the condition $c_n=0$ requires that the boundaries of the pieces be correctly detectable by the binary tree partitioning rule with a given split-net $\mathcal Z$; that is, $\mathfrak X_0^\ast$ should be a $\mathcal Z$-tree partition. 
Except for the case $R=1$, this limitation is too restrictive and impractical, as the locations of the boundaries are unknown (every split-net $\mathcal Z$ is $(\mathfrak X_0^\ast,0)$-dense if $R=1$). The assertion in \ref{lmm0:sup2} relaxes this limitation by means of the continuity restriction. We will use \ref{lmm0:sup1} and \ref{lmm0:sup2} for a density estimation problem in Section~\ref{sec:density}.

The assertions in \ref{lmm0:Leb1} and \ref{lmm0:Leb2} are with respect to the $L_v$-norm, $v\ge 1$, which is useful in many statistical setups. We note that, despite the continuity restriction, the condition for $c_n$ of \ref{lmm0:Leb2} is not always milder than that of \ref{lmm0:Leb1}. Indeed, the former is milder than the latter only if $\lambda|S_0^\ast|c_n^{\min_{r,j}\alpha_{rj}}\lesssim \lVert f_0\rVert_\infty$, which is often satisfied, as the left-hand side is prone to be decreasing with a suitably chosen $c_n$.
We will use the results in \ref{lmm0:Leb1} and \ref{lmm0:Leb2} for nonparametric regression and binary classification with random design in Sections~\ref{sec:nonran} and \ref{sec:classification}.

The assertion in \ref{lmm0:emp1} is particularly useful in regression setups with $\mathcal Z$ chosen by fixed covariates; see Sections~\ref{sec:fixeddesign} and \ref{sec:addnonpreg}. Note that \ref{lmm0:emp1} only explicitly requires the regularity of a split-net $\mathcal Z$, and an upper bound for $c_n$ is not specified. This is because the closeness between $f_0$ and $\hat f_0$ is measured only at points in $\mathcal Z$, and the boundary detection needs to be performed much loosely compared with the other metrics. Although not explicitly stated, \ref{lmm0:emp1} still requires a dense split-net in an implicit way. Indeed, every assertion in Theorem~\ref{thm:approx2} necessitates a condition on $\mathcal T^\ast$ imposed implicitly by the regularity with $\Xi_r^\ast$; for $\mathcal Z$ to be regular for every $\Xi_r^\ast$, it must be sufficiently evenly distributed and hence suitably dense.

As stated above,  if $R=1$, i.e., the global anisotropic case, we always obtain the strongest result in \ref{lmm0:sup1} as soon as a split-net is suitably regular.
If $R>1$, a split-net should also be suitably dense except for the case of the empirical $\lVert \cdot\rVert_{v,P_{\mathcal Z}}$-norm in \ref{lmm0:emp1}.
As the conditions on $c_n$ depend on unknown model specification, e.g., $A_{\bar\alpha}\in\mathcal A_{\bar\alpha}^{R,d}$, $\lambda$, and $|S_0^\ast|$, more practical conditions can be obtained by plugging in reasonable bounds of the unknown components.
For example, we cannot hope for better than $\bar\epsilon_n\gtrsim (\lambda d R (\log n )/n)^{1/3} $ owing to the fundamental limitation of piecewise constant learners. We can also assume that $\min_{r,j}\mathsf{len}([\Xi_r^\ast]_j)$ is bounded away from zero or decreases at most polynomially. To establish the posterior contraction rate, we will eventually assume  $\lVert f_0\rVert_\infty\lesssim \sqrt{\log n}$ (see \ref{asm:fsup} below).
Because the necessary conditions $d/\bar\alpha\ll \log n$ and $\lambda^{\bar\alpha/d}R\ll n$ are required for consistent estimation (see the rate in \eqref{eqn:rate} below), making mild assumptions on $d$ and $\lambda$ is not prohibitive (note that $|S_0^\ast|\le d$).
Putting everything together, the conditions on $c_n$ can be easily satisfied if $c_n$ is a decreasing polynomial in $n$ with a suitable exponent.
The results are formalized in the following corollary.

	\begin{mycorollary}[Approximation with $L_\infty$ and $L_v$ when $R>1$]
		Under the setup of Theorem~\ref{thm:approx2} with $L=L_0$, suppose that $R>1$ and $d\lesssim \log n$. Then, the following assertions hold.
		\begin{enumerate}[label=\rm(\roman*)]
			\item \label{lmm0:sup2-1} Suppose that $\min_{r,j}\alpha_{rj}\ge a_1$ and $\lambda\lesssim n^{a_2}$ for some constants $a_1>0$ and $a_2\ge 0$.
If $c_n\lesssim n^{-(1+2a_2)/(3a_1)}(\log n)^{-1/(3a_1)}$, then for every $f_0\in\Gamma^{A_{\bar\alpha},d,p}_{\lambda}(\mathfrak X_0)\cap \mathcal C([0,1]^p)$, there exists  $\hat f_0\in\mathcal F_{\widehat{\mathcal T}}$ such that 	$\lVert f_0 - \hat f_0\rVert_\infty\lesssim \bar \epsilon_n$.
			\item \label{lmm0:Leb1-1} Suppose that $\lVert f_0\rVert_\infty\lesssim \sqrt{\log n}$ and 			$\min_{r,j}\mathsf{len}([\Xi_r^\ast]_j)\gtrsim n^{-a_3}$ for some constant $a_3\ge 0$.
			Fix any $v\ge 1$. If $c_n\lesssim n^{-(v/3+a_3)}(\log n)^{-(\max\{0,1-v/3\}+v/6)}$, then for every $f_0\in\Gamma^{A_{\bar\alpha},d,p}_{\lambda}(\mathfrak X_0)$, there exists $\hat f_0\in\mathcal F_{\widehat{\mathcal T}}$ such that $\lVert f_0 - \hat f_0\rVert_{v}\lesssim \bar \epsilon_n$.
		\end{enumerate}
		\label{cor:approx}
	\end{mycorollary}

\begin{myproof}
	See Section~\ref{sec:prooflmm1} in Appendix.
\end{myproof}

Corollary~\ref{cor:approx} implies that the target approximation error is attained with both the $L_\infty$- and $L_v$-norms as soon as $c_n$ decreases polynomially. The assertion in \ref{lmm0:sup2-1} provides the stronger result with the aid of the continuous restriction. It also requires a constant lower bound of the minimum smoothness parameter $\min_{r,j}\alpha_{rj}$, causing $\bar\alpha$ to be bounded away from zero. In contrast, \ref{lmm0:Leb1-1} removes such a restriction at the expense of a tighter upper bound. 
 In general, the conditions for \ref{lmm0:Leb1-1} are much milder, yielding a relatively weaker but still useful result in many statistical setups.

 \begin{myremark}
 	No upper bounds for $b_n$ and $b_j(\mathcal Z)$ are made for Theorem~\ref{thm:approx2}; the approximation results are more easily achieved with larger values of $b_j(\mathcal Z)$, $j=1,\dots, p$.
 	However, values increasing too fast may harm the contraction rate as they escalate the model complexity.
 	In Section~\ref{sec:nonreg}, we will see that our main results on the optimal posterior contraction require that $\max_{1\le j\le p}\log b_j(\mathcal Z)\lesssim \log n$.
 	We are ultimately interested in well-balanced split-nets.
 	\label{lmk:upper}
 \end{myremark}

 \begin{myremark} Our approximation theory is presented with the error $\bar\epsilon_n$ motivated by our target rate $\epsilon_n$ in \eqref{eqn:rate}. However, what we really need is the weaker approximation error $\epsilon_n$, which is identical to the posterior contraction rate (see Sections~\ref{sec:nonreg}--\ref{sec:furapp}). Although the latter slightly relaxes the required conditions, we stick to the approximation result with $\bar\epsilon_n$ because such generalization complicates the technical details too much for a small gain. 
  \end{myremark}	 

 \begin{myremark} The assertion in \ref{lmm0:emp1} requires $\mathcal Z$ to be regular over $[0,1]^p$. Because the assertion is with respect to $L_v(\mathcal P_{\mathcal Z})$-norm, one may anticipate the regularity over $[0,1]^p$ to be relaxed into a smaller subset. Indeed, we can restrict our attention to a subset of $[0,1]^p$ and the technical details require the regularity $\mathcal Z$ only over such a smaller subset. We do not consider such an extension so that the $L_v(\mathcal P_{\mathcal Z})$-consistency can be interpreted as an approximate result for the $L_v$-norm, which is more appealing in the usual sense.
\end{myremark}	 

\subsection{Examples of Split-Nets for Approximation}
\label{sec:exmspnet}

Although the notion of dense and regular split-nets is crucial in characterizing the approximation theory in Section~\ref{sec:app}, how to obtain such a good split-net in practice remains unsolved.
Clearly, a split-net attains the suitable density and regularity more easily with larger $b_n$. As mentioned in Remark~\ref{lmk:upper}, however, we will see that a split-net must satisfy $\max_{1\le j\le p}\log b_j(\mathcal Z)\lesssim \log n$ to establish the optimal posterior contraction rate. Accordingly, our primary concern is examining split-nets that are suitably dense and regular under the restriction on $\log b_j(\mathcal Z)$. In this subsection, we show that the two split-nets described in Section~\ref{sec:partition} are dense and regular as required, and hence fulfill the requirements of Theorem~\ref{thm:approx2} and Corollary~\ref{cor:approx}.

\subsubsection{Regular Grid}
\label{sec:grid}

We first consider a regular grid $\mathcal Z = \{(i-1/2)/b_n^{1/p},i=1,\dots,b_n^{1/p}\}^p$ for $b_n$ such that $b_n^{1/p}$ is an integer.
This simplest example is a split-net according to Definition~\ref{def:split-net}.
We will see that a regular grid can be useful for density estimation, binary classification, and nonparametric regression with random design, but it also has the potential to be used for many other statistical models.
A two-dimensional example is illustrated in Figure~\ref{fig:spn1}.
The following lemma shows that, with an appropriately chosen $b_n$, a regular grid is suitably dense and regular under mild conditions.

\begin{mylemma}[Regular grid]
	Consider a regular grid $\mathcal Z$ with $b_n=n^{cp}$ for a constant $c\ge 1$. If $\min_{r,j}\mathsf{len}([\Xi_r^\ast]_j)\ge n^{-c}$ and $\lambda d /\min_{r,j}\mathsf{len}([\Xi_r^\ast]_j)^{\bar\alpha/d+1/2}\lesssim n^{c \bar\alpha /d+(c-1)/2}\sqrt{R \log n}$, then $\mathcal Z$ is $(\mathfrak X_0^\ast, c_n)$-dense and $(\Xi_r^\ast,\alpha_r,L_0,S_0)$-regular for $r=1,\dots,R$, where $c_n=n^{-c}\mathbbm 1(R>1)$.
	\label{lmm:grid2}
\end{mylemma}
\begin{myproof}
	See Section~\ref{sec:prooflmm2} in Appendix.
\end{myproof}

The second condition is replaced by $\lambda d /\sqrt{\min_{r,j}\mathsf{len}([\Xi_r^\ast]_j)}\lesssim n^{(c-1)/2}\sqrt{R\log n}$ if we consider the worst-case scenario $\bar\alpha\rightarrow 0$ with the upper bound $\bar\alpha/d\le 1$.
Combined with the necessary conditions $d/\bar\alpha\ll \log n$ and $\lambda^{\bar\alpha/d}R\ll n$ for consistent estimation (see \eqref{eqn:rate}), the conditions are very mild as soon as $c$ is suitably large.
 The choice $c=1$ may even be sufficient with stronger boundedness conditions, i.e., $\lambda\lesssim 1$, $d\lesssim \sqrt{\log n}$, and $\min_{r,j}\mathsf{len}([\Xi_r^\ast]_j)\gtrsim 1$. In particular, the first condition is trivially satisfied if $R=1$, i.e., $\mathfrak X_0^\ast=\{[0,1]^p\}$.
 In this case, we obtain the strongest result in \ref{lmm0:sup1} of Theorem~\ref{thm:approx2} as soon as the second condition is satisfied (recall that $\Gamma^{A_{\bar\alpha},d,p}_{\lambda}(\mathfrak X_0)=\Gamma^{A_{\bar\alpha},d,p}_{\lambda}(\mathfrak X_0)\cap \mathcal C([0,1]^p)$ if $R=1$).
 If $R>1$, $c_n$ is a decreasing polynomial in $n$ with our choice of $b_n$. This concludes that, with a suitably large $c$, the assertions in \ref{lmm0:Leb1} and \ref{lmm0:sup2} of Theorem~\ref{thm:approx2} (or the assertions in \ref{lmm0:sup2-1} and \ref{lmm0:Leb1-1} of Corollary~\ref{cor:approx}) hold. Note that \ref{lmm0:emp1} of Theorem~\ref{thm:approx2} also holds trivially with this $\mathcal Z$.

As $\max_{1\le j\le p}\log b_j(\mathcal Z)=p^{-1}\log b_n\lesssim \log n$, a regular grid satisfies the condition for the optimal posterior contraction specified in Section~\ref{sec:nonreg} (see Remark~\ref{lmk:upper}). This makes a regular grid very appealing for practical use given its simplicity, and there is little benefit of considering more complicated split-nets. The only exception is a set of fixed design points commonly used in the literature of BART \citep{chipman2010bart,rockova2020posterior}.

A regular grid can easily be extended to an irregular rectangular grid with boxes of different sizes.
If every mesh-size of an irregular checkerboard is asymptotically proportional to $1/b_n^{1/p}$,  the above results still hold with minor modification.
This extension is particularly interesting in a regression setup where the distribution of covariates is explicitly available. For example,  it allows us to use the quantiles for grid points, which is a natural way to generate a weakly balanced system \citep{castillo2021uncertainty}.

\begin{myremark}
	Lemma~\ref{lmm:grid2} indicates that a large value of $c$ is preferred in the sense of making the required conditions mild. Furthermore, a large $c$ does not harm the posterior contraction rate, as the boundedness condition $\max_{1\le j\le p}\log b_j(\mathcal Z)\lesssim \log n$ is satisfied for any $c>0$. Nonetheless, the empirical performance is affected by the size of $c$; an extremely large $c$ produces unnecessarily many split-candidates, making the algorithm inefficient. Consequently, we want to choose a suitable but not extremely large $c$.
	A good choice of $c$ is model-specific. In Section~\ref{sec:density}, we will see that density estimation requires approximation with respect to the $L_\infty$-norm, which can be fulfilled by \ref{lmm0:sup2-1} of Corollary~\ref{cor:approx} with the continuity assumption on $f_0$. If $\lambda\lesssim 1$, $d\lesssim \sqrt{\log n}$, $\min_{r,j}\mathsf{len}([\Xi_r^\ast]_j)\gtrsim 1$, and $\min_{r,j}\alpha_{rj}>1/3$, then $c=1$ and the corresponding $c_n$ satisfy the requirements for \ref{lmm0:sup2-1} of Corollary~\ref{cor:approx} and Lemma~\ref{lmm:grid2}. The most disappointing assumption is the lower bound for the minimum smoothness parameter, $\min_{r,j}\alpha_{rj}>1/3$. Although we recommend $c=1$ as the default choice by assuming such requirements, increasing $c$ is recommended if the density function is thought to be less smooth.\footnote{A careful examination of the proof indicates that the isotropy assumption eliminates the condition $\min_{r,j}\alpha_{rj}>1/3$, so $c=1$ works for all smoothness levels. This is because isotropy causes $\min_{r,j}\alpha_{rj}=\bar\alpha$, and there is enough cancellation in simplifying \ref{lmm0:sup2} of Theorem~\ref{thm:approx2}. To maintain anisotropy throughout the paper, we do not investigate such a particular situation in greater detail.} In contrast, nonparametric regression with random design and binary classification require approximation with respect to the $L_2$-norm (see Sections~\ref{sec:nonran} and \ref{sec:classification}), which is obtained by \ref{lmm0:Leb1-1} of Corollary~\ref{cor:approx}. One can easily verify that, if $\lambda\lesssim 1$, $d\lesssim \sqrt{\log n}$, and $\min_{r,j}\mathsf{len}([\Xi_r^\ast]_j)\gtrsim 1$, then $c=1$ and the corresponding $c_n$ satisfy the conditions for \ref{lmm0:Leb1-1} of Corollary~\ref{cor:approx} and Lemma~\ref{lmm:grid2}, and hence $c=1$ is the default choice.
\label{rmk:practical.c}
\end{myremark}

\subsubsection{Fixed Design Points}
\label{sec:fixeddesign}

Now we focus on a fixed design regression setup, where observed covariate values are readily available. In this case, using fixed design points is particularly appealing in that \ref{lmm0:emp1} of Theorem~\ref{thm:approx2} (coupled with this split-net) gives an approximation error relative to the empirical probability measure as soon as it is suitably regular (the assertion does not require a further bound on $c_n$).
The strategy is conventional in the literature of Bayesian CART and BART \citep{chipman1998bayesian,denison1998bayesian,chipman2010bart}. 

Suppose that a split-net $\mathcal Z=\{z_i\in[0,1]^p ,~i=1,\dots,n\}$ consists of the observed covariate values in a regression setup.
We need to assume that the design points are sufficiently evenly distributed in $S_0$.
The required assumption is formalized as follows.
\begin{enumerate}[leftmargin=2.0\parindent,label=\rm(F)]
	\item \label{asm:fixeddesign} For every $\alpha\in(0,1]^d$ and every box $\Psi\subseteq[0,1]^p$ with $n P_{\mathcal Z}(\Psi)\gg 1$,  $\mathcal Z$ is $(\Psi,\alpha,L,S_0)$-regular with $L = \lfloor\log_2 (c n P_{\mathcal Z}(\Psi))\rfloor$ for some constant $c>0$. 
\end{enumerate}
Although assumption \ref{asm:fixeddesign} may appear nontrivial, it is actually not restrictive. As $P_{\mathcal Z}$ is defined as $P_{\mathcal Z}(\cdot)=b_n^{-1} \sum_{i=1}^{b_n} \delta_{z_i}(\cdot)$, for $\mathcal Z$ chosen above, $n P_{\mathcal Z}(\Psi)$ denotes the number of split-candidates contained in $\Psi$. Hence, the condition $n P_{\mathcal Z}(\Psi)\gg 1$ implies that the number of design points in $\Psi$ increases with $n$, which is a certainly mild assumption.
As noted in Remark~\ref{rmk:reg}, if $\mathcal Z$ is balanced very well in $S_0$ and there are no ties so that splits can occur $n P_{\mathcal Z}(\Psi)$ times, then  $\mathcal Z$ is $(\Psi,\alpha,L,S_0)$-regular for $L=\lfloor\log_2 (n P_{\mathcal Z}(\Psi)+1)\rfloor$. Our requirement in \ref{asm:fixeddesign} is milder with the aid of the constant $c$.

\begin{mylemma}[Fixed design points]
	Consider fixed design points $\mathcal Z=\{z_i,i=1,\dots,n\}$ satisfying assumption \ref{asm:fixeddesign}.
	If $\lambda d \lesssim (n/R)^{\bar\alpha/d} \sqrt{\log n}$, $\min_r P_{\mathcal Z}(\Xi_r^\ast)\gtrsim R^{-1}$, and $R\ll n$, then $\mathcal Z$ is $(\Xi_r^\ast,\alpha_r,L_0,S_0)$-regular for $r=1,\dots,R$.
	\label{lmm:fixed2}
\end{mylemma}
\begin{myproof}
	See Section~\ref{sec:prooflmm2} in Appendix.
\end{myproof}

As $n P_{\mathcal Z}(\Xi_r^\ast)$ is the number of split-candidates in $\Xi_r^\ast$, the condition $\min_r P_{\mathcal Z}(\Xi_r^\ast)\gtrsim R^{-1}$ implies that the number of split-candidates should be balanced well among the $R$ boxes.
Our condition $\lambda d \lesssim (n/R)^{\bar\alpha/d} \sqrt{\log n}$ slightly relaxes the condition $\lambda d \lesssim \sqrt{\log n}$ of Theorem~4.1 in \citet{rockova2020posterior} (for the case of global isotropy). The latter is obtained if we consider the worst-case scenario $\bar\alpha\rightarrow 0$.
We see that \ref{lmm0:emp1} of Theorem~\ref{thm:approx2} directly follows from this lemma. As the design points are used as $\mathcal Z$, the term $\lVert f_0 - \hat f_0\rVert_{v,P_{\mathcal Z}}$ is translated into the approximation error relative to the empirical probability measure. In regression setups, this fact makes fixed design points much more attractive than other split-nets in the previous sections.
We also note that the requirement $\max_{1\le j\le p}\log b_j(\mathcal Z)\lesssim \log n$ for the optimal posterior contraction is trivially satisfied.

\section{BART in Nonparametric Regression}
\label{sec:nonreg}

\subsection{Posterior Contraction Rates}
\label{sec:nonregcon}
BART is an archetypal example of Bayesian forests \citep{chipman1998bayesian,denison1998bayesian,chipman2010bart}.
For a fixed design Gaussian nonparametric regression,
\citet{rockova2020posterior} and \citet{rockova2019theory} established $L_2$ rate-optimal posterior contraction of BART for high-dimensional isotropic regression functions.
Our investigation goes beyond these studies in three aspects: (i)  we treat the variance parameter $\sigma^2$ as unknown with a prior; (ii) we consider both fixed and random regression design; and, most importantly, (iii) the true function  is assumed to be in the piecewise heterogeneous anisotropic space introduced earlier. The last point significantly enlarges the optimality scope of BART.

We separately deal with fixed and random designs. This section is focused on the fixed design case, while the random design case will be considered in Section~\ref{sec:nonran}.
The fixed design regression model writes as
\begin{align}
	Y_i=f_0(x_i)+\varepsilon_i,\quad \varepsilon_i\sim \text{N}(0,\sigma_0^2),\quad i=1,\dots,n,
	\label{eqn:modelreg}
\end{align}
where $x_i=(x_{i1},\dots,x_{ip})^\top\in[0,1]^p$, $i=1,\dots, n$, are fixed. The model is independent but not identically distributed, and hence the asymptotic studies are established under the product measure for the $n$ observations.
The general theory of posterior contraction requires an exponentially powerful test function of a semimetric under this product measure  \citep{ghosal2017fundamentals}.
In nonparametric regression with fixed design, such a good test function can be directly constructed for the empirical $L_2$-distance even when the noise error is unknown \citep{ning2020bayesian,jeong2021unified,lim2023synergizing}.
The general theory also requires desirable properties of the prior. We show that the tree priors in Section~\ref{sec:prior} satisfy those conditions.

We impose the following assumptions on the true parameters $f_0$ and $\sigma_0^2$.
\begin{enumerate}[leftmargin=2.0\parindent,label=\rm(A\arabic*)]
	\item \label{asm:truef} For $d>0$, $\lambda>0$, $R>0$, $\mathfrak X_0=\{\Xi_1,\dots,\Xi_R\}$, and $A_{\bar\alpha}\in\mathcal A_{\bar\alpha}^{R,d}$ with $\bar\alpha\in(0,1]$, the true function satisfies $f_0\in\Gamma^{A_{\bar\alpha},d,p}_{\lambda}(\mathfrak X_0)$ or $f_0\in\Gamma^{A_{\bar\alpha},d,p}_{\lambda}(\mathfrak X_0)\cap \mathcal C([0,1]^p)$.
	\item \label{asm:dp} It is assumed that $d$, $p$, $\lambda$, $R$, and $\bar\alpha$ satisfy $\epsilon_n\ll 1$, where
	\begin{align}
		\epsilon_n=\sqrt{\frac{d\log p}{n}}+(\lambda d)^{d/(2\bar\alpha+d)}\left(\frac{R\log n}{n}\right)^{\bar\alpha/(2\bar\alpha+d)}.
		\label{eqn:rate}
	\end{align}
	\item \label{asm:fsup} The true function satisfies $\lVert f_0 \rVert_\infty\lesssim \sqrt{\log n}$.
	\item \label{asm:sig} The true variance parameter satisfies $\sigma_0^2\in[C_0^{-1},C_0]$ for a sufficiently large $C_0>1$.
\end{enumerate}
Assumption \ref{asm:truef} means that the true regression function $f_0$ lies on a sparse piecewise heterogeneous anisotropic space.
If the continuity assumption is further imposed, the approximation results in Theorem~\ref{thm:approx2} are obtained under milder conditions.
Assumption \ref{asm:dp} is required to make our target rate $\epsilon_n$ tend zero. 
The boundedness condition in \ref{asm:fsup} is made to guarantee a sufficient prior concentration under the normal prior on the step-heights specified in \ref{pri:normal} below.
Although the Gaussian prior can be replaced by a thick-tailed prior \citep[e.g.,][]{rockova2019semi}, we only consider the Gaussian prior to leverage its semi-conjugacy. Assumption \ref{asm:sig} allows one to assign a standard prior to $\sigma^2$, e.g., an inverse gamma distribution.

It is also important to choose a suitable split-net so that Theorem~\ref{thm:approx2} can be deployed.
For regression with fixed design, we need an approximation result with respect to the empirical $L_2$-norm  $\lVert\cdot\rVert_n$ defined as $\lVert f \rVert_n^2=n^{-1}\sum_{i=1}^n |f(x_i)|^2$.
We make the following assumptions on the split-net $\mathcal Z$. The notation $\mathsf{dep}$ means the depth of a node, the number of nodes along the path from the root node down to that node.

\begin{enumerate}[leftmargin=2.0\parindent,label=\rm(A\arabic*),resume]
	\item \label{asm:spnet0} The split-net $\mathcal Z$ satisfies $\max_{1\le j\le p}\log b_j(\mathcal Z)\lesssim\log n$.
	\item \label{asm:spnet} The split-net $\mathcal Z$  is suitably dense and regular to construct a $\mathcal Z$-tree partition $\widehat {\mathcal T}$ such that there exists $\hat f_0\in\mathcal F_{\widehat {\mathcal T}}$ satisfying $\lVert f_0-\hat f_0 \rVert_n\lesssim \bar\epsilon_n$ by Theorem~\ref{thm:approx2}. 
	\item \label{asm:spnet1} The $\mathcal Z$-tree partition $\mathcal T^\ast=\{\Omega_1^\ast,\dots,\Omega_R^\ast\}$ approximating $\mathfrak X_0^\ast$ satisfies $\max_r\mathsf{dep}(\Omega_r^\ast)\lesssim \log n$.
\end{enumerate}
Assumption \ref{asm:spnet0} is required for a suitable bound of the entropy and a good prior concentration (see Lemma~\ref{lmm:priorcon}).
Assumption \ref{asm:spnet} provides the desired approximation error with respect to the $\lVert\cdot\rVert_n$-distance. 
Owing to \ref{lmm0:emp1} of Theorem~\ref{thm:approx2} and Lemma~\ref{lmm:fixed2}, using fixed design points as $\mathcal Z$ is of particular interest, as $\lVert\cdot\rVert_{2,P_{\mathcal Z}}$ is equivalent to the empirical $L_2$-norm $\lVert\cdot\rVert_n$ in this case.
Assumption \ref{asm:spnet1} is a technical requirement which is certainly mild. This condition is trivially satisfied if $R$ is bounded.

Lastly, careful prior specification is required to obtain the optimal posterior contraction. We consider the following prior distributions discussed in Section~\ref{sec:prior}.
\begin{enumerate}[leftmargin=2.0\parindent,label=\rm(P\arabic*)]
	\item \label{pri:tree} For a fixed $T>0$, each tree $\mathcal T^t$, $t=1,\dots,T$, is independently assigned a tree prior with Dirichlet sparsity.
	\item \label{pri:normal} The step-heights $B$ are assigned a normal prior with a zero-mean and a covariance matrix whose eigenvalues are bounded below and above.
	\item \label{pri:invgam} The variance parameter $\sigma^2$ is assigned an inverse gamma prior.
\end{enumerate}

Under the above assumptions and priors, the following theorem formalizes the posterior contraction rate of model \eqref{eqn:modelreg}.

\begin{mytheorem}[Nonparametric regression, fixed design]
	Consider model \eqref{eqn:modelreg} with Assumptions \ref{asm:truef}--\ref{asm:spnet1} and the prior assigned through \ref{pri:tree}--\ref{pri:invgam}.
	Then, there exists a constant $M>0$ such that for $\epsilon_n$ in \eqref{eqn:rate},
	\begin{align*}
		\mathbb E_0 \Pi\Big\{(f,\sigma^2):\lVert f-f_0\rVert_n+|\sigma^2-\sigma_0^2|>M \epsilon_n \,\big|\, Y_1,\dots,Y_n\Big\}\rightarrow 0.
	\end{align*}
	\label{thm:nonreg}
\end{mytheorem}
\begin{myproof}
	See Section~\ref{sec:proofthm1-2} in Appendix.
\end{myproof}

Intuitively, the rate in \eqref{eqn:rate} resembles a near-minimax rate of estimation of high-dimensional anisotropic functions.
{The first part in \eqref{eqn:rate} is the near-minimax risk of the penalty for not knowing the subset $S_0$ \citep{raskutti2011minimax}.} The second part in \eqref{eqn:rate} is incurred by anisotropic regression function estimation.
Although $\lambda$ and $R$ can be a polynomial in $n$ with a suitably small power to satisfy $\epsilon_n\rightarrow 0$, a particularly interesting case is when both are at most $\log^c n$ for some $c>0$. The second term then corresponds to the near-minimax rate of anisotropic function estimation \citep{hoffman2002random}. Whether the rate in \eqref{eqn:rate} is in fact the actual (near) minimax rate remains to be established. The answer to this question is provided  in the following subsection, where we formally derive the minimax lower bound with respect to the $L_2$-risk.

\begin{myremark}
	In isotropic regression using BART, \citet{rockova2020posterior} assumed that the first part of the rate in \eqref{eqn:rate} is dominated by the second part, whereby the resulting rate is simplified such that it only depends on the risk of  function estimation. 
	As this restriction is not required, we keep the rate in the form of \eqref{eqn:rate}.
\end{myremark}

\subsection{Minimax Lower Bound}
\label{sec:minimax}
In Section \ref{sec:nonregcon}, we established the posterior contraction rate of BART under relaxed smoothness assumptions. Although the rate in \eqref{eqn:rate} consists of two logical components (a penalty for variable selection uncertainty and a rate of anisotropic function estimation), it is not guaranteed that the {\em whole rate} is (nearly) minimax optimal. While the minimax rates in high-dimensional {\em isotropic} function estimation were studied exhaustively in \citet{yang2015minimax},
extensions to  (piecewise) {\em anisotropic} functions {\em have not} been obtained in the literature. We fill this gap by deriving a minimax lower bound in our general smoothness setup. These results will certify that the rates  obtained in Section~\ref{sec:nonregcon} are indeed minimax optimal (with respect to the $L_2$-risk) up to a logarithmic factor.

To deploy the conventional minimax theory, we consider the model with random design given by
\begin{align}
	Y_i=f_0(X_i)+\varepsilon_i,\quad X_i\sim Q,\quad\varepsilon_i\sim \text{N}(0,\sigma_0^2),\quad i=1,\dots,n,
	\label{eqn:modelregrd}
\end{align}
where $X_i=(X_{i1},\dots,X_{ip})$, $i=1,\dots,n$, are $p$-dimensional random covariates and $Q$ is a probability measure such that ${\rm supp}(Q)\subseteq [0,1]^p$. We assume (without loss of generality) that $\sigma_0^2$ is fixed to $1$. To obtain a lower bound of the minimax rate, we use the Le Cam equation \citep{birge1993rates,wong1995probability,barron1999risk}.
Now the density $q$ of $Q$ is assumed to satisfy the following assumption under which the $L_2(Q)$-norm is replaced by the $L_2$-norm.

\begin{enumerate}[leftmargin=2.0\parindent,label=\rm(M),start=8]
	\item \label{asm:density} There exist constants $0<\underline q\le\overline q\le\infty$ such that the density $q$ satisfies $\underline q\le \inf_{x} q(x)\le\sup_{x} q(x)\le \overline q$.
\end{enumerate}
We define the $L_2$-minimax risk for any function space $\mathcal F\in \mathcal L_2$ as
\begin{align}
	r_n^2 (\mathcal F) = \inf_{\hat f \in\mathcal B_n} \sup_{f_0\in \mathcal F}\mathbb E_{f_0,Q}\lVert \hat f-f_0 \rVert_{2}^2,
	\label{eqn:minimaxrisk}
\end{align}
where $\mathcal B_n$ is the space of all $\mathcal L_2$-measurable function estimators and $\mathbb E_{f,Q}$ is the expectation operator under the model with $f$ and $Q$.
The Le Cam equation requires suitable upper and lower bounds of the metric entropy of the target function space. We thus define the bounded function space $\overline{\Gamma}^{A_{\bar\alpha},d,p}_{\lambda,M}(\mathfrak X_0)=\{f\in{\Gamma}^{A_{\bar\alpha},d,p}_{\lambda}(\mathfrak X_0) : \lVert f \rVert_\infty\le M\lambda\}$ for any $M>0$.
As our contraction rate is the same for both ${\Gamma}^{A_{\bar\alpha},d,p}_{\lambda}(\mathfrak X_0)$ and ${\Gamma}^{A_{\bar\alpha},d,p}_{\lambda}(\mathfrak X_0)\cap {\mathcal C}([0,1]^p)$, we aim to construct a lower bound of $r_n \big(\overline{\Gamma}^{A_{\bar\alpha},d,p}_{\lambda,M}(\mathfrak X_0)\cap {\mathcal C}([0,1]^p)\big)$ close enough to $\epsilon_n$.

\begin{mytheorem}[Minimax lower bound]
	Consider model \eqref{eqn:modelregrd} for $\sigma_0^2=1$ with Assumption \ref{asm:density}.
	For $d>0$, $\lambda>0$, $R>0$, a partition $\mathfrak X_0=\{\Xi_1,\dots,\Xi_R\}$ of $[0,1]^d$, and a smoothness parameter $A_{\bar\alpha}\in \mathcal A_{\bar\alpha}^{R,d}$ for $\bar\alpha\in(0,1]$ such that
	$\log\mathsf{len}([\Xi_r]_j)\gtrsim -1/\alpha_{rj}$, $1\le r\le R$, $1\le j\le d$,  there exists $M_d>0$ depending only on $d$ such that 
	\begin{align*}
		r_n \big(\overline{\Gamma}^{A_{\bar\alpha},d,p}_{\lambda,M}(\mathfrak X_0)\cap {\mathcal C}([0,1]^p)\big)\gtrsim \sqrt{\frac{1}{n}\log\binom{p}{d}}+M_d\left(\frac{\lambda^{d/\bar\alpha}}{n}\right)^{\bar\alpha/(2\bar\alpha+d)}.
	\end{align*}
	\label{thm:minimax}
\end{mytheorem}

\begin{myproof}
	See Section~\ref{sec:proofthm3} in Appendix.
\end{myproof}

As $M_d$ can be dependent on $d$, the correct interpretation of the result is with a bounded $d$. 
Also, our contraction rate $\epsilon_n$ is derived under the condition $\lVert f_0\rVert_\infty\lesssim \sqrt{\log n}$, and hence we assume that $\lambda\lesssim \sqrt{\log n}$ to match the two spaces.
One can easily verify that the condition $\log\mathsf{len}([\Xi_r]_j)\gtrsim -1/\alpha_{rj}$, $1\le r\le R$, $1\le j\le d$, leads to the restriction $\log R\lesssim d/\bar\alpha$, which removes the term $R$ from our rate $\epsilon_n$ in \eqref{eqn:rate}. Putting the bounds together, $\epsilon_n$ matches the lower bound up to a logarithmic factor.

\subsection{Numerical Study}
In this section, we conduct a numerical study that shows the successful performance of BART  with a variety of multivariate functions. For competitors we consider Gaussian process (GP) prior regression, gradient boosting (GB), random forest (RF), and neural network (NN) models with the rectified linear unit (ReLU) activation function. GP prior regression is widely exploited for multiple nonparametric regression and ensures theoretical optimality for smooth functions \citep{van2008rates}.
GB is expected to work similarly to BART.
RF is expected to satisfactorily detect discontinuous boundaries along the coordinates, as it is based on the additive tree ensembles.
We know that NN models adapt well to complicated function classes with the guaranteed optimal properties \citep[e.g.,][]{petersen2018optimal,imaizumi2019deep,schmidt2020nonparametric,hayakawa2020minimax}. Our numerical study shows that BART outperforms these competitors in adapting to complicated smoothness structures.

\begin{figure}[t!]
	\centering
	\begin{subfigure}[b]{2.2in}
		\includegraphics[width=\linewidth]{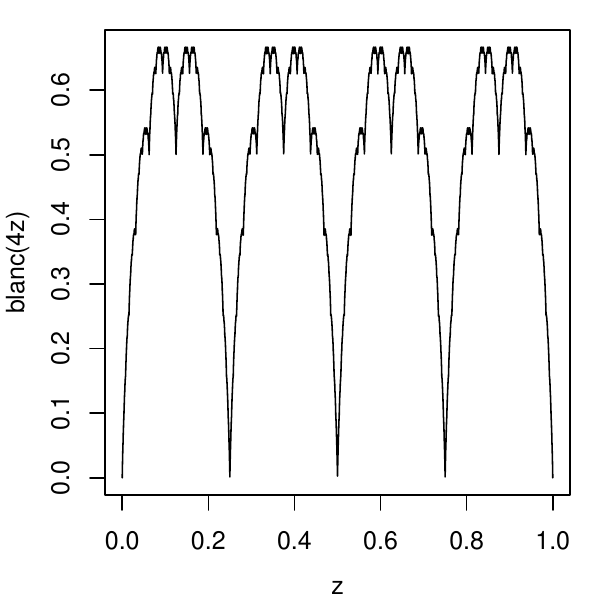}	
		\caption{Blancmange function $\mathsf{blanc}(4z)$}
	\end{subfigure}
	\begin{subfigure}[b]{2.2in}
	\includegraphics[width=\linewidth]{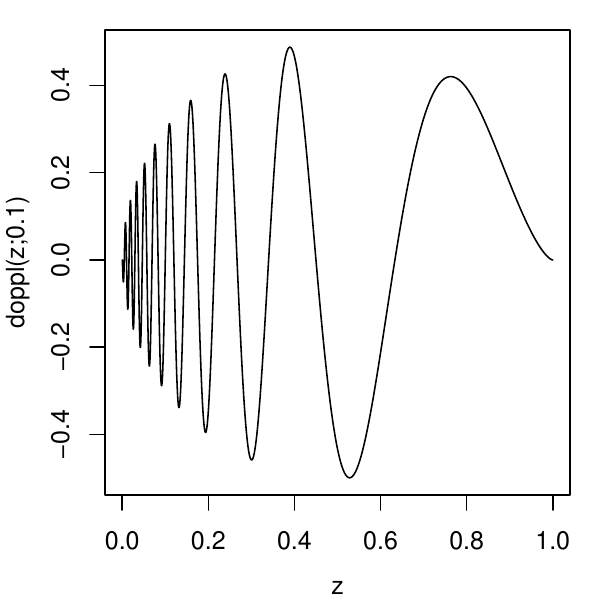}	
	\caption{Doppler function $\mathsf{doppl}(z;0.1)$}
\end{subfigure}
	\caption{Functions that exhibit non-Lipschitz continuity and spatially varying smoothness.}
	\label{fig:aniso}
\end{figure}

Our synthetic datasets are generated from model \eqref{eqn:modelreg} with a few different functions $f_0:[0,1]^p\rightarrow \mathbb R$. To specify the simulation setups, we first introduce the following functions that maps $[0,1]^p$ to $\mathbb R$:
\begin{align*}
	\mathsf{base}_p: (x_1,\dots,x_p) &\mapsto \sin\bigg(\frac{10}{\sqrt{p}}\bigg\{\sum_{j=1}^p (x_j-0.5)^2-\frac{p}{12}\bigg\}\bigg),\\
	 \mathsf{discont1}: (x_1,\dots,x_p) &\mapsto \mathbbm 1(x_1\le 0.5, x_2>0.5)+\mathbbm 1(x_1>0.5, x_2\le0.5),\\
	 \mathsf{discont2}_p:(x_1,\dots,x_p) &\mapsto\mathbbm 1\bigg(\sum_{j=1}^p(x_j-0.5)\le0, \sum_{j=1}^p (-1)^j(x_j-0.5)>0\bigg) \\
	 &\quad+  \mathbbm 1\bigg(\sum_{j=1}^p(x_j-0.5)> 0, \sum_{j=1}^p (-1)^j(x_j-0.5)\le 0\bigg).
\end{align*}
The function $\mathsf{base}_p$ is viewed as having an isotropic smoothness and is used as the base component for $f_0$.\footnote{The argument of the sine function is chosen so that it is centered at zero and has a reasonable scale for every $p$, allowing the period of the sine function to be roughly maintained with $p$.
In particular, if $X_j$ has a uniform distribution on $[0,1]$ independently, one can easily see that $(10/\sqrt{p})\{\sum_{j=1}^p (X_j-0.5)^2-p/12\}$ weakly converges to $\text{N}(0,5/9)$ as $p\rightarrow \infty$.} The functions $\mathsf{discont1}$ and $\mathsf{discont2}_p$ render discontinuous jumps along hyperplanes in different directions.
To account for non-Lipschitz continuity and spatially varying smoothness, we also define the blancmange function and the Doppler function as,
\begin{align*}
	 \mathsf{blanc}(z)&=\sum_{k=0}^\infty \frac{|2^k z-\lfloor 2^k z+0.5 \rfloor|}{2^k},\quad z\in[0,1], \\
\mathsf{doppl}(z;a)&=\sqrt{z(1-z)}\sin\!\left(\frac{2\pi(1+a)}{z+a}\right),\quad z\in[0,1],
\end{align*}
which are illustrated in Figure~\ref{fig:aniso}.

\begin{figure}[t!]
	\centering
	\begin{subfigure}[b]{1.95in}
		\includegraphics[width=1.95in]{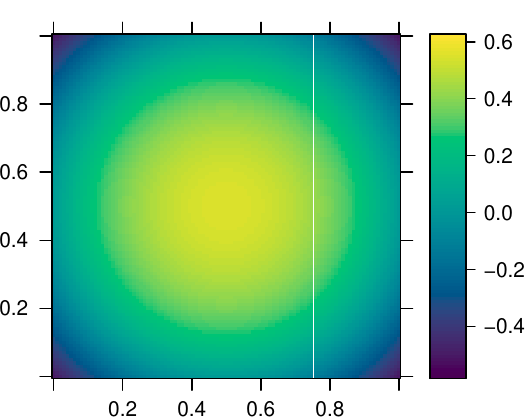}
		\caption{$f_0^{(1)}$}
	\end{subfigure}
	\begin{subfigure}[b]{1.95in}
	\includegraphics[width=1.95in]{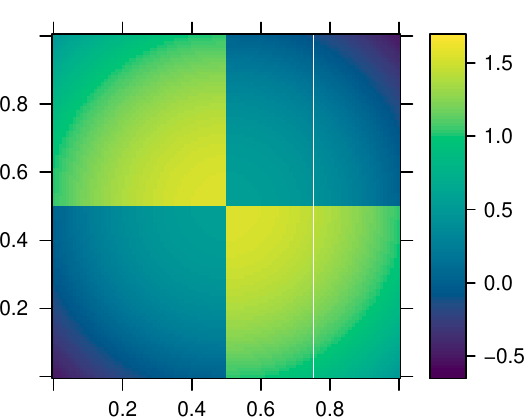}
	\caption{$f_0^{(2)}$}
\end{subfigure}
	\begin{subfigure}[b]{1.95in}
	\includegraphics[width=1.95in]{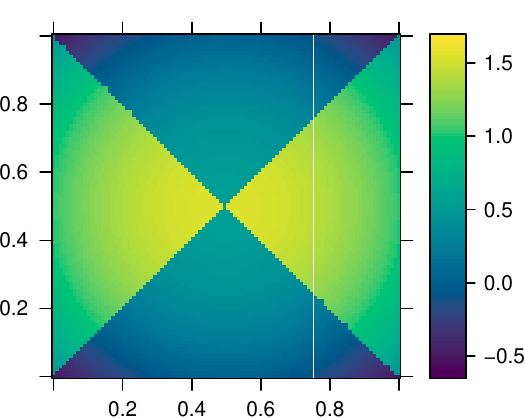}
	\caption{$f_0^{(3)}$}
\end{subfigure}
\\ \vspace*{10pt}
	\begin{subfigure}[b]{1.95in}
	\includegraphics[width=1.95in]{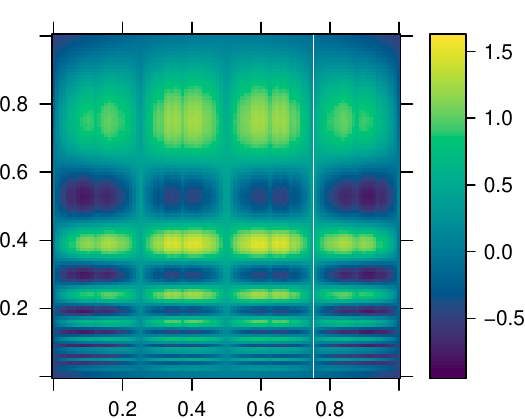}
	\caption{$f_0^{(4)}$}
\end{subfigure}
	\begin{subfigure}[b]{1.95in}
	\includegraphics[width=1.95in]{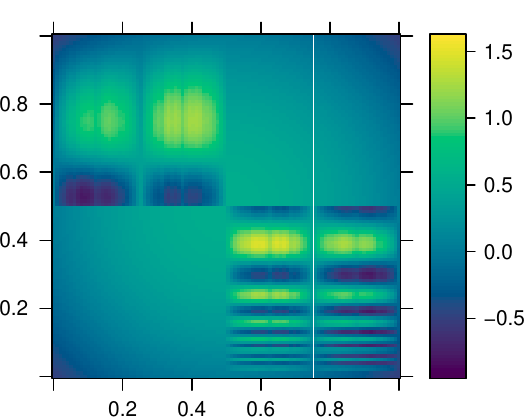}
	\caption{$f_0^{(5)}$}
\end{subfigure}
	\begin{subfigure}[b]{1.95in}
	\includegraphics[width=1.95in]{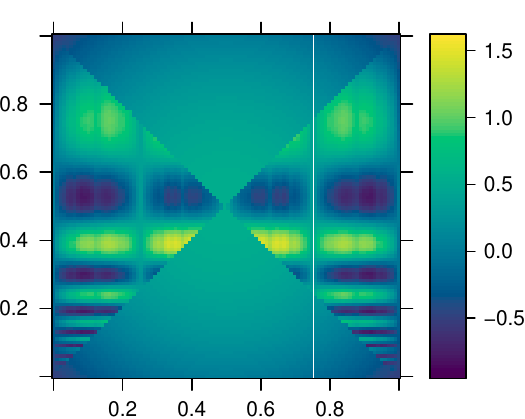}
	\caption{$f_0^{(6)}$}
\end{subfigure}
	\caption{Level plots of the true functions $f_0^{(k)}$, $k=1,\dots,6$, for $p=2$.}
	\label{fig:true}
\end{figure}

Using the above functions, we describe six simulation scenarios. Specifically, Scenario~$k$ is defined by model \eqref{eqn:modelreg} with $f_0=f_0^{(k)}$, $k=1,\dots,6$, where the true functions $f_0^{(k)}:[0,1]^p\rightarrow \mathbb R$ are defined as
\begin{align*}
f_0^{(1)}:(x_1,\dots,x_p) &\mapsto  \mathsf{base}_p(x_1,\dots,x_p), \\
f_0^{(2)}:(x_1,\dots,x_p) &\mapsto \mathsf{base}_p(x_1,\dots,x_p) +  \mathsf{discont1}(x_1,\dots,x_p), \\
f_0^{(3)}:(x_1,\dots,x_p) &\mapsto  \mathsf{base}_p(x_1,\dots,x_p)  +  \mathsf{discont2}_p(x_1,\dots,x_p), \\
f_0^{(4)}:(x_1,\dots,x_p) &\mapsto \mathsf{base}_p(x_1,\dots,x_p)+3\mathsf{blanc}(4x_1)\mathsf{doppl}(x_2;0.1), \\
f_0^{(5)}:(x_1,\dots,x_p) &\mapsto \mathsf{base}_p(x_1,\dots,x_p)+3\mathsf{blanc}(4x_1)\mathsf{doppl}(x_2;0.1)\mathsf{discont1}(x_1,\dots,x_p), \\
f_0^{(6)}:(x_1,\dots,x_p) &\mapsto \mathsf{base}_p(x_1,\dots,x_p)+3\mathsf{blanc}(4x_1)\mathsf{doppl}(x_2;0.1)\mathsf{discont2}_p(x_1,\dots,x_p).
\end{align*}
The functions $f_0^{(1)}$ and $f_0^{(4)}$ represent globally isotropic and anisotropic functions, respectively. The other functions produce discontinuous jumps that are either parallel or oblique to the coordinate system.
Specifically, $f_0^{(2)}$ and $f_0^{(5)}$ are regarded as piecewise isotropic and anisotropic functions, respectively, as defined in Definition~\ref{def:genhol}. The remaining functions $f_0^{(3)}$ and $f_0^{(6)}$ are similarly piecewise isotropic and anisotropic, but they differ from Definition~\ref{def:genhol} in that the jumps are not parallel to the coordinates.
The two-dimensional case of each $f_0^{(k)}$ is visualized in Figure~\ref{fig:true}.

We generate the synthetic datasets under Scenarios~1--6. For each scenario, we consider two sample sizes $n\in\{1000,5000\}$ and five dimension values $p\in\{2,5,10,20,50\}$, while fixing $\sigma_0^2=0.5^2$ for reasonable signal to noise ratios. Therefore, each scenario has 10 synthetic datasets generated with all possible combinations of $n$ and $p$. For given predictor variables $X_i$ generated uniformly on $[0,1]^p$, the response variable $Y_i$ is generated from model \eqref{eqn:modelreg}, $i=1,\dots,n$.

All datasets are fitted by BART and the other competitors. 
For a fair comparison to the other methods, we do not use the Dirichlet sparse prior in \eqref{eqn:dir} for BART. Instead, we assign a uniform prior that corresponds to the Dirichlet prior with concentration parameter $1$, with a priori assumption that all predictor variables contribute equally to the observations. 
We fit BART with 200 trees using the prior that splits a node at depth $\ell$ with probability $\alpha (\ell+1)^{-\beta}$ for $\alpha\in(0,1)$ and $\beta\in[0,\infty]$, the original construction by \citet{chipman2010bart}, which is implemented in the R package \texttt{BART}. However, as our theory resorts to the exponentially decaying prior for splits as mentioned in Section~\ref{sec:treeprior}, we also consider BART with the prior that splits a node at depth $\ell$ with probability $\nu^{\ell+1}$ for $\nu\in(0,1/2)$. We choose $\alpha=0.3$, $\beta=2$, and $\nu=0.3$ to make the two priors roughly similar for small $\ell$. We will see that the two priors exhibit similar empirical behavior.
For GP prior regression, the squared exponential covariance kernel $k(x,x')=\tau^2\exp(\lVert x-x'\rVert^2/l^2)$ is employed with half normal priors $\tau\sim \text{N}_+(0,1)$ and $l\sim \text{N}_+(0,1)$. Optimizing other parameters in the posterior distribution, the posterior mode of $f_0$ is obtained in a closed-form expression (we also tried other informative priors for $\tau$ and $l$ and observed no significant difference).
GB is trained by the \texttt{gbm} package with trees of five splits and the number of trees determined via cross validation (CV).
RF is fitted by the \texttt{randomForest} package with 200 trees and the maximal node size $5$ or $50$ for each tree.
The NN models are trained by TensorFlow with the Keras interface. We consider two NN models with two and four hidden layers with $(64,32)$ and $(256,128,64,32)$ hidden units. All hidden units take the ReLU activation function with the dropout of rate $0.3$ for regularization.
The description of the methods is summarized in Table~\ref{tab:method}.

\begin{table}[!t]
	\caption{The description of the methods for simulation.}
	\centering
			\begin{tabular}{ll}
\toprule
				Method &  Description \\
\midrule
				BART1  & BART with 200 trees \\
				& Node at depth $\ell$ is split with prior probability $\alpha (\ell+1)^{-\beta}$, $\alpha=0.3$, $\beta=2$\\
				BART2   & BART with 200 trees \\
				& Node at depth $\ell$ is split with prior probability $\nu^{\ell+1}$, $\nu=0.3$  \\
				GP  &  GP prior regression with the squared exponential covariance kernel \\
				GB & GB with trees of five splits and the number of trees determined via CV\\
				RF1 &  RF of 200 trees with maximal node size 5 for each tree  \\
		 		RF2 &  RF of 200 trees with maximal node size 50 for each tree  \\
		 		NN1 &  NN model with two hidden layers and $(64,32)$ hidden units \\
		 		NN2 &  NN model with four hidden layers and $(256,128,64,32)$ hidden units	\\
\bottomrule
			\end{tabular}
			\label{tab:method}  
\end{table}

\begin{figure}[t!]
	\centering
	\includegraphics[width=6in]{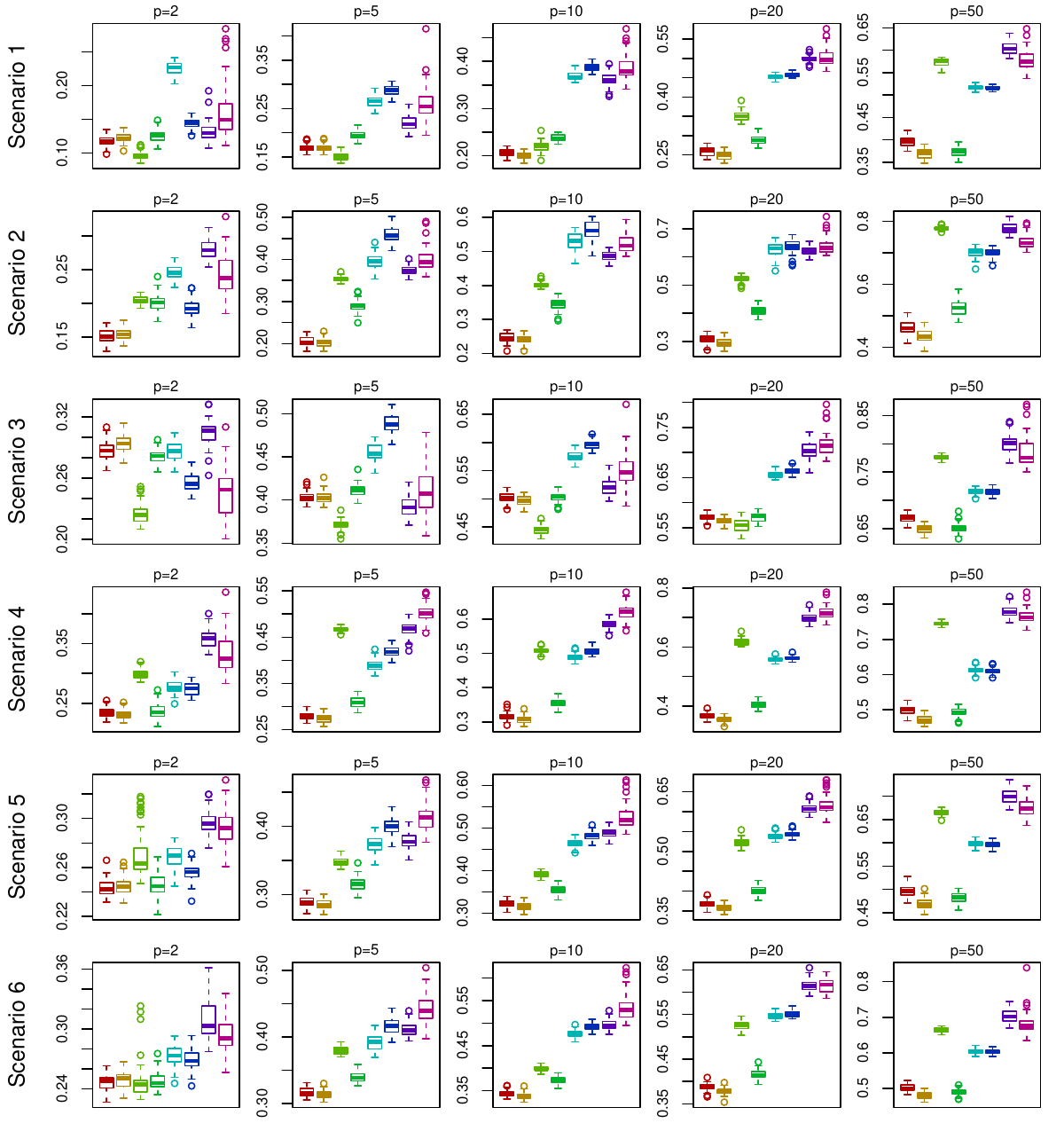}\\
	\includegraphics[width=4.5in]{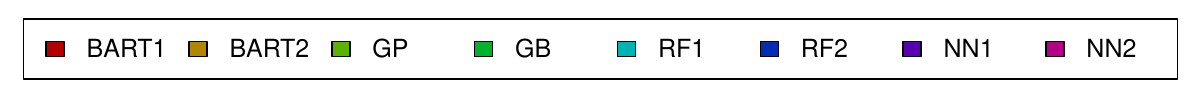}	
	\caption{RMSPEs obtained from 50 replicated datasets of size $n=1000$.}
	\label{fig:sim1}
\end{figure}

\begin{figure}[t!]
	\centering
		\includegraphics[width=6in]{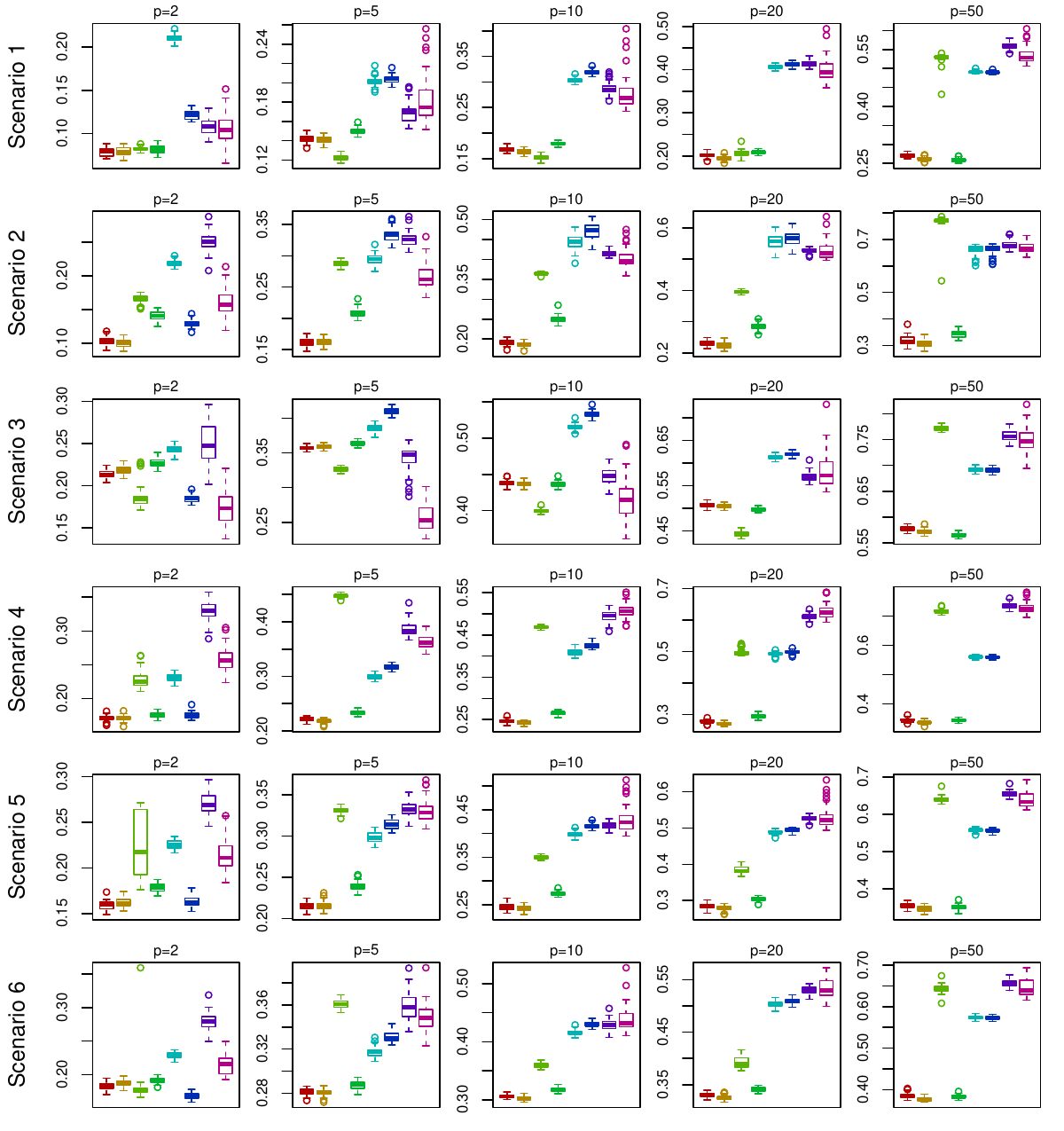}\\
	\includegraphics[width=4.5in]{sim-legend.pdf}	
	\caption{RMSPEs obtained from 50 replicated datasets of size $n=5000$.}
	\label{fig:sim2}
\end{figure}

Figures~\ref{fig:sim1} and \ref{fig:sim2} show the root mean squared prediction error (RMSPE) obtained by the methods described in Table~\ref{tab:method}. The RMSPEs are estimated by randomly drawn out-of-samples. For Scenario~1 with the global isotropic function $f_0^{(1)}$, BART, GP regression, and GB perform similarly well in relatively lower dimensions ($p=2,5,10$), but the performance of GP degrades as $p$ increases. For Scenario~2 with the piecewise isotropic function $f_0^{(2)}$, BART clearly outperforms the other methods as expected. 
Interestingly, GB performs substantially worse than BART in this situation, implying that BART detects discontinuous jumps along the coordinates better.
RF falls behind BART and GB although it is also based on binary tree ensembles. For Scenario~3 with $f_0^{(3)}$, GP and NN perform better than BART and GB in lower dimensions; this makes sense given that BART cannot detect such discontinuous jumps efficiently using the coordinate parallel splitting rule. However, the performance of GP and NN deteriorates as $p$ increases, and BART and GB beat the competition in higher dimensions ($p=50$). The interpretation of the results is similar for the remaining scenarios. The major difference is that BART produces the best prediction error in almost all cases of Scenarios~4--6. Given that BART is designed to capture local anisotropy very effectively, this finding appears to be a natural consequence. Overall, GB performs slightly worse than BART.
As well as the setups used in our simulation, we also tested many other tuning parameter setups and network structures for GB, RF, and NN, but found no clear improvement.

Based on Figures~\ref{fig:sim1} and \ref{fig:sim2}, we can also compare the performance of the two BART priors.
BART with the polynomially decaying prior (the original BART prior by \citet{chipman2010bart}) works slightly better in lower dimensions ($p=2,5$), whereas the exponentially decaying prior is marginally preferred in higher dimensions ($p=10,20,50$). However, because the difference is not significant, we conclude that there are no substantial differences in empirical behavior between the two BART priors.

\section{Further Applications}
\label{sec:furapp}
Section~\ref{sec:nonreg} establishes the posterior contraction rate of BART for the nonparametric regression model and justifies its near-minimax optimality. As our approximation theory only requires conditions on a split-net, the results can be extended to statistical models beyond nonparametric regression with fixed design. In this section, we consider other applications such as nonparametric regression with random design, density estimation, and nonparametric binary classification. Moreover, as the technical results in Section~\ref{sec:nonreg} hold even with the single tree model ($T=1$), one can find no theoretical advantages of BART over Bayesian CART. 
A theoretical advantage of BART can be recognized if the true function has an additive structure \citep{linero2018bayesian-jrssb,rockova2020posterior}. Such an extension is also considered in this section.

\subsection{Nonparametric Regression with Random Design}
\label{sec:nonran}

Theorem~\ref{thm:nonreg} quantifies the posterior contraction rate of nonparametric regression with fixed design where the predictor variables are not random variables. Now we consider a random design regression in \eqref{eqn:modelregrd} in which the model is treated as independent and identically distributed.
We establish the posterior contraction rate of BART for the random design model in \eqref{eqn:modelregrd}. The main advantage of considering random design is that it provides the $L_{2,Q}$-contraction rate without empirical process theory, where $Q$ is a probability measure for $X_i$, whereas fixed design essentially provides the contraction rate with respect to the empirical $L_2$-norm as in Section~\ref{sec:nonreg}.
The random design assumption is also often necessary in certain statistical models, for example, in measurement error models \citep{tuo2015efficient} or causal inference models \citep{hahn2020bayesian,ray2020semiparametric}. Note that fixed design points in Section~\ref{sec:fixeddesign} cannot be used for a split-net, as the procedure is not truly Bayesian if the prior is dependent on the data ($X_i$ is now considered a part of the observation.). Instead, a regular grid in Section~\ref{sec:grid} can be useful for this framework.

We consider model \eqref{eqn:modelregrd} for $Q$ a probability measure that satisfies ${\rm supp}(Q)\subseteq [0,1]^p$ with a bounded density.
Unlike model \eqref{eqn:modelreg}, model \eqref{eqn:modelregrd} is independent and identically distributed.
The well-known fact that exponentially powerful tests exist with respect to the Hellinger metric $\rho_{\rm H}(\cdot,\cdot)$ allows one to establish the contraction rate for the corresponding metric \citep{ghosal2000convergence}.
However, in normal models, the Hellinger distance is matched to the $L_2$-type metric only when $\lVert f\rVert_\infty$ and $|\log\sigma^2|$ are bounded in the entire parameter space, not only for the true values \citep[e.g.,][]{xie2018adaptive}.
Unlike in Theorem~\ref{thm:nonreg}, this restriction requires that $f_0$ be uniformly bounded and a prior be appropriately truncated.
Note also that we need a good approximation error with respect to the integrated $L_2$-norm.
We summarize the required modifications of \ref{asm:fsup}, \ref{asm:spnet}, \ref{pri:normal}, and \ref{pri:invgam}.
\begin{enumerate}[leftmargin=2.0\parindent,label=\rm(A\arabic*)]
	\item[\mylabel{asm:fsup2}{\rm(A3{$^\ast$})}] The true function $f_0$ satisfies $\lVert f_0 \rVert_\infty\le C_0^\ast$ for some sufficiently large $C_0^\ast>0$.
	\item[\mylabel{asm:spnet2}{\rm(A6{$^\ast$})}] The split-net $\mathcal Z$  is suitably dense and regular to construct a $\mathcal Z$-tree partition $\widehat {\mathcal T}$ such that there exists $\hat f_0\in\mathcal F_{\widehat {\mathcal T}}$ satisfying $\lVert f_0-\hat f_0 \rVert_2\lesssim \bar\epsilon_n$ by Theorem~\ref{thm:approx2}.
	\item[\mylabel{pri:truncb}{\rm(P2{$^\ast$})}] A prior on the compact support $[-\overline C_1,\overline C_1]$ is assigned to the step-heights $B$ for some $\overline C_1>C_0^\ast$.
	\item[\mylabel{pri:trucsig}{\rm(P3{$^\ast$})}] A prior on the compact support $[\overline C_2^{-1},\overline C_2]$ is assigned to $\sigma^2$ for some $\overline C_2>C_0$.
\end{enumerate}
Assumption \ref{asm:spnet2} requires good approximability with respect to the $L_2$-norm.
Owing to \ref{lmm0:Leb1-1} of Corollary~\ref{cor:approx} and Lemma~\ref{lmm:grid2}, a regular grid in Section~\ref{sec:grid} can be useful to meet this requirement (see Remark~\ref{rmk:practical.c}).
We wrap up this section with a theorem that formalizes the posterior contraction of BART for model~\eqref{eqn:modelregrd}. 

\begin{mytheorem}[Nonparametric regression, random design]
	Consider model \eqref{eqn:modelregrd} with 
	Assumptions \ref{asm:truef}, \ref{asm:dp}, \ref{asm:fsup2}, \ref{asm:sig}, \ref{asm:spnet0}, \ref{asm:spnet2}, and \ref{asm:spnet1}, and the prior assigned through \ref{pri:tree}, \ref{pri:truncb}, and \ref{pri:trucsig}.
	Then, there exists a constant $M>0$ such that for $\epsilon_n$ in \eqref{eqn:rate},
	\begin{align*}
		\mathbb E_0 \Pi\Big\{(f,\sigma^2):\lVert f-f_0\rVert_{2,Q}+|\sigma^2-\sigma_0^2|>M \epsilon_n \,\big|\, (X_1,Y_1),\dots,(X_n,Y_n)\Big\}\rightarrow 0.
	\end{align*}
	\label{thm:nonregrd}
\end{mytheorem}
\begin{myproof}
	See Section~\ref{sec:proofthm4-5} in Appendix.
\end{myproof}

\subsection{Density Estimation}
\label{sec:density}
In addition to classical nonparametric regression, density estimation is an interesting branch of nonparametric inference. There exist a few studies employing the Bayesian tree ensembles for density regression \citep{orlandi2021density,li2022adaptive}. Here we consider a more traditional density estimation problem. With the Bayesian tree ensembles, we only provide a theoretical flavor for density estimation rather than practical implementation. It may be difficult to develop an efficient algorithm for the setup considered here.

For some probability measure $P$ that satisfies ${\rm supp}(P)\subseteq [0,1]^p$,
suppose $n$ independent observations $X_i$, $i=1,\dots,n$, are drawn from $P$, i.e.,
\begin{align}
	X_i \sim P,\quad i=1,\dots, n.
	\label{eqn:density}
\end{align}
Assume that $P$ is absolutely continuous with respect to the Lebesgue measure $\mu$ with the true density $p_0$.
We assign a prior on $p_f$ indexed by $f$ such that
$p_f=e^f/\int_{[0,1]^p} e^f d\mu$ with $f$ assigned the forest priors in Section~\ref{sec:prior}.
We write $f_0=\log p_0$ while assuming \ref{asm:truef}--\ref{asm:fsup}.
That is, our $d$-sparsity for density estimation implies that the remaining $p-d$ variables are independent and uniformly distributed on $[0,1]^{p-d}$. This sparsity setup is useful in high dimensions because the density cannot be estimated effectively without a stronger assumption for $p>n$, such as isotropy. A similar sparsity structure was also imposed in \citet{liu2007sparse} for high-dimensional density estimation.	
We leverage the existence of an exponentially powerful test for the Hellinger metric $\rho_{\rm H}(\cdot,\cdot)$.
Owing to the relationship between Hellinger balls and $L_\infty$ balls in density estimation with the exponential link, we need an approximation result with respect to the $L_\infty$-norm. This is obtained by \ref{lmm0:sup2} of Theorem~\ref{thm:approx2} with the continuity restriction on the true function. As \ref{lmm0:sup2-1} of Corollary~\ref{cor:approx} and Lemma~\ref{lmm:grid2} show, a regular grid in Section~\ref{sec:grid} is useful to obtain the $L_\infty$-approximation (see Remark~\ref{rmk:practical.c}).
We make the following assumptions to satisfy this requirement.

\begin{enumerate}[leftmargin=2.0\parindent,label=\rm(A\arabic*)]
	\item [\mylabel{asm:truef2}{\rm(A1{$^{\ddag}$})}] For $d>0$, $\lambda>0$, $R>0$, $\mathfrak X_0=\{\Xi_1,\dots,\Xi_R\}$, and $A_{\bar\alpha}\in\mathcal A_{\bar\alpha}^{R,d}$ with $\bar\alpha\in(0,1]$, the true function satisfies  $f_0\in\Gamma^{A_{\bar\alpha},d,p}_{\lambda}(\mathfrak X_0)\cap \mathcal C([0,1]^p)$.
	\item[\mylabel{asm:spnet3}{\rm(A6{$^{\ddag}$})}] The split-net $\mathcal Z$  is suitably dense and regular to construct a $\mathcal Z$-tree partition $\widehat {\mathcal T}$ such that there exists $\hat f_0\in\mathcal F_{\widehat {\mathcal T}}$ satisfying $\lVert f_0-\hat f_0 \rVert_{\infty}\lesssim \bar\epsilon_n$ by Theorem~\ref{thm:approx2}.
\end{enumerate}

We assign the tree prior with Dirichlet sparsity and a normal prior on the step-heights.
Under suitable assumptions, the following theorem provides the posterior contraction rate for $p_f$ with respect to the Hellinger distance.

\begin{mytheorem}[Density estimation]	
	Consider model \eqref{eqn:density} with Assumptions \ref{asm:truef2}, \ref{asm:dp}--\ref{asm:fsup}, \ref{asm:spnet0}, \ref{asm:spnet3}, and \ref{asm:spnet1}, and the prior assigned through \ref{pri:tree}--\ref{pri:normal}.
	Then, there exists a constant $M>0$ such that for $\epsilon_n$ in \eqref{eqn:rate},
	\begin{align*}
		\mathbb E_0 \Pi\Big\{f:\rho_{\rm H}(p_f,p_0)>M \epsilon_n \,\big|\, X_1,\dots,X_n\Big\}\rightarrow 0.
	\end{align*}
	\label{thm:density}
\end{mytheorem} 

\begin{myproof}
	See Section~\ref{sec:proofthm4-5} in Appendix.
\end{myproof}

As mentioned in Section~\ref{sec:nonregcon}, the normal prior in \ref{pri:normal} is not necessary and a heavy-tailed prior can relax the assumption on $\lVert f \rVert_\infty$. As normal priors are not conjugate to the model likelihood in the density estimation example, there is no clear benefit of adopting \ref{pri:normal} anymore.
This is also the case in the example of binary classification given in the next subsection.
Nevertheless, we employ \ref{pri:normal} for the sake of simplicity. 

\begin{myremark}
	As previously stated, the practical implementation of the density estimation problem here is not as straightforward as the Gaussian regression case. 
	We do not believe that there is a highly efficient algorithm for density estimation with the type of Bayesian forest considered here. One possible option is employing the idea of  reversible jump moves \citep{green1995reversible}, as in \citet{linero2022generalized} for generalized BART for exponential family models. 
\end{myremark}

\subsection{Nonparametric Binary Classification}
\label{sec:classification}

Nonparametric classification is useful for modeling categorical response variables. In the original work by \citet{chipman2010bart}, BART for Gaussian regression was readily adapted to probit regression using the latent variable expression \citep{albert1993bayesian}. Later, \citet{kindo2016multinomial} devised a BART algorithm for multi-category response variables using multinomial probit models. Although the probit models are particularly simple to implement, we consider nonparametric binary classification with the logistic link function to make use of the classical theory \citep{van2008rates}. The computation is still straightforward owing to the latent variable expression with a P\'olya-gamma distribution \citep{polson2013bayesian}.

For a binary response $Y_i\in\{0,1\}$ and a random covariate $X_i\in\mathbb R^p$, assume that we have $n$ independent observations $(X_1,Y_1),\dots(X_n,Y_n)$ from  
the binary classification model,
\begin{align}
	\mathbb E_0[\mathbbm  1(Y_i=1)|X_i=x] = \varphi_0(x), \quad X_i \sim Q,\quad i=1,\dots, n,
	\label{eqn:binary}
\end{align}
for some $\varphi_0: [0,1]^p\rightarrow [0,1]$ and some probability measure $Q$ such that ${\rm supp}(Q)\subseteq [0,1]^p$ with a bounded density. We thus consider a binary classification problem with random design.
We parameterize the probability function using the logistic link function $H:\mathbb R\rightarrow [0,1]$ such that
$\varphi_f=H(f)$ for $f$ on which the forest priors in Section~\ref{sec:prior} are assigned.
For true function $\varphi_0$,
we write $f_0=H^{-1}(\varphi_0)$ while assuming \ref{asm:truef}--\ref{asm:fsup} as in the density estimation problem.
The proof shows that the Hellinger metric is bounded by the $L_2(Q)$-distance in this example, and hence \ref{asm:spnet2} is assumed. Similar to Section~\ref{sec:nonran}, fixed design points are not available for a split-net, but a regular grid in Section~\ref{sec:grid} can be useful.
The following theorem formalizes the posterior contraction rate with respect to the $L_2(Q)$-distance.

\begin{mytheorem}[Binary classification]
	Consider model \eqref{eqn:binary} with Assumptions \ref{asm:truef}--\ref{asm:fsup}, \ref{asm:spnet0}, \ref{asm:spnet2}, and \ref{asm:spnet1}, and the prior assigned through \ref{pri:tree}--\ref{pri:normal}.
	Then, there exists a constant $M>0$ such that for $\epsilon_n$ in \eqref{eqn:rate},
	\begin{align*}
		\mathbb E_0 \Pi\Big\{f:\lVert H(f)- H (f_0)\rVert_{2,Q}>M \epsilon_n \,\big|\, (X_1,Y_1),\dots,(X_n,Y_n)\Big\}\rightarrow 0.
	\end{align*}
	\label{thm:binary}
\end{mytheorem}
\begin{myproof}
	See Section~\ref{sec:proofthm4-5} in Appendix.
\end{myproof}

\subsection{Additive Nonparametric Regression}
\label{sec:addnonpreg}

Thus far we have considered statistical models with the true function $f_0$ that belongs to the piecewise heterogeneous anisotropic H\"older space with sparsity. As Theorems~\ref{thm:nonreg}-\ref{thm:binary} hold even with the single tree model ($T=1$), the empirical success of BART is not well explained by the previous examples, although the empirical performance of BART should be attributed to its fast mixing to some extent.
However, \citet{linero2018bayesian-jrssb} and \citet{rockova2020posterior} observed that BART optimally adapts to a larger class of additive functions which single tree models do not adapt to. In this section, we consider additive nonparametric regression to show theoretical advantages of BART over Bayesian CART.

We consider the nonparametric regression model with fixed design in \eqref{eqn:modelreg}, but the true function $f_0$ is assumed to have an additive structure with $T_0$ components, $f_0=\sum_{t=1}^{T_0} f_{0t}$, where each $f_{0t}$ belongs to the piecewise heterogeneous anisotropic H\"older space with sparsity. We also need suitable conditions on a split-net $\mathcal Z$ such that the approximation theory works for every additive component. We thus make the following modifications of the conditions used in Section~\ref{sec:nonregcon}. In what follows, the subscript or superscript $t$ stands for additive component-specific extensions of the model elements used in Section~\ref{sec:nonregcon}.

\begin{enumerate}[leftmargin=2.0\parindent,label=\rm(A\arabic*)]
	\item [\mylabel{asm:truefadd}{\rm(A1{$^{\mathsection}$})}]  For $d_t>0$, $\lambda_t>0$, $R_t>0$, ${\mathfrak X_0}_t=\{\Xi_{t1},\dots,\Xi_{tR}\}$, and $A_{t,\bar\alpha_t}\in\mathcal A_{\bar\alpha_t}^{R_t,d_t}$ with $\bar\alpha_t\in(0,1]$, $t=1,\dots ,T_0$, the true function satisfies $f_0=\sum_{t=1}^{T_0} f_{0t}$ for $f_{0t}\in\Gamma^{A_{t,\bar\alpha_t},d_t,p}_{\lambda_t}({\mathfrak X_0}_t)$ or  $f_{0t}\in\Gamma^{A_{t,\bar\alpha_t},d_t,p}_{\lambda_t}({\mathfrak X_0}_t)\cap \mathcal C([0,1]^p)$.
	\item [\mylabel{asm:dpadd}{\rm(A2{$^{\mathsection}$})}] It is assumed that $d_t$, $p_t$, $\lambda_t$, $R_t$, and $\bar\alpha_t$ satisfy $\epsilon_{t,n}\ll 1$, where 	$\epsilon_{t,n}=\sqrt{(d_t\log p)/{n}}+(\lambda_t d_t)^{d_t/(2\bar\alpha_t+d_t)}\left(({R_t\log n})/{n}\right)^{\bar\alpha_t/(2\bar\alpha_t+d_t)}$.
	\item [\mylabel{asm:spnetadd}{\rm(A6{$^{\mathsection}$})}] The split-net $\mathcal Z$  is suitably dense and regular to construct a $\mathcal Z$-tree partition $\widehat {\mathcal T}^t$ such that for $\bar\epsilon_{t,n}=(\lambda_t d_t)^{d_t/(2\bar\alpha_t+d_t)}\left(({R_t\log n})/{n}\right)^{\bar\alpha_t/(2\bar\alpha_t+d_t)}$, there exists $\hat f_{0t}\in\mathcal F_{\widehat {\mathcal T}^t}$ satisfying $\lVert f_{0t}-\hat f_{0t} \rVert_n\lesssim \bar\epsilon_{t,n}$ by Theorem~\ref{thm:approx2}, $t=1,\dots ,T_0$. 
	\item [\mylabel{asm:spnet1add}{\rm(A7{$^{\mathsection}$})}] The $\mathcal Z$-tree partition $\mathcal T_t^\ast=\{\Omega_{t1}^\ast,\dots,\Omega_{tR}^\ast\}$ approximating ${\mathfrak X_0}_t^\ast$ satisfies $\max_r\mathsf{dep}(\Omega_{tr}^\ast)\lesssim \log n$, $t=1,\dots ,T_0$.
\end{enumerate}
These simply mean that the assumptions in Section~\ref{sec:nonregcon} hold for every additive component $f_{0t}$. 
It is worth noting that we do not need to modify the prior distribution for additive regression, which makes BART very appealing in that the procedure truly adapts to the unknown true function. This is owing to the use of the Dirichlet prior in \eqref{eqn:dir}; the spike-and-slab prior does not yield such a nice property \citep{rockova2020posterior}. The next theorem provides the posterior contraction rate for the additive regression model.

\begin{mytheorem}[Additive nonparametric regression]
	Consider model \eqref{eqn:modelreg} with Assumptions \ref{asm:truefadd}--\ref{asm:dpadd}, \ref{asm:fsup}--\ref{asm:spnet0}, and \ref{asm:spnetadd}--\ref{asm:spnet1add} and the prior assigned through \ref{pri:tree}--\ref{pri:invgam}.
	If $T_0\le T$, there exists a constant $M>0$ such that for $\epsilon_n^\ast=\sqrt{\sum_{t=1}^{T_0}\epsilon_{t,n}^2}$,
	\begin{align*}
		\mathbb E_0 \Pi\Big\{(f,\sigma^2):\lVert f-f_0\rVert_n+|\sigma^2-\sigma_0^2|>M \epsilon_n^\ast \,\big|\, Y_1,\dots,Y_n\Big\}\rightarrow 0.
	\end{align*}
	\label{thm:addreg}
\end{mytheorem}

\begin{myproof}
	See Section~\ref{sec:proofthm4-5} in Appendix.
\end{myproof}

Theorem~\ref{thm:addreg} shows that the posterior contraction rate for additive regression is the sum of the rates for the additive components. If the function space is reduced to a high-dimensional isotropic class, then our rate $\epsilon_n^\ast$ matches the minimax rate for high-dimensional additive regression \citep{yang2015minimax}. We believe that $\epsilon_n^\ast$ is indeed near-minimax optimal, which can be formally justified by combining the proof technique of our Theorem~\ref{thm:minimax} and the tools for additive scenarios developed in \citet{yang2015minimax}. Considering the length of the paper, we do not pursue this direction in this study.

\section{Discussion}
\label{sec:disc}
In this study, we enlarged the scope of theoretical understanding of Bayesian forests in the context of function estimation by considering relaxed smoothness assumptions. We introduced a new class of piecewise anisotropic sparse functions, which form a blend of anisotropy and spatial inhomogeneity. We derived a minimax lower bound for estimation of these functions in high-dimensional regression setups, extending existing results obtained earlier {\em only} for isotropic functions. We formalized that Bayesian forests attain the near-optimal posterior concentration rate for these general function classes without any need for prior modification.

Our results are extended to a general class of estimation problems including nonparametric regression with a fixed and random design, binary classification, and density estimation.
Although we do not consider further nonparametric statistical models with BART priors in view of the length of the work, there are many other possible directions, such as mean-variance function estimation \citep{pratola2020heteroscedastic} and causal inference \citep{hahn2020bayesian}. Refer to \citet{linero2017review}, \citet{tan2019bayesian}, and \citet{hill2020bayesian} for extensive surveys of the application of BART to various nonparametric models. Because our Lemmas~\ref{lmm:priorcon}--\ref{lmm:cond3} enjoy a model-free framework, they will also be useful in investigating the posterior contraction rates for other statistical models.


\acks{We are grateful to two referees and the associate editor for valuable comments and constructive suggestions. We also thank Qurie Moon for sharing the code for BART with the exponentially decaying prior distribution.	
	 Seonghyun Jeong was supported by the Yonsei University Research Fund of 2021-22-0032 and by the National Research Foundation of Korea (NRF) grant funded by the Korea government (MSIT) (NRF-2022R1C1C1006735). 
	Veronika Ro{\v{c}}kov{\'a} gratefully acknowledges support from the James S. Kemper Foundation Faculty Research Fund at the University of Chicago Booth School of Business and the National Science Foundation (DMS:1944740).}

\appendix
\section*{Appendix}

\section{Technical Proofs}

\subsection{Proof of Theorem~\ref{thm:approx2} and Corollary~\ref{cor:approx}}
\label{sec:prooflmm1}

We can prove Theorem~\ref{thm:approx2} and Corollary~\ref{cor:approx} using a suitably chosen approximator $\hat f_0$ of $f_0$. The following proof shows that an approximator can be constructed with the step-heights evaluated at any $y_{rk}\in\Omega_{rk}^\circ\cap \Xi_r^\ast$, $r=1,\dots,R$, $k=1,\dots,2^{L}$.

\bigskip

\begin{myprooftitle}[Proof of Theorem~\ref{thm:approx2}]
	As $\mathcal T^\ast$ is chosen as in \eqref{eqn:treestar} and every $\Xi_r^\ast$ is regular, note that $\Omega_{rk}^\circ\cap \Xi_r^\ast$ is not empty for every $r$ and $k$.
	We fix any $y_{rk}\in\Omega_{rk}^\circ\cap \Xi_r^\ast$ and let $\hat f_0(x)=\sum_{r,k}\mathbbm 1_{\Omega_{rk}^\circ}(x) \beta_{r k,0}$ for $\beta_{rk,0}=f_0(y_{rk})$, so that 
	\begin{align*}
		f_0(x)-\hat f_0(x)&=\sum_{r=1}^R \sum_{k=1}^{2^{L}} \mathbbm 1_{\Omega_{rk}^\circ}(x) ( f_0(x) - f_0(y_{rk}) ) .
	\end{align*}
	In what follows, we write $S_0=\{s_{0,1},\dots,s_{0,d}\}$. 
	We verify the assertion for each of the given metrics.

	{\it Verification of \ref{lmm0:sup1} and \ref{lmm0:sup2}}: We first prove \ref{lmm0:sup2}. Fix $r$ and $k$.
	For any $x\in\Omega_{rk}^\circ$, define
	\begin{align}
		x\mapsto x^\ast: x^\ast=\underset{z\in \mathsf{cl}(\Omega_{rk}^\circ\cap \Xi_r^\ast)}{\text{arg\,min}}\lVert x-z\rVert_1,
		\label{eqn:xstar}
	\end{align}
	where $\mathsf{cl}(\cdot)$ denotes the closure of a set.
	If $x\in \mathsf{cl}(\Omega_{rk}^\circ\cap \Xi_r^\ast)$, it is trivial that $x^\ast=x$, which gives $|f_0(x)-f_0(x^\ast)|=0$. 
	If $x\notin \mathsf{cl}(\Omega_{rk}^\circ\cap \Xi_r^\ast)$, there exists $r'\ne r$ such that $x\in \Xi_{r'}^\ast$ for $\Xi_{r'}^\ast$ that is contiguous to $\Xi_r^\ast$, and hence, $x^\ast\in \mathsf{cl}(\Xi_r^\ast)\cap \mathsf{cl}(\Xi_{r'}^\ast)$. In this case, we have $x\ne x^\ast$ but $x_j=x_j^\ast$ for $j\notin S_0^\ast\subseteq S_0$, where $x_j$ and $x_j^\ast$ are the $j$th entries of $x$ and $x^\ast$, respectively.
	As $f_0$ is continuous and $x,x^\ast\in\mathsf{cl}(\Xi_{r'}^\ast)$, we obtain
	$|f_0(x)-f_0(x^\ast)|=|h_0(x_{S_0})-h_0(x_{S_0}^\ast)|\le\lambda \sum_{j=1}^d  |x_{s_{0,j}}-x_{s_{0,j}}^\ast|^{\alpha_{r'j}}$.
	It follows that, for any $x\in\Omega_{rk}^\circ$ with given $r$ and $k$, 
	\begin{align}
		|f_0(x)-f_0(x^\ast)|\le\lambda \sum_{j=1}^d  |x_{s_{0,j}}-x_{s_{0,j}}^\ast|^{\min_{r,j}\alpha_{rj}}\le \lambda |S_0^\ast| c_n^{\min_{r,j}\alpha_{rj}},
		\label{eqn:qaaq}
	\end{align} 
	since $\lVert x-x^\ast \rVert_\infty\le c_n$ and $x_j=x_j^\ast$ for $j\notin S_0^\ast$.
	Hence, by the triangle inequality, for any $x\in\Omega_{rk}^\circ$,
	\begin{align}
		|f_0(x)-f_0(y_{rk})| \le  \lambda |S_0^\ast| c_n^{\min_{r,j}\alpha_{rj}}  + |f_0(x^\ast)-f_0(y_{rk})|.
		\label{eqn:appsup}
	\end{align}
	Let  $\{\tilde \Omega_{r1}^\circ,\dots, \tilde \Omega_{r2^L}^\circ\}$ be the tree partition and $(l_{r1},\dots l_{rd})^\top$ be the counter vector 
	 returned by $\mathsf{Akd}(\Xi_r^\ast;\mathcal Z,\alpha,L,S_0)$ such that $L=\sum_{j=1}^d l_{rj}$, $r=1,\dots,R$.
	As $x^\ast,y_{rk}\in\mathsf{cl}(\Omega_{rk}^\circ\cap \Xi_r^\ast)\subseteq \mathsf{cl}(\Xi_r^\ast)$ and $f_0$ is continuous,
	\begin{align}
		|f_0(x^\ast)-f_0(y_{rk})| \le  \lambda\sum_{j=1}^d \mathsf{len}([\tilde\Omega_{rk}^\circ]_{s_{0,j}})^{\alpha_{rj}}  \lesssim \lambda\sum_{j=1}^d 2^{-\alpha_{rj} l_{rj} }.
		\label{eqn:approbo0}
	\end{align}
	Let $\tilde l_{rj}=L\bar\alpha/(d\alpha_{rj})$ for $r=1,\dots, R$, $j=1,\dots, d$, such that $\alpha_{r1} \tilde l_{r1}=\dots=\alpha_{rd} \tilde l_{rd}$ and  $L=\sum_{j=1}^d \tilde l_{rj}$ for every $r$ (note that $\tilde l_{rj}$ may not be integers).
	Then, it can be easily seen that  $l_{rj}> \tilde l_{rj}-1$ for every $r,j$, and hence
	\begin{align}
		\lambda\sum_{j=1}^d 2^{-\alpha_{rj} l_{rj} }\le 2 \lambda\sum_{j=1}^d 2^{-\alpha_{rj} \tilde l_{rj} }
		\le 2 \lambda d 2^{-\bar\alpha L/d }.
		\label{eqn:approbo}
	\end{align}
	Putting the bounds together for every $r$ and $k$, we obtain
	\begin{align*}
		\lVert f_0-\hat f_0\rVert_\infty \lesssim   \lambda |S_0^\ast| c_n^{\min_{r,j}\alpha_{rj}} +  \lambda d 2^{-\bar\alpha L/d}.
	\end{align*}
	This verifies \ref{lmm0:sup2}.
	
	Now, to prove \ref{lmm0:sup1}, note that $c_n=0$ implies $\mathsf{cl}(\Omega_r^\ast)=\mathsf{cl}(\Xi_r^\ast)$ although possibly $\Omega_r^\ast\ne \Xi_r^\ast$. That is, $(\Omega_r^\ast\cup \Xi_r^{\ast})\cap(\Omega_r^\ast\cap \Xi_r^{\ast})^c$ is a null set with measure zero for every $r$. Therefore, in evaluating the $L_\infty$-norm $\lVert f_0-\hat f_0\rVert_\infty$ with the essential supremum, we can ignore such a null set and focus on $\Omega_r^\ast\cap \Xi_r^{\ast}$.
	If $x\in\Omega_{rk}^\circ\cap \Xi_r^\ast$, then similar to \eqref{eqn:approbo0} and \eqref{eqn:approbo}, we obtain that 
	$
		|f_0(x)-f_0(y_{rk})| \lesssim \lambda\sum_{j=1}^d 2^{-\alpha_{rj} l_{rj} }\lesssim \lambda d 2^{-\bar\alpha L/d}
	$
since $x,y\in\Xi_r^\ast$. 	Putting the bounds together for every $r$ and $k$, we conclude the assertion.

	{\it Verification of \ref{lmm0:Leb1} and \ref{lmm0:Leb2}}: To verify \ref{lmm0:Leb1}, we first show that when $f_0\in\Gamma^{A_{\bar\alpha},d,p}_{\lambda}(\mathfrak X_0)$, for any finite measure $\mu$ and any fixed $v\ge1$, 
	\begin{align}
		\lVert f_0-\hat f_0\rVert_{v,\mu}\lesssim \tilde\epsilon_n, \quad\text{if}\quad\sum_{r=1}^{R} \mu(\Omega_r^\ast\cap \Xi_r^{\ast c})\lesssim (\tilde\epsilon_n/\lVert f_0\rVert_\infty)^v.
		\label{eqn:genmu}
	\end{align}
	Observe that
	\begin{align}
		\int |f_0(x)-\hat f_0(x)|^v d\mu(x)&=\sum_{r=1}^R \sum_{k=1}^{2^{L}} \int_{\Omega_{rk}^\circ} | f_0(x) -f_0(y_{rk}) |^v d \mu(x).
		\label{eqn:appl2}
	\end{align}
	The integral term in each summand is bounded by
	\begin{align}
		\int_{\Omega_{rk}^\circ\cap\Xi_r^\ast}  | f_0(x)-f_0(y_{rk})|^v d\mu(x)+    \mu(\Omega_{rk}^\circ\cap \Xi_r^{\ast c}) (2\lVert f_0\rVert_\infty)^v.
		\label{eqn:appl3}
	\end{align}
	Using \eqref{eqn:approbo0} and \eqref{eqn:approbo}, observe that, for every $x,y_{rk}\in \Omega_{rk}^\circ\cap\Xi_r^\ast$,
	\begin{align*}
		|f_0(x)-f_0(y_{rk})| \lesssim  \lambda d 2^{-\bar\alpha L/d }.
	\end{align*}
	The first term of \eqref{eqn:appl3} is thus bounded by a constant multiple of $\mu(\Omega_{rk}^\circ\cap \Xi_r^\ast)(\lambda d 2^{-\bar\alpha L/d })^v$.
	Note also that $\sum_k \mu(\Omega_{rk}^\circ\cap \Xi_r^{\ast c}) = \mu(\Omega_r^\ast\cap \Xi_r^{\ast c})$.
	Therefore,
	\begin{align}
		\begin{split}
			\lVert f_0-\hat f_0\rVert_{v,\mu}^v&\lesssim \sum_{r=1}^R \sum_{k=1}^{2^{L}} \left\{\mu(\Omega_{rk}^\circ\cap \Xi_r^\ast)\left(\lambda d 2^{-\bar\alpha L/d }\right)^v  +  \mu(\Omega_{rk}^\circ\cap \Xi_r^{\ast c}) \lVert f_0\rVert_\infty^v \right\}\\
			&\le \mu([0,1]^p)\left(\lambda d 2^{-\bar\alpha L/d }\right)^v + \lVert f_0\rVert_\infty^v\sum_{r=1}^{R} \mu(\Omega_r^\ast\cap \Xi_r^{\ast c}).
		\end{split}
		\label{eqn:app11}
	\end{align}
	This leads to the assertion in \eqref{eqn:genmu}. 
	Now, to verify the first part of \ref{lmm0:Leb1}, it suffices to show that $\sum_{r=1}^{R} \text{Leb}_p(\Omega_r^\ast\cap \Xi_r^{\ast c})\lesssim (\tilde\epsilon_n/\lVert f_0\rVert_\infty)^v$ for $\text{Leb}_p$, the Lebesgue measure on a $p$-dimensional space.
	For each $r$, we only need to consider the case $\Xi_r^\ast\subsetneq\Omega_{r}^\ast$, as $\text{Leb}_p(\Omega_{r}^\ast\cap \Xi_r^{\ast c})$ is maximized in this case.
	Then, $\Omega_{r}^\ast\cap \Xi_r^{\ast c}$ is not a box but a $p$-dimensional orthogonal polyhedron (for example, with a rectangular hole).
	One can easily see that 
	\begin{align*}
		\text{Leb}_p(\Omega_{r}^\ast\cap \Xi_r^{\ast c})
		\le \sum_{j=1}^p \text{Leb}_1([\Omega_r^\ast \cap \Xi_r^{\ast c}]_j)\prod_{k\ne j} \mathsf{len}([\Omega_r^\ast]_k).
	\end{align*}	
	It should be noticed that $[\Omega_r^\ast \cap \Xi_r^{\ast c}]_j$ may not be an interval but can be an empty set or a union of two isolated intervals. As $\text{Leb}_1([\Omega_r^\ast \cap \Xi_r^{\ast c}]_j)=0$ for $j\notin S_0^\ast\subseteq S_0$ and $\max_j\text{Leb}_1([\Omega_r^\ast \cap \Xi_r^{\ast c}]_j)\le 2 c_n$, the last expression is bounded by
	\begin{align*}
		|S_0^\ast|\max_j\left\{\text{Leb}_1([\Omega_r^\ast \cap \Xi_r^{\ast c}]_j)\prod_{k\ne j} \mathsf{len}([\Omega_r^\ast]_k)\right\}\le 
		\frac{2c_n|S_0^\ast|\mathsf{vol}(\Omega_r^\ast)}{\min_j\mathsf{len}([\Omega_r^\ast]_j)},
	\end{align*}
	where we use the notation $\mathsf{vol}(\cdot)$ to denote the volume of a box.
	As $\mathsf{len}([\Omega_r^\ast]_j)\ge \mathsf{len}([\Xi_r^\ast]_j)-2 c_n$ for every $j$,
	\begin{align}
		\sum_{r=1}^R \text{Leb}_p(\Omega_r^\ast \cap \Xi_r^{\ast c}) \le \frac{2c_n|S_0^\ast|}{\min_{r,j}\mathsf{len}([\Xi_r^\ast]_j)-2c_n} \le \frac{3c_n|S_0^\ast|}{\min_{r,j}\mathsf{len}([\Xi_r^\ast]_j)},
		\label{eqn:lebbound}
	\end{align}
	for every small $c_n>0$.
	It follows from this that $\sum_{r=1}^{R} \text{Leb}_p(\Omega_r^\ast\cap \Xi_r^{\ast c})\lesssim (\tilde\epsilon_n/\lVert f_0\rVert_\infty)^v$ if $c_n\lesssim ({\tilde\epsilon_n }/{\lVert  f_0 \rVert_\infty})^v {\min_{r,j}\mathsf{len}([\Xi_r^\ast]_j)/|S_0^\ast|}$.
	The first part of \ref{lmm0:Leb1} is verified. 
	
	We now verify \ref{lmm0:Leb2}.
	Similar to \eqref{eqn:genmu}, we first show that when $f_0\in\Gamma^{A_{\bar\alpha},d,p}_{\lambda}(\mathfrak X_0)\cap\mathcal C([0,1]^p)$, for any finite measure $\mu$ and any $v\ge1$,  
	\begin{align}
		\lVert f_0-\hat f_0\rVert_{v,\mu}\lesssim \tilde\epsilon_n, \quad\text{if}\quad c_n^{v\min_{r,j}\alpha_{rj}}\sum_{r=1}^{R} \mu(\Omega_r^\ast\cap \Xi_r^{\ast c})\lesssim (\tilde\epsilon_n/(\lambda |S_0^\ast|))^v.
		\label{eqn:genmu2}
	\end{align}
	We start from the identity in \eqref{eqn:appl2}. Similar to the above, one can observe that the integral term in \eqref{eqn:appl2} is bounded by
	\begin{align}
		\begin{split}
			&\mu(\Omega_{rk}^\circ\cap \Xi_r^\ast)\left( \lambda d 2^{-\bar\alpha L/d }\right)^v +    \int_{\Omega_{rk}^\circ\cap\Xi_r^{\ast c} }   | f_0(x)-f_0(y_{rk})|^v d\mu(x) .
		\end{split}
		\label{eqn:appl222}
	\end{align}
	Using $x^\ast\in \mathsf{cl}(\Omega_{rk}^\circ\cap \Xi_r^\ast)$ in \eqref{eqn:xstar}, the second term of \eqref{eqn:appl222} is bounded by
	\begin{align*}
		&2^{v-1}  \int_{\Omega_{rk}^\circ\cap\Xi_r^{\ast c}} ( | f_0(x)-f_0(x^\ast)|^v+| f_0(x^\ast)-f_0(y_{rk})|^v) d\mu(x) \\
		&\quad \le   2^{v-1} \mu(\Omega_{rk}^\circ\cap \Xi_r^{\ast c}) \left\{  
		\left(\lambda |S_0^\ast| c_n^{\min_{r,j}\alpha_{rj}}\right)^v +\left( \lambda d 2^{-\bar\alpha L/d } \right)^v\right\},
	\end{align*}
	where the inequality holds by \eqref{eqn:qaaq} combined with the fact that $x,x^\ast\in \mathsf{cl}(\Xi_{r'}^\ast)$ and $x^\ast,y_{rk}\in \mathsf{cl}(\Xi_r^\ast)$ for some $r'\ne r$. Hence, \eqref{eqn:appl222} is further bounded by a constant multiple of
	\begin{align*}
		\mu(\Omega_{rk}^\circ)\left(\lambda d 2^{-\bar\alpha L/d }\right)^v +   \mu(\Omega_{rk}^\circ\cap \Xi_r^{\ast c})\left(\lambda |S_0^\ast| c_n^{\min_{r,j}\alpha_{rj}}\right)^v,
	\end{align*}
	and we obtain that
	\begin{align}
		\begin{split}
			\lVert f_0-\hat f_0\rVert_{v,\mu}^v&\lesssim \sum_{r=1}^R \sum_{k=1}^{2^{L}} \left\{\mu(\Omega_{rk}^\circ)\left(\lambda d 2^{-\bar\alpha L/d }\right)^v  +    \mu(\Omega_{rk}^\circ\cap \Xi_r^{\ast c})\left(\lambda |S_0^\ast| c_n^{\min_{r,j}\alpha_{rj}}\right)^v \right\}\\
			&\le \mu([0,1]^p)\left(\lambda d 2^{-\bar\alpha L/d }\right)^v + \left(\lambda |S_0^\ast| c_n^{\min_{r,j}\alpha_{rj}}\right)^v \sum_{r=1}^R \mu(\Omega_r^\ast\cap \Xi_r^{\ast c}).
		\end{split}
		\label{eqn:app22}
	\end{align}
	This leads to \eqref{eqn:genmu2}. Now, to verify the second part of \ref{lmm0:Leb1}, we take the Lebesgue measure for $\mu$.
	Then using the bound in \eqref{eqn:lebbound}, we have that $(\lambda |S_0^\ast|c_n^{\min_{r,j}\alpha_{rj}})^v \sum_{r=1}^R \text{Leb}_p(\Omega_{r}^\ast\cap \Xi_{r}^{\ast c})\lesssim \tilde\epsilon_n^v$ if $c_n^{1+v\min_{r,j}\alpha_{rj}}\lesssim ({\tilde\epsilon_n }/\lambda)^v {\min_{r,j}\mathsf{len}([\Xi_r^\ast]_j)}/|S_0^\ast|^{v+1}$. This proves the assertion.

	{\it Verification of \ref{lmm0:emp1}}:	We again use the result in \eqref{eqn:genmu}.
	Take $P_{\mathcal Z}$ for $\mu$. Then, it can be seen that split-points can be picked up such that there are no $z_i$ on $\cup_r(\Omega_r^\ast\cap \Xi_r^{\ast c})$ by choosing the points closest to the boundaries in every split. As we have $\sum_{r=1}^R P_{\mathcal Z}(\Omega_r^\ast \cap \Xi_r^{\ast c})=0$ in this case, \ref{lmm0:emp1} easily follows.
\end{myprooftitle}

\begin{myprooftitle}[Proof of Corollary~\ref{cor:approx}]
{\it Verification of \ref{lmm0:sup2-1}}:
As $\bar\epsilon_n\gtrsim (\lambda d R (\log n )/n)^{1/3} $ and $|S_0^\ast|\le d$, we obtain $\bar\epsilon_n /(\lambda|S_0^\ast|)\gtrsim (n^{-1}R\log n)^{1/3} (\lambda d)^{-2/3}\ge n^{-(1+2a_2)/3}\log^{-1/3} n$. The assertion in \ref{lmm0:sup2-1} follows by combining \ref{lmm0:sup2} of Theorem~\ref{thm:approx2} and the bound $\min_{r,j}\alpha_{rj}\ge a_1$.

{\it Verification of \ref{lmm0:Leb1-1}}:
Similar to above, we obtain
\begin{align*}
(\bar\epsilon_n/\lVert f_0\rVert_\infty)^v  \min_{r,j}\mathsf{len}([\Xi_r^\ast]_j)/|S_0^\ast|&\gtrsim (n^{-1}\lambda dR\log n)^{v/3}(\log n)^{-v/2}n^{-a_3}/d\\
&\gtrsim \lambda^{v/3} n^{-(v/3+a_3)}(\log n)^{-(\max\{0,1-v/3\}+v/6)}.
\end{align*}
As $\lambda\gtrsim 1$, we can verify the assertion in \ref{lmm0:Leb1-1} using \ref{lmm0:Leb1} of Theorem~\ref{thm:approx2}.
\end{myprooftitle}

\subsection{Proof of Lemmas~\ref{lmm:grid2}--\ref{lmm:fixed2}}
\label{sec:prooflmm2}

To prove Lemma~\ref{lmm:grid2}, we first provide the following lemma, which shows that a regular grid is dense and regular for arbitrary inputs under mild conditions. 

\begin{mylemma}[Regular grid, general case]
	For a regular grid $\mathcal Z$, we have the following assertions.
	\begin{enumerate}[label=\rm(\roman*)]
		\item \label{lmm2:1} For any $S\subseteq\{1,\dots,p\}$ and any $S$-chopped flexible tree partition $\mathfrak Y=\{\Psi_1,\dots,\Psi_J\}$ with $J\ge 2$, $\mathcal Z$ is $(\mathfrak Y,1/b_n^{1/p})$-dense if $\min_{r,j}\mathsf{len}([\Psi_r]_j)\ge b_n^{-1/p}$.
		\item \label{lmm2:2} For any $S\subseteq\{1,\dots,p\}$, $\alpha\in(0,1]^d$, $\Psi\subseteq[0,1]^p$, and  $L = \lfloor\log_2 ( b_n^{1/p}\min_j\mathsf{len}([\Psi]_j)-1)\rfloor$, $\mathcal Z$ is $(\Psi,\alpha,L,S)$-regular if $\min_{r,j}\mathsf{len}([\Psi_r]_j)\ge 3b_n^{-1/p}$.
	\end{enumerate}
	\label{lmm:grid}
\end{mylemma}

\begin{myproof}
	{\it Verification of \ref{lmm2:1}}:
	Consider a $p$-dimensional checkerboard $\prod_{j=1}^p [(i_j-1)/b_n^{1/p},i_j/b_n^{1/p}]$, $i_j=1,\dots, b_n$.
	Note that each point $z_i$ in $\mathcal Z$ is located at the center of each box of this checkerboard. As the mesh-size of the checkerboard is $1/b_n^{1/p}$, there exists an $S$-chopped $\mathcal Z$-tree partition $\mathcal T$ such that $\Upsilon(\mathfrak Y,\mathcal T)\le 1/b_n^{1/p}$ if $\min_{r,j}\mathsf{len}([\Psi_r]_j)\ge 1/b_n^{1/p}$. The assertion easily follows.
	
	{\it Verification of \ref{lmm2:2}}:
	The condition $\min_{r,j}\mathsf{len}([\Psi_r]_j)\ge 3b_n^{-1/p}$ is made to ensure that there is at least one split-point that is sufficiently far away from the boundaries of $\Psi$ in every coordinate.
	Observe that for any box $\Psi\subseteq [0,1]^p$, we obtain
	\begin{align}
		\tilde b_j(\mathcal Z,\Psi)\le b_n^{1/p}\mathsf{len}([\Psi]_j) \le\tilde b_j(\mathcal Z,\Psi)+1,\quad j=1,\dots,p.
		\label{eqn:bbb}
	\end{align}
	Thus, in every coordinate, midpoint-splits can occur $\lfloor b_n^{1/p}\min_j\mathsf{len}([\Psi]_j)\rfloor-1$ times without choosing the leftmost and rightmost split-points (these two points may produce too small cells). This allows us to 
	choose 	$L = \lfloor\log_2 ( b_n^{1/p}\min_j\mathsf{len}([\Psi]_j)-1)\rfloor$ for an anisotropic $k$-d tree (note that $\lfloor\log\lfloor x \rfloor \rfloor=\lfloor\log x\rfloor$, $x>0$).
	
	For any $\Psi\subseteq [0,1]^p$ and $j\in S$, a mid-point split chooses $\lceil \tilde b_j(\mathcal Z,\Psi)/2 \rceil$th split-candidate in $[\mathcal Z]_j\cap {\mathsf{int}}([\Psi]_j)$ as a split-point $\tau_j$. The resulting two cells have at most $\lfloor\tilde b_j(\mathcal Z,\Psi)/2\rfloor$ split-points in coordinate $j$. 
	Therefore, using \eqref{eqn:bbb},
	\begin{align*}
		\max_k\mathsf{len}([\Omega_k^\circ]_{s_j})&\le \frac{\tilde b_j(\mathcal Z,\Psi)2^{-l_j}+1}{b_n^{1/p}}\\
		&\le \mathsf{len}([\Psi]_{s_j}) 2^{-l_j}+1/b_n^{1/p}\\
		&\le \mathsf{len}([\Psi]_{s_j}) \left(2^{-l_j}+\frac{1}{b_n^{1/p}\min_{r,j}\mathsf{len}([\Psi_r]_j)}\right).
	\end{align*}
	As $L\le\log_2 ( b_n^{1/p}\min_j\mathsf{len}([\Psi]_j)-1)\le\log_2 ( b_n^{1/p}\min_j\mathsf{len}([\Psi]_j))-1$ and $l_j\le L$ for every $j=1,\dots,d$, the last expression is bounded by
	\begin{align*}
		\mathsf{len}([\Psi]_{s_j})(2^{-l_j}+2^{1-L})\le 3\mathsf{len}([\Psi]_{s_j})2^{-l_j}.
	\end{align*}
	This leads to the assertion.
\end{myproof}

\begin{myprooftitle}[Proof of Lemma~\ref{lmm:grid2}]
	If $R=1$, it is obvious that $\mathcal Z$ is $(\mathfrak X_0^\ast, 0)$-dense.
	If $R>1$, by \ref{lmm2:1} of Lemma~\ref{lmm:grid}, $\mathcal Z$ is $(\mathfrak X_0^\ast, 1/b_n^{1/p})$-dense since $b_n^{1/p}\min_{r,j}\mathsf{len}([\Xi_r^\ast]_j)\gg 1$.
	Also, \ref{lmm2:2} of Lemma~\ref{lmm:grid} shows that  $\mathcal Z$ is $(\Xi_r^\ast,\alpha_r,L_r,S_0)$-regular for  $L_r=\lfloor\log_2 ( b_n^{1/p}\min_j\mathsf{len}([\Xi_r^\ast]_j)-1)\rfloor$, $r=1,\dots, R$.
	To conclude that  $\mathcal Z$ is $(\Omega_r^\ast,\alpha_r,L_0,S_0)$-regular for $r=1,\dots, R$, we only need to show that  $L_0\le\min_r L_r$. 
	As $2^{L_0}\asymp (n(\lambda d)^2/(R\log n))^{d/(2\bar\alpha+d)}$, $L_0$ can be chosen to be $2^{L_0} \le C_1(n(\lambda d)^2/(R\log n))^{d/(2\bar\alpha+d)}$ for small enough $C_1>0$ as desired. Therefore, a sufficient condition for $L_0\le \min_r L_r$ is $(n(\lambda d)^2/(R\log n))^{d/(2\bar\alpha+d)}\lesssim b_n^{1/p}\min_j\mathsf{len}([\Xi_r^\ast]_j)$. Plugging in $b_n=n^{cp}$, the conditions in the lemma are obtained.
\end{myprooftitle}

\begin{myprooftitle}[Proof of Lemma~\ref{lmm:fixed2}]
Assumption \ref{asm:fixeddesign} implies that $\mathcal Z$ is $(\Xi_r^\ast,\alpha_r,L_r,S_0)$-regular for $L_r = \lfloor\log_2(C_1 n P_{\mathcal Z}(\Xi_r^\ast))\rfloor$, for some $C_1>0$, $r=1,\dots,R$. It remains to show that $L_0\le \min_r L_r$.
	Recall that $2^{L_0}\asymp (n(\lambda d)^2/(R\log n))^{d/(2\bar\alpha+d)}$. As $L_0$ can be chosen to be $2^{L_0} \le C_2(n(\lambda d)^2/(R\log n))^{d/(2\bar\alpha+d)}$ for small enough $C_2>0$ as desired, a sufficient condition for $L_0\le \min_r L_r$ is given by $(n(\lambda d)^2/(R\log n))^{d/(2\bar\alpha+d)}\lesssim n P_{\mathcal Z}(\Xi_r^\ast)$ no matter what $C_1$ is. Using that $P_{\mathcal Z}(\Xi_r^\ast)\gtrsim R^{-1}$, the inequality is translated into $\lambda d \lesssim (n/R)^{\bar\alpha/d} \sqrt{\log n}$.
\end{myprooftitle}

\subsection{Proof of Theorem~\ref{thm:nonreg}}
\label{sec:proofthm1-2}

We deploy the standard theory on posterior contraction  \citep{ghosal2000convergence,ghosal2007convergence}. The required conditions for the general theory are deferred to Lemmas~\ref{lmm:priorcon}--\ref{lmm:cond3}.

\bigskip

\begin{myprooftitle}[Proof of Theorem~\ref{thm:nonreg}]
	As $\sigma_0^2$ is bounded below and above, $|\sigma^2-\sigma_0^2|$ and $|\sigma-\sigma_0|$ have the same rate. We will work with the latter for convenience. 
	We write $\rho_n^2((f_1,\sigma_1),(f_2,\sigma_2))=\lVert f_1-f_2\rVert_n^2+|\sigma_1-\sigma_2|^2$ for any $f_1,f_2:\mathbb R^p\rightarrow \mathbb R$ and any $\sigma_1,\sigma_2\in(0,\infty)$.
	(Observe that  $\sqrt{\lVert\cdot\rVert_n^2+|\cdot|^2}$ and $\lVert\cdot\rVert_n+|\cdot|$ have the same order.)
	By Lemma~1 of \citet{lim2023synergizing},
for every $\epsilon>0$ and $(f_1,\sigma_1)$ with $\lVert f_1-f_0\rVert_n^2 + |\sigma_1-\sigma_0|^2\ge \epsilon^2$, there exists a test $\phi_n$ such that, for a universal constant $K>0$,
	\begin{align*}
		\mathbb E_0 \phi_n \le e^{-Kn\epsilon^2},\quad 	\sup_{(f,\sigma^2):\lVert f-f_1\rVert_n^2 + |\sigma-\sigma_1|^2\le \epsilon^2/36}\mathbb E_{f,\sigma^2} (1-\phi_n) \le e^{-Kn\epsilon^2}.
	\end{align*}
	
	We write $\mathcal F_\ast=\cup_{\mathcal E}\mathcal F_{\mathcal E}$, where the union is taken over all $\mathcal E$ generated by a given $\mathcal Z$.
	For the Kullback-Leibler (KL) divergence $K(p_1,p_2)=\int \log(p_1/p_2)p_1$ and its second order variation $V(p_1,p_2)=\int |\log(p_1/p_2)-K(p_1,p_2)|^2 p_1$, define 
	\begin{align*}
	B_n=\left\{ (f,\sigma): \sum_{i=1}^n K(p_{0,i},p_{f,\sigma, i})\le n\epsilon_n^2 , \, \sum_{i=1}^n V(p_{0,i},p_{f,\sigma,i})\le n\epsilon_n^2\right\}.
	\end{align*}
	By Theorem~8.19 of \citet{ghosal2017fundamentals}, we only need to verify that there exists a sieve $\Theta_n\subseteq \mathcal F\times (0,\infty)$ such that 
	for some $\bar c>0$ and a sufficiently large $\bar c'>0$,
	\begin{align}
		\label{eqn:con1}
		\Pi(B_n)&\ge e^{-\bar c n\epsilon_n^2 },\\
		\label{eqn:con2}
		\log N(\epsilon_n, \Theta_n, \rho_n)&\lesssim n\epsilon_n^2,\\
		\label{eqn:con3}
		\Pi((f,\sigma)\notin \Theta_n)&\ll e^{-\bar c' n\epsilon_n^2 },
	\end{align}
	
	We first verify \eqref{eqn:con1}.
	By direct calculations, 
	\begin{align*}
		\frac{1}{n}\sum_{i=1}^n K(p_{0,i},p_{f,\sigma,i})&=\frac{1}{2} \log \left( \frac{\sigma^2}{\sigma_0^2} \right)-\frac{1}{2}\left(1-\frac{\sigma_0^2}{\sigma^2}\right)+\frac{\lVert f-f_0 \rVert_n^2}{2\sigma^2},\\
		\frac{1}{n}\sum_{i=1}^n V(p_{0,i},p_{f,\sigma,i}) &=\frac{1}{2}\left(1-\frac{\sigma_0^2}{\sigma^2}\right)^2+\frac{\sigma_0^2\lVert f-f_0 \rVert_n^2}{\sigma^2}.
	\end{align*}
	Using the Taylor expansion, it is easy to see that, for any $\epsilon_n\rightarrow 0$,
	there exists a constant $C_1>0$ such that
	\begin{align*}
		B_n&\supseteq \{(f,\sigma): \lVert f-f_0\rVert_n\le C_1 \epsilon_n , | \sigma-\sigma_0|\le C_1 \epsilon_n \}.
	\end{align*}
	First, note that $\log \Pi(\sigma^2: |\sigma-\sigma_0|\le C_1\epsilon_n)\gtrsim -\log n$ if $\sigma_0$ lies on a compact subset of $(0,\infty)$.
	We will construct a good approximating ensemble denoted by $\widehat{\mathcal E}=(\widehat{\mathcal T}^1,\dots,\widehat{\mathcal T}^T)$.
	By restricting the function space to the one constructed by $\widehat{\mathcal E}$, we obtain
	\begin{align}
		&	\Pi(f\in\mathcal F_\ast:\lVert f-f_0\rVert_n\le C_1\epsilon_n)\ge\Pi(\widehat{\mathcal E})\Pi(f\in\mathcal F_{\widehat{\mathcal E}}:\lVert f-f_0\rVert_n \le C_1\epsilon_n).
		\label{eqn:lnnb}
	\end{align}	
	Assumption \ref{asm:spnet} states that, for a given split-net $\mathcal Z$ there exists a $\mathcal Z$-tree partition $\widehat {\mathcal T}$ producing $\hat f_0\in\mathcal F_{\widehat {\mathcal T}}$ satisfying $\lVert f_0-\hat f_0 \rVert_n\lesssim \bar\epsilon_n$.
	An approximating ensemble $\widehat{\mathcal E}$ can be constructed by setting $\widehat{\mathcal T}^1$ to be $\widehat{\mathcal T}$ and $\widehat{\mathcal T}^t$, $t=2,\dots, T$, to be root nodes with no splits, i.e., $\widehat{\mathcal T}^t=\{[0,1]^p\}$, $t=2,\dots, T$.
	Then,
	\begin{align*}
		\log \Pi(\widehat{\mathcal E})= \sum_{t=1}^T\log \Pi(\widehat{\mathcal T}^t)= \log \Pi(\widehat{\mathcal T}^1)+(T-1) \log (1-\nu)
		\gtrsim  -n\epsilon_n^2,
	\end{align*}
	by Lemma~\ref{lmm:priorcon}.
	It remains to bound the second term of \eqref{eqn:lnnb}.
	By \ref{asm:spnet}, we have $\lVert f-f_0 \rVert_n
	\lesssim  \lVert f-\hat f_0 \rVert_\infty +\epsilon_n$ for some $\hat f_0\in\mathcal F_{\widehat{\mathcal T}}$. 	We can construct $\hat f_0$ as in the proof of Theorem~\ref{thm:approx2}. We denote this $\hat f_0$ by $f_{0,\widehat{\mathcal T},\widehat\beta}$, where $\widehat\beta$ is the corresponding step-heights, to emphasize the dependence on $\widehat{\mathcal T}$ and $\widehat\beta$.
	We shall now express $f_{0,\widehat{\mathcal T},\widehat\beta}$ using the approximating ensemble $\widehat{\mathcal E}$ with corresponding step-heights $\widehat B$. As all trees in $\widehat{\mathcal E}$ are the root nodes except for the first one $\widehat{\mathcal T}^1$, every step-heights vector $B$ for $\widehat{\mathcal E}$ has the form $B=({\beta^{1\top}},\beta^2,\dots,\beta^T)^\top\in\mathbb R^{\widehat K+T-1}$ with $\beta^1\in\mathbb R^{\widehat K}$ and $\beta^t\in\mathbb R$, $t=2,\dots, T$, where $\widehat K$ is the size of $\widehat{\mathcal T}$. Hence, letting $\widehat B=({\widehat\beta}^\top,0,\dots,0)^\top$, we can write $f_{0,\widehat{\mathcal T},\widehat\beta}=f_{0,\widehat{\mathcal E},\widehat B}$ for $f_{0,\widehat{\mathcal E},\widehat B}$ defined with the ensemble components $(\widehat{\mathcal E},\widehat B)$.
	Putting the bounds together, for some $C_2>0$,
	\begin{align*}
		\Pi\left(f\in\mathcal F_{\widehat{\mathcal E}}:\lVert f-f_0\rVert_n\le C_1\epsilon_n\right)&\ge \Pi\left(f\in\mathcal F_{\widehat{\mathcal E}}:\lVert f-f_{0,\widehat{\mathcal E},\widehat B} \rVert_\infty\le C_2\epsilon_n\right).
	\end{align*}	
By Lemma~\ref{lmm:cond1}, the right-hand side is bounded below as desired.
	Putting everything together, we conclude that there exists a constant $\bar c$ such that $\Pi(B_n)\ge e^{-\bar c n\epsilon_n^2 }$.
	
	Next, we verify the entropy condition \eqref{eqn:con2}.
	We denote by $\mathscr E_{S,K^1,\dots,K^T}$ the collection of $\mathcal E=\{\mathcal T^1,\dots,\mathcal T^T\}$ with given $S,K^1,\dots,K^T$; that is, each $\mathcal T^t$ is an $S$-chopped $\mathcal Z$-tree partition of size $K^t$. 
	With given $\mathcal E$ and $M>0$, we first define the function spaces $\mathcal F_{\mathcal E,M}^{(1)}=\{f_{\mathcal E, B}\in\mathcal F_{\mathcal E}: \lVert B\rVert_\infty\le M\}$ and $\mathcal F_{\mathcal E,M}^{(2)}=\{f_{\mathcal E, B}\in\mathcal F_{\mathcal E}: \lVert B\rVert_\infty > M\}$ such that $\mathcal F_{\mathcal E,M}^{(1)}\cup \mathcal F_{\mathcal E,M}^{(2)} = \mathcal F_{\mathcal E}$.
	We also define 
	\begin{align}
		\mathcal F_{\bar s_n,\bar K_n,M}^{(\ell)}:=\bigcup_{\mathcal E \in \mathscr E_{S,K^1,\dots,K^T}:  |S| \le \bar s_n, K^t\le \bar K_n,t=1,\dots,T}\mathcal F_{\mathcal E,M}^{(\ell)},\quad \ell=1,2,
		\label{eqn:fundef}
	\end{align}
	for $\bar K_n\asymp n\epsilon_n^2/\log n$ and $\bar s_n\asymp n\epsilon_n^2/\log p$. 
	That is, $\mathcal F_{\bar s_n,\bar K_n,M}^{(\ell)}$ is the collection of all $\mathcal F_{\mathcal E,M}^{(\ell)}$ such that $K^t\le \bar K_n$ and $|S|\le \bar s_n$.
	We take $\Theta_n=\mathcal F_{\bar s_n,\bar K_n,n^{M_1}}^{(1)}\times(n^{-M_2},e^{M_2 n\epsilon_n^2})$  for large $M_1,M_2>0$.
	It is easy to see that $\log N(\epsilon_n,(n^{-M_2},e^{M_2 n\epsilon_n^2}),{\lvert\cdot\rvert})\lesssim n\epsilon_n^2 $. Combining this with Lemma~\ref{lmm:cond2}, we conclude that \eqref{eqn:con2} is verified.

	Lastly, we verify \eqref{eqn:con3}. First,
	it is easy to see that $\Pi(\sigma^2\notin(n^{-2M_2},e^{2M_2 n\epsilon_n^2}))e^{\bar c' n\epsilon_n^2 }\rightarrow 0$ if $M_2$ is large enough, using the tail probabilities of inverse gamma distributions. Choose $\bar K_n=\lfloor M_3 n\epsilon_n^2/\log n \rfloor$ and $\bar s_n=\lfloor M_3 n\epsilon_n^2/\log p \rfloor$ for a sufficiently large $M_3>0$.
	As we have 
	$\Pi(\mathcal F_\ast\setminus \mathcal F_{\bar s_n,\bar K_n, n^{M_1}}^{(1)})e^{\bar c' n\epsilon_n^2 }\rightarrow 0$ by Lemma~\ref{lmm:cond3}, the condition is verified.
\end{myprooftitle}

\begin{mylemma}[Prior concentration of tree sizes]
	Let $\widehat{\mathcal T}$ be the $\mathcal Z$-tree partition defined in \eqref{eqn:treeapp}. Under Assumptions \ref{asm:spnet0} and \ref{asm:spnet1}, $\log \Pi( \widehat{\mathcal T})\gtrsim - \widehat K \log n - d\log p$.
	\label{lmm:priorcon}
\end{mylemma}

\begin{myproof}
	We will obtain a lower bound of $\Pi(\widehat{\mathcal T})$. As this depends on splitting proportions drawn from a Dirichlet prior, we first restrict the proportions to the set
	\begin{align*}
		V_1=\left\{\eta\in\mathbb S^p:\eta_j\ge \frac{1}{2d}, j\in S_0 , \sum_{j\notin S_0}\eta_j \le \frac{1}{2d}\right\}.
	\end{align*}
	Fix $\eta^\ast=(\eta_1^\ast,\dots,\eta_p^\ast)^\top\in\mathbb S^p$ such that $\eta_j^\ast=1/d$, $j\in S_0$, and $\eta_j^\ast=0$, $j\notin S_0$. It can be easily shown that $V_1\supseteq \{\eta\in\mathbb S^p: \lVert\eta-\eta^\ast\rVert_1\le 1/(2d)\}$. By \eqref{eqn:dirineq1} of Lemma~\ref{lmm:dir}, it follows that $\Pi(V_1)\ge e^{-C_1 d\log p}$ for some $C_1>0$.
	Recall that the first $R-1$ splits of $\widehat{\mathcal T}$ form $\mathcal T^\ast=\{\Omega_1^\ast,\dots,\Omega_R^\ast\}$, the approximating tree partition of $\mathfrak X_0^\ast$, and the remaining splits generate $\mathcal T_r^\circ$, the tree partition of $\Omega_r^\ast$ constructed by an anisotropic $k$-d tree, $r=1,\dots,R$. 
	Hence, we can write
	\begin{align*}
		\Pi(\widehat{\mathcal T})
		&\ge e^{-C_1 d\log p}\Pi(\widehat{\mathcal T}|V_1)=e^{-C_1 d\log p}\Pi(\text{$\mathcal T^\ast$ is a pruned tree of $\widehat{\mathcal T}$}|V_1)\prod_{r=1}^R \Pi(\mathcal T_r^\circ|\Omega_r^\ast,V_1).
	\end{align*}
	
	We first focus on the prior probability $\Pi(\text{\text{$\mathcal T^\ast$ is a pruned tree of $\widehat{\mathcal T}$}}|V_1)$. To generate $\mathcal T^\ast$, the root node is subdivided $R-1$ times in a top-down manner.
	As each node splits with probability $\nu^{\ell+1}$ for depth $\ell$, this occurs with probability at least $\nu^{(R-1)\max_r \mathsf{dep}(\Omega_r^\ast)}$ no matter what the partition is. Note also that, for every split, there are at most $\max_{1\le j\le p}b_j(\mathcal Z)$ splitting points and a splitting coordinate $j$ is chosen by $\eta_j$, $j\in S_0$, which is at least $1/(2d)$ on $V_1$. Hence the prior probability of choosing the correct split is bounded below by $1/(2d \max_{1\le j\le p}b_j(\mathcal Z))$ for every split. This gives us a lower bound:
	\begin{align*}
		\Pi(\text{$\mathcal T^\ast$ is a pruned tree of $\widehat{\mathcal T}$}|V_1)\ge \frac{ \nu^{(R-1)\max_r \mathsf{dep}(\Omega_r^\ast)} }{(2d\max_{1\le j\le p}b_j(\mathcal Z) )^{R-1}}.
	\end{align*}
	It follows that
	$\log\Pi(\text{$\mathcal T^\ast$ is a pruned tree of $\widehat{\mathcal T}$}|V_1) \gtrsim
	-R\log n$ as $\log (2d\max_{1\le j\le p}b_j(\mathcal Z) )\lesssim \log n$ by \ref{asm:spnet0} and $\max_r \mathsf{dep}(\Omega_r^\ast)\lesssim \log n$ by \ref{asm:spnet1}.

	We now obtain a lower bound of $\Pi(\mathcal T_r^\circ|\Omega_r^\ast,V_1)$.
	In splitting each $\Omega_r^\ast$, observe that $2^k$ cells split at depth $k=0,\dots,L_0-1$, and each cell splits with probability $\nu^{\mathsf{dep}(\Omega_r^\ast)+k+1}$ at depth $k$.
	Note that closing each of the terminal nodes is of probability at least $1-\nu$ and there are $2^{L_0}$ terminal nodes.
	Hence, similar to the above, 
	\begin{align*}
		\Pi(\mathcal T_r^\circ|\Omega_r^\ast,V_1)&\ge (1-\nu)^{2^{L_0}} \prod_{k=0}^{L_0-1} \left(\frac{\nu^{\mathsf{dep}(\Omega_r^\ast)+k+1}}{2d\max_{1\le j\le p}b_j(\mathcal Z)}\right)^{2^k} \\
		&= (1-\nu)^{2^{L_0}} \frac{\nu^{(\mathsf{dep}(\Omega_r^\ast)+1)(2^{L_0}-1)+(L_0-2)2^{L_0}+2} }{(2d\max_{1\le j\le p}b_j(\mathcal Z))^{2^{L_0}-1}},
	\end{align*}
	where we used the formulae
	$\sum_{k=0}^{a-1} 2^k=2^{a}-1$ and $\sum_{k=0}^{a-1} k 2^k= (a-2)2^a+2$.
	This yields
	$\sum_{r=1}^R\log \Pi(\mathcal T_r^\circ|\Omega_r^\ast,V_1) \gtrsim   - R 2^{L_0} \log (2d\max_{1\le j\le p}b_j(\mathcal Z) )- R 2^{L_0} \max_r\mathsf{dep}(\Omega_r^\ast) - RL_0 2^{L_0}\gtrsim - R2^{L_0}\log n $ since $L_0\lesssim \log n$.
	
	Putting everything together, we thus obtain $\log\Pi(\widehat{\mathcal T})\gtrsim -R2^{L_0}\log n-d\log p$.
	As $\widehat K=R2^{L_0}$, this verifies the assertion.
\end{myproof}

\begin{mylemma}[Prior concentration of tree learners] Define $\widehat{\mathcal E}$ and $\widehat B$ as in the proof of Theorem~\ref{thm:nonreg}.
Under \ref{asm:fsup} and \ref{pri:normal}, for any $C>0$,
\begin{align*}
 - \log \Pi\left(f\in\mathcal F_{\widehat{\mathcal E}}:\lVert f-f_{0,\widehat{\mathcal E},\widehat B} \rVert_\infty\le C\epsilon_n\right)\lesssim n\epsilon_n^2.
\end{align*}
\label{lmm:cond1}
\end{mylemma}

\begin{myproof}
	For any step-heights $B_1=({\beta_1^{1\top}},\beta_1^2,\dots,\beta_1^T)^\top,B_2=({\beta_2^{1\top}},\beta_2^2,\dots,\beta_2^T)^\top\in\mathbb R^{\widehat K+T-1}$ with $\widehat{\mathcal E}$, we write $f_{\widehat{\mathcal E},B_1}$, $f_{\widehat{\mathcal E},B_2}\in\mathcal F_{\widehat{\mathcal E}}$ to denote two additive tree functions that lie on the same partition ensemble $\widehat{\mathcal E}$. Evidently,
\begin{align*}
	\lVert f_{\widehat{\mathcal E},B_1}-f_{\widehat{\mathcal E},B_2} \rVert_\infty  &=\left\lVert\sum_{t=1}^T \beta_{1}^t-\sum_{t=1}^T \beta_{2}^t\right\rVert_\infty 
	\le \lVert \beta_1^1-\beta_2^1 \rVert_1+\sum_{t=2}^T|\beta_1^t-\beta_2^t| 
	\le \lVert B_1-B_2\rVert_2\sqrt{\widehat K_\ast},
\end{align*}	
where $\widehat K_\ast=\widehat K+T-1$.
It follows that, for some $C_1>0$,
\begin{align*}
\Pi\left(f\in\mathcal F_{\widehat{\mathcal E}}:\lVert f-f_{0,\widehat{\mathcal E},\widehat B} \rVert_\infty\le C\epsilon_n\right)\ge	\Pi\left( B \in\mathbb R^{\widehat K_\ast} :\lVert B - \widehat B\rVert_2\le {C_1\epsilon_n}/\sqrt{\widehat K_\ast}\right).
\end{align*}
Recall that the eigenvalues of the covariance matrix for a normal prior is bounded below and above.
This means that there exists an invertible matrix $D\in\mathbb R^{\widehat K_\ast\times\widehat K_\ast}$ such that $DB$ has a product of independent standard normal priors. Following the computations in page 216 of \citet{ghosal2007convergence}, the last display is further bounded below by
\begin{align}
	\begin{split}
		&\Pi\left(B \in\mathbb R^{\widehat K_\ast} :\lVert D(B-\widehat B)\rVert_2\le C_1\epsilon_n \sigma_{\max}^{-1}(D^{-1}) /\sqrt{\widehat K_\ast}\right)\\
		&\quad\ge 2^{-\widehat K_\ast/2}e^{-\Vert D \widehat B\rVert_2^2}\Pi\left(B\in\mathbb R^{\widehat K_\ast} :\lVert D B\rVert_2\le C_1\epsilon_n\sigma_{\max}^{-1}(D^{-1}) /\sqrt{2 \widehat K_\ast}\right),
		\label{eqn:norbo}
	\end{split}
\end{align}	
where $\sigma_{\max}(D^{-1})$ is the spectral norm of $D^{-1}$, which is bounded by the assumption. As the induced prior for $\lVert D B\rVert_2^2$ is a chi-squared distribution with degree of freedom $\widehat K_\ast$, we obtain that for $\upsilon_n=\epsilon_n\sigma_{\max}^{-1}(D^{-1})/\sqrt{\widehat K_\ast}\lesssim \epsilon_n$,
\begin{align*}
	\Pi(B\in\mathbb R^{\widehat K_\ast} :\lVert D B\rVert_2\le C_1\upsilon_n/\sqrt{2})\ge \frac{2/\widehat K_\ast}{2^{\widehat K_\ast}\Gamma(\widehat K_\ast/2)} (C_1\upsilon_n)^{\widehat K_\ast} e^{-C_1^2\upsilon_n^2/4}.
\end{align*}			
The logarithm of the right-hand side is bounded below by a constant multiple of $-(\widehat K+T)\log n-\upsilon_n^2\gtrsim -n\epsilon_n^2$. 
It only remains to bound $e^{-\Vert D \widehat B\rVert_2^2}$ in \eqref{eqn:norbo}. Observe that $\lVert\widehat\beta\rVert_\infty=\Vert f_{0,\widehat{\mathcal T}, \widehat \beta}\rVert_\infty\le\lVert f_0\rVert_\infty $, where the inequality follows from our choice of $\hat f_0=f_{0,\widehat{\mathcal T}, \widehat \beta}$ (see the proof of Theorem~\ref{thm:approx2}). Therefore,
\begin{align*}
	\Vert D \widehat B\rVert_2^2\le \sigma_{\max}^2(D)\Vert \widehat \beta\rVert_2^2 \le \sigma_{\max}^2(D) \widehat K \Vert \widehat \beta\rVert_\infty^2 \lesssim \widehat K\log n,
\end{align*}	
as soon as $\lVert f_0\rVert_\infty\lesssim \sqrt{\log n }$.
\end{myproof}

\begin{mylemma}[Metric entropy] 
	Let $\bar K_n\asymp n\epsilon_n^2/\log n$ and $\bar s_n\asymp n\epsilon_n^2/\log p$. 
	Define $\mathcal F_{\bar s_n,\bar K_n,M}^{(1)}$ for $M>0$ as in \eqref{eqn:fundef}. Under \ref{asm:spnet0}, for any $C>0$,
$$
\log N\left(\epsilon_n,\mathcal F_{\bar s_n,\bar K_n,n^{C}}^{(1)},\lVert \cdot \rVert_n\right)\lesssim n\epsilon_n^2.
$$
\label{lmm:cond2}
\end{mylemma}

\begin{myproof}
	Observe that the exponential of the left-hand side is bounded by
\begin{align}
	\begin{split}
		 \sum_{S:|S|\le\bar s_n} \sum_{(K^1,\dots,K^T):K^t\le \bar K_n,t=1,\dots,T} \sum_{\mathcal E \in \mathscr E_{S,K^1,\dots,K^T} } N\left(\epsilon_n,\mathcal F_{\mathcal E,n^{C}}^{(1)},\lVert \cdot \rVert_\infty\right).
		\label{eqn:convering}
	\end{split}
\end{align}
For any given $\mathcal E$ and $B_1,B_2\in\mathbb R^{\sum_{t=1}^T K^t}$, 
\begin{align*}
	\lVert f_{\mathcal E,B_1}-f_{\mathcal E,B_2} \rVert_\infty =\sup_{x\in[0,1]^p}\left|\sum_{t=1}^T \sum_{k=1}^{K^t} (\beta_{1k}^t-\beta_{2k}^t) \mathbbm 1(x\in\Omega_k^t)\right|\le \left(\sum_{t=1}^T K^t\right)\lVert B_1-B_2\rVert_\infty.
\end{align*}
Observe that  the cardinality of the set $\mathscr E_{S,\bar K^1,\dots,\bar K^T}$ is equal to $\prod_{t=1}^T |\mathscr T_{S, K^t,\mathcal Z}| \le |\mathscr T_{S, \bar K_n,\mathcal Z}|^T$. Hence, \eqref{eqn:convering} is further bounded by
\begin{align}
	(\bar K_n)^T\times N\left(\frac{\epsilon_n}{T\bar K_n},\left\{ B \in\mathbb R^{T\bar K_n}:\lVert  B \rVert_\infty\le n^{C}\right\},\lVert \cdot \rVert_\infty\right) 	  \sum_{S:|S|\le \bar s_n}|\mathscr T_{S,\bar K_n,\mathcal Z}|^T.
	\label{eqn:entropybound}
\end{align}
Observe that $|\mathscr T_{S,\bar K_n,\mathcal Z}|\le (|S| \max_{1\le j\le p} b_j)^{\bar K_n}$, as all splits are restricted to $S$ and each one has at most $\max_{1\le j\le p} b_j$ split points. It follows that
\begin{align*}
	\sum_{S:|S|\le \bar s_n}|\mathscr T_{S,\bar K_n,\mathcal Z}|^T
	\le \sum_{s=1}^{\bar s_n}\binom{p}{s} \left(s  \max_{1\le j\le p} b_j\right)^{T\bar K_n}\le \bar s_n p^{\bar s_n} \left(\bar s_n  \max_{1\le j\le p} b_j\right)^{T\bar K_n}.
\end{align*}
Therefore, \eqref{eqn:entropybound} is further bounded by
$ (\bar K_n)^T s_n p^{\bar s_n}(\bar s_n \max_{1\le j\le p} b_j)^{T\bar K_n}  ({3 T  \bar K_n n^{C}}/{\epsilon_n})^{T \bar K_n}$.
The logarithm is bounded by a constant multiple of $\bar s_n\log p+\bar K_n\log n\lesssim n\epsilon_n^2$ as soon as $\max_{1\le j\le p} b_j\lesssim \log n$. 
\end{myproof}

\begin{mylemma}[Prior mass of sieve] Let $\bar K_n=\lfloor M' n\epsilon_n^2/\log n \rfloor$, and $\bar s_n=\lfloor M' n\epsilon_n^2/\log p \rfloor$ for a sufficiently large $M'>0$.
	Define $\mathcal F_{\bar s_n,\bar K_n,M}^{(1)}$ for $M>0$ as in \eqref{eqn:fundef}.  Under \ref{pri:tree} and \ref{pri:normal}, for any $C>1$ and $C'>0$,
$$
\Pi(\mathcal F_\ast\setminus \mathcal F_{\bar s_n,\bar K_n, n^{C}}^{(1)})\ll e^{-C' n\epsilon_n^2 }.
$$
\label{lmm:cond3}
\end{mylemma}

\begin{myproof}
Note that $\mathcal F_\ast\setminus \mathcal F_{\bar s_n,\bar K_n, n^{C}}^{(1)} = \mathcal F_{\bar s_n,\bar K_n, n^{C}}^{(2)} \cup (\mathcal F_\ast\setminus (\mathcal F_{\bar s_n,\bar K_n, n^{M_1}}^{(1)}\cup \mathcal F_{\bar s_n,\bar K_n, n^{M_1}}^{(2)}))$. We will give a union bound. First, observe that
\begin{align*}
	\Pi(\mathcal F_{\mathcal E,M}^{(2)}) & = \Pi(B\in\mathbb R^{\sum_{t=1}^T K^t}\! \!: \lVert B\rVert_\infty >M)\le  \Pi\Bigg(B\in\mathbb R^{\sum_{t=1}^T K^t} \! \!: \lVert D B\rVert_\infty >\frac{M\sigma_{\max}^{-1}(D^{-1})}{\sqrt{\sum_{t=1}^T K^t}}\Bigg),
\end{align*}
where $D$ is the matrix with bounded singular values that makes the prior for $DB$ the standard normal distribution.
	Using the tail probability of normal distributions,
\begin{align*}
	\Pi\big(\mathcal F_{\bar s_n,\bar K_n, n^{C}}^{(2)}\big)
	&\le  \sum_{S:|S|\le\bar s_n} \sum_{(K^1,\dots,K^T):K^t\le \bar K_n,t=1,\dots,T} \sum_{\mathcal E \in \mathscr E_{S,K^1,\dots,K^T} } \Pi\big(\mathcal F_{\mathcal E,n^C}^{(2)}\big)\\
	&\le  (\bar K_n)^T \bar s_n p^{\bar s_n} \left(\bar s_n  \max_{1\le j\le p} b_j\right)^{T\bar K_n} 2T\bar K_n e^{-\sigma_{\max}^{-2}(D^{-1})n^{2{C}}/{(2T\bar K_n)}}.
\end{align*}
Since $T\bar K_n\lesssim n\epsilon_n^2/\log n\ll n $ and $\sigma_{\max}(D^{-1})$ is bounded, if $\max_{1\le j\le p} b_j\lesssim \log n$ and ${C}>1$, the right most side of the expression is $o(e^{-C' n\epsilon_n^2 })$ for any $C'>0$. Now observe that
\begin{align}
	\begin{split}
		&\Pi\Big(\mathcal F_\ast\setminus (\mathcal F_{\bar s_n,\bar K_n, n^{M_1}}^{(1)}\cup \mathcal F_{\bar s_n,\bar K_n, n^{M_1}}^{(2)})\Big)\\
		&\quad\le \sum_{t=1}^T\Pi(K^t>\bar K_n)+\Pi(S:s > \bar s_n|K^t\le\bar K_n,t=1,\dots, T).
		\label{eqn:priorsieve}
	\end{split}
\end{align}
The prior satisfies $\log \Pi(K^t>\bar K_n)\lesssim -\bar K_n\log \bar K_n$ for every $t=1,\dots,T$ (see Lemma 5.1 and Corollary 5.2 of \citet{rockova2019theory}). Using that $\bar K_n\asymp n\epsilon_n^2/\log n$ and $n\epsilon_n^2 \gtrsim n^{d/(2\bar\alpha+d)}\ge n^{1/3}$,
we obtain $	-\bar K_n\log \bar K_n \lesssim -\bar K_n \log n$.
To bound the second term of the right-hand side of \eqref{eqn:priorsieve}, we define the set
\begin{align*}
	V_2=\left\{\eta\in\mathbb S^p: \min_{S:|S|=\bar s_n}\sum_{j\notin S} \eta_j \ge \kappa_n\right\},
\end{align*}
for $\kappa_n$ specified below.
By \eqref{eqn:dirineq2} of Lemma~\ref{lmm:dir}, we show that the prior satisfies $\Pi(V_2)\le e^{-C_1(\xi-1)\bar s_n\log p-\log \kappa_n}$ for some $C_1>0$. Hence,
\begin{align*}
	\Pi(S:s > \bar s_n|K^t\le\bar K_n,t=1,\dots, T)&\le e^{-C_1(\xi-1)\bar s_n\log p-  \log \kappa_n}\\
	&\quad + \Pi(S:s > \bar s_n|K^t\le\bar K_n,t=1,\dots, T,V_2^c).
\end{align*}
The term $\Pi(S:s > \bar s_n|K^t\le\bar K_n,t=1,\dots, T,V_2^c)$ is interpreted as the prior probability that splits occur along more than $\bar s_n$ coordinates with at most $T \bar K_n$ splits given $V_2^c$. If $\eta$ is available, this probability is
\begin{align*}
	1-\sum_{S: |S|\le\bar s_n}\left(\sum_{j\in S} \eta_j \right)^{T\bar K_n}\le 1-\left(\max_{S: |S|\le\bar s_n}\sum_{j\in S} \eta_j \right)^{T\bar K_n}.
\end{align*}
Conditional on $V_2^c$, the last expression is further bounded by $1-(1-\kappa_n)^{T\bar K_n}\le \kappa_n T\bar K_n$. Choosing $\kappa_n=e^{-(C'+1)n\epsilon_n^2}$, the resulting bound for \eqref{eqn:priorsieve} gives $\Pi(\mathcal F_\ast\setminus \mathcal F_{\bar s_n,\bar K_n, n^{C}}^{(1)})\ll e^{-C' n\epsilon_n^2 }$ as $M'$ is sufficiently large.
\end{myproof}

\subsection{Proof of Theorem~\ref{thm:minimax}}
\label{sec:proofthm3}
Our proof is similar to the proof of Theorem 3.1 in \citet{yang2015minimax}, which is based on the Le Cam equation  \citep{birge1993rates,wong1995probability,barron1999risk}.
A minimax lower bound of nonparametric regression can be obtained by solving the Le Cam equation with the metric entropy of the target function space \citep{yang1999information}.
We first formalize this result in the following lemma, which is a corollary induced by Theorem 1 of \citet{yang1999information}.

\begin{mylemma}[Minimax lower bound in nonparametric regression]
	For a function space $\mathcal F\subset \mathcal L_2(Q)$, suppose there are upper and lower bounds of the metric entropies as
	\begin{align}
		\begin{split}
	\log N(\epsilon,\mathcal F,\lVert\cdot\rVert_{2,Q})&\le V^\ast(\epsilon),\\
\log D(\epsilon,\mathcal F,\lVert\cdot\rVert_{2,Q})&\ge V_\ast(\epsilon).
		\end{split}
	\label{eqn:minimaxentropybounds}
	\end{align}
	Suppose that $\bar\gamma_n$ is the solution to $V^\ast(\bar\gamma_n)\asymp n\bar\gamma_n^2 $. Then, for the nonparametric regression model in \eqref{eqn:modelregrd}, the sequence $\gamma_n$ such that $V_\ast(\gamma_n)\asymp n\bar\gamma_n^2$ satisfies 
	\begin{align*}
		r_n(\mathcal F,Q)\gtrsim \gamma_n,
	\end{align*}
where $r_n$ is the $L_2(Q)$-minimax risk defined as $r_n(\mathcal F,Q)= \inf_{\hat f \in\mathcal B_n} \sup_{f_0\in \mathcal F}\mathbb E_{f_0,Q}\lVert \hat f-f_0 \rVert_{2,Q}^2$ with $\mathcal B_n$ the space of all $L_2(Q)$-measurable function estimators.
	\label{lmm:yangbarron}
\end{mylemma}

\begin{myproof}
	By Theorem 1 of \citet{yang1999information}, the assertion holds for every statistical model if \eqref{eqn:minimaxentropybounds} is replaced by
	\begin{align*}
		\log N(\epsilon,\mathcal F,K^{1/2})&\le V^\ast(\epsilon),\\
		\log D(\epsilon,\mathcal F,\lVert\cdot\rVert_{2,Q})&\ge V_\ast(\epsilon),
\end{align*}
for the KL divergence $K$.
Let $p_f(x,y)=(2\pi\sigma_0^2)^{-1/2} \exp\{-(y-f(x))^2/(2\sigma_0^2)\} q(x)$. One can easily observe that $	K(p_{f_1},p_{f_2}) =(2\sigma_0^2)^{-1} \lVert f_1-f_2\rVert_{2,Q}^2$. The assertion in the lemma follows immediately.
\end{myproof}

The key to obtaining a sharp minimax lower bound $\gamma_n$ is to establish the bounds $ V^\ast$ and $ V_\ast$ as tight as possible. In the Lemmas~\ref{lmm:entropy}--\ref{lmm:funcex} below, we provide entropy estimates for the $d$-dimensional (non-sparse) piecewise heterogeneous anisotropic H\"older space.
While an upper bound of the metric entropy is well known for isotropic classes (e.g., Theorem~2.7.1 of \citet{van1996weak}), we believe that there is no available result on more complicated function space in the literature, even for the simple anisotropic classes in Definition~\ref{def:hol}. Lemma~\ref{lmm:packbound} concatenates the results to obtain entropy bounds for the sparse function space.

Below we write $\overline{\mathcal H}_{\lambda,M}^{A_{\bar\alpha},d}(\mathfrak X_0)=\{h\in {\mathcal H}_\lambda^{A_{\bar\alpha},d}(\mathfrak X_0):\lVert h\rVert_\infty\le M\lambda \}$ for $M>0$. For the upper bound of the metric entropy, we consider a bound for the space $\overline{\mathcal H}_{1,M}^{A_{\bar\alpha},d}(\mathfrak X_0)$, which is not worse than that for $\overline{\mathcal H}_{1,M}^{A_{\bar\alpha},d}(\mathfrak X_0)\cap {\mathcal C}([0,1]^d)$. This implies that the Le Cam equation gives the same minimax lower bound for the two spaces.

\begin{mylemma}[Covering number, upper bound]
	For $d>0$, $R>0$, a partition $\mathfrak X_0=\{\Xi_1,\dots,\Xi_R\}$ of $[0,1]^d$, and a smoothness parameter $A_{\bar\alpha}\in \mathcal A_{\bar\alpha}^{R,d}$ for $\bar\alpha\in(0,1]$ such that
	$\log\mathsf{len}([\Xi_r]_j)\gtrsim -1/\alpha_{rj}$, $1\le r\le R$, $1\le j\le d$, there exist constants $\epsilon_0>0$ and $M_0>1$ such that for any $\epsilon<\epsilon_0$,
	\begin{align}
		\label{eqn:entropy2}
		\log N\big(\epsilon,\overline{\mathcal H}_{1,M}^{A_{\bar\alpha},d}(\mathfrak X_0),\lVert\cdot\rVert_\infty\big) &\le  (M_0 d/\epsilon)^{d/\bar\alpha}.
	\end{align}
	\label{lmm:entropy}
\end{mylemma}

\begin{myproof}
	To express the assumption more explicitly, let $C_1>0$ be a constant such that $\log\mathsf{len}([\Xi_r]_j)\ge -C_1/\alpha_{rj}$ for every $r$ and $j$.
	For a sufficiently small $C_2>0$, choose $\delta_d\in(0,\min\{e^{-C_1},C_2/d\})$ such that $\min_{r,j}\mathsf{len}([\Xi_r]_j)\delta_d^{-1/\alpha_{rj}}>1$. 
	On each box $\mathsf{cl}(\Xi_r)$, consider a Cartesian product of grid points,
	\begin{align*}
		\tilde{\mathcal G}_r:=\prod_{j=1}^d \left \{ I_{rj}^L, I_{rj}^L+u_{rj} ,I_{rj}^L+2u_{rj},\dots,  I_{rj}^L+\mathsf{len}([\Xi_r]_j)\right\},
	\end{align*}
	where $u_{rj}=\mathsf{len}([\Xi_r]_j)/\lceil \mathsf{len}([\Xi_r]_j)\delta_d^{-1/\alpha_{rj}}\rceil$ is the mesh-size and $I_{rj}^L$ is the left-boundary of $\Xi_r$ in coordinate $j$. 
	Observe that
	\begin{align}
		\tilde m_r := |\tilde{\mathcal G}_r| =\prod_{j=1}^d (1+\lceil \mathsf{len}([\Xi_r]_j)  \delta_d^{-1/\alpha_{rj}}\rceil)\le\prod_{j=1}^d (2+ \mathsf{len}([\Xi_r]_j)  \delta_d^{-1/\alpha_{rj}})\le \mathsf{vol}(\Xi_r) 3^d\delta_d^{-d/\bar\alpha}.
		\label{eqn:mbo1}
	\end{align}
We write the elements of $\tilde{\mathcal G}_r$ as $x_{r}^\ell=(x_{r1}^\ell,\dots,x_{rd}^\ell)^\top$, i.e., $x_r^\ell\in \tilde{\mathcal G}_r$, $\ell=1,\dots,\tilde m_r$, $r=1,\dots, R$.
	For every $h\in\overline{\mathcal H}_{1,M}^{A_{\bar\alpha},d}(\mathfrak X_0)$, we define the vector
	\begin{align*}
		G h = \left(\lfloor h(x_1^1)/\delta_d \rfloor,\dots,\lfloor h(x_1^{\tilde m_1})/\delta_d \rfloor,\dots,\lfloor h(x_R^1)/\delta_d \rfloor,\dots,\lfloor h(x_R^{\tilde m_R})/\delta_d \rfloor\right)^\top.
	\end{align*}
	Because mesh-size satisfies $u_{rj}\le \delta_d^{1/\alpha_{rj}}$, for every $x=(x_1,\dots, x_d)^\top\in \Xi_r$ with given $r$, there exists a point $x_{r}^\ell\in \tilde{\mathcal G}_r$ such that $\sum_{j=1}^d |x_j-x_{rj}^\ell|^{\alpha_{rj}}\le d\delta_d$. Hence, for every such $x$ and $x_r^\ell$, all functions $h_1,h_2\in\overline{\mathcal H}_{1,M}^{A_{\bar\alpha},d}(\mathfrak X_0)$ such that $Gh_1=Gh_2$ satisfy
	\begin{align*}
		|h_1(x)-h_2(x)| \le |h_1(x_r^\ell)-h_2(x_r^\ell)|+2\sum_{j=1}^d|x_j-x_{rj}^\ell|^{\alpha_{rj}}\le \delta_d+ 2d\delta_d.
	\end{align*}
	As this holds for every $1\le r\le R$, it follows that $\lVert h_1-h_2 \rVert_\infty\le  3 d \delta_d$ for any $h_1,h_2$ such that $Gh_1=Gh_2$.
	This means that, whenever $3 d \delta_d<\epsilon_0$ for some small constant $\epsilon_0>0$, the covering number $N(3 d \delta_d,\overline{\mathcal H}_{1,M}^{A_{\bar\alpha},d}(\mathfrak X_0),\lVert\cdot\rVert_\infty)$ is bounded by the number of possible vectors $Gh$ for $h$ that ranges over $\overline{\mathcal H}_{1,M}^{A_{\bar\alpha},d}(\mathfrak X_0)$.
	
	Without loss of generality, we now assume that $(x_r^\ell)_{\ell=1}^{\tilde m_r}$ in $\tilde{\mathcal G}_r$ are appropriately sorted so that every two successive values differ in only one coordinate by mesh-size; that is, for every $\ell>1$, there exists $\ell' < \ell$ such that $\sum_{j=1}^d |x_{rj}^{\ell'}-x_{rj}^\ell|^{\alpha_{rj}}=u_{rj'}^{\alpha_{rj'}}\le \delta_d$ for some $j'$.
	For the enumeration, we begin with the first element of $Gh$, which is defined with $x_1^1\in\tilde{\mathcal G}_1$. As $\lVert h \rVert_\infty\le M$, the number of possible values of $\lfloor h(x_1^1)/\delta_d\rfloor$ does not exceed $2M/\delta_d+1$. For every remainder defined with $x_1^\ell\in\tilde{\mathcal G}_1$, $2\le\ell\le \tilde m_1$, there exists $\ell' < \ell$ such that
	\begin{align*}
		&|\lfloor h(x_1^{\ell'})/\delta_d \rfloor- \lfloor h(x_1^\ell)/\delta_d\rfloor |\\
		&\quad\le \delta_d^{-1} |h(x_1^{\ell'})-h(x_1^\ell)| + |h(x_1^{\ell'})/\delta_d- \lfloor h(x_1^{\ell'})/\delta_d\rfloor|+|h(x_1^\ell)/\delta_d- \lfloor h(x_1^\ell)/\delta_d\rfloor|\\
		&\quad\le \delta_d^{-1}\sum_{j=1}^d|x_j^{\ell'}-x_j^\ell|^{\alpha_{rj}} + 2\le 3.
	\end{align*}
	It follows that, for a given $\lfloor h(x_1^{\ell'})/\delta_d\rfloor$, the number of possible values of $\lfloor h(x_1^\ell)/\delta_d\rfloor$ is at most $7$, which is the case for every $\ell>1$.
	Putting the bounds together, the number of possible values of the first $m_1$ elements of $Gh$ is bounded by $(2M/\delta_d+1) 7^{\tilde m_1 -1}$. Next, because $h\in\overline{\mathcal H}_{1,M}^{A_{\bar\alpha},d}(\mathfrak X_0)$ can be discontinuous at the boundaries of the pieces of $\mathfrak X_0$, the $(\tilde m_1+1)$th element of $Gh$, defined with $x_2^1\in\tilde{\mathcal G}_2$, has no restriction. Similar to the case with $r=1$ above, the number of possible values of $\lfloor h(x_2^1)/\delta_d\rfloor$ at most $2M/\delta_d+1$, and the number of possible values of $\lfloor h(x_2^\ell)/\delta_d\rfloor$ is at most $7$ for every $2\le\ell\le \tilde m_2$. This concludes that the number of possible values of the next $\tilde m_2$ elements of $Gh$ is bounded by $(2M/\delta_d+1) 7^{\tilde m_2 -1}$. Concatenating this for all $r$,  the number of possible vectors $Gh$ is clearly at most $\prod_{r=1}^R (2M/\delta_d+1) 7^{\tilde m_r -1}=(2M/\delta_d+1)^R 7^{\tilde m -R}$, where $\tilde m=\sum_{r=1}^R \tilde m_r$. Using \eqref{eqn:mbo1}, it is evident that $\tilde m\le 3^d\delta_d^{-d/\bar\alpha}$ because $\sum_{r=1}^R \mathsf{vol}(\Xi_r) =1$.
	Taking $\epsilon = 3d\delta_d$,
	\begin{align*}
		\log N(\epsilon,\overline{\mathcal H}_{1,M}^{A_{\bar\alpha},d}(\mathfrak X_0)\cap {\mathcal C}([0,1]^d),\lVert\cdot\rVert_\infty)\le R\log(6Md/\epsilon+1)+
		3^d(3d/\epsilon)^{d/\bar\alpha}\log 7.
	\end{align*}
	As $ \log(6Md/\epsilon+1)\lesssim (6Md/\epsilon)^{d/\bar\alpha}$ and $\log R\lesssim d/\bar\alpha$ (by the condition $\log\mathsf{len}([\Xi_r]_j)\gtrsim -1/\alpha_{rj}$), the last expression is bounded by $(M_0d/\epsilon)^{d/\bar\alpha}$ for some $M_0>0$. To complete the proof, we must now show that there exists a small constant $\epsilon_0>0$ such that $\epsilon=3d\delta_d<\epsilon_0$. This is achieved by a sufficiently small $C_2>0$ since $\delta_d<C_2/d$.
\end{myproof}

\begin{mylemma}[Packing number, lower bound]
	For $d>0$, $R>0$, a partition $\mathfrak X_0=\{\Xi_1,\dots,\Xi_R\}$ of $[0,1]^d$, and a smoothness parameter $A_{\bar\alpha}\in \mathcal A_{\bar\alpha}^{R,d}$ for $\bar\alpha\in(0,1]$ such that
	$\log\mathsf{len}([\Xi_r]_j)\gtrsim -1/\alpha_{rj}$, $1\le r\le R$, $1\le j\le d$, there exist constants $\epsilon_1>0$ and $M_1>1$ such that for any $\epsilon<\epsilon_1^d$, there are $N\ge \exp\{1/(M_1^d\epsilon)^{d/\bar\alpha}\}$ functions $h_i\in \overline{\mathcal H}_{1,M}^{A_{\bar\alpha},d}(\mathfrak X_0)\cap {\mathcal C}([0,1]^d)$, $i=1,\dots,N$, and $h_0=0$ satisfying
	\begin{align}
		\label{eqn:lmm01}	\int_{[0,1]^d} h_i(x)dx_j&=0,\quad 0\le i\le N,\quad 1\le j\le d, \\
		\label{eqn:lmm02} \lVert h_i -h_k \rVert_2&\ge \epsilon,\quad  0\le i\le k\le N.
	\end{align}
	\label{lmm:funcex}
\end{mylemma}

\begin{myproof}
	Similar to above, let $C_1\ge \log 8$ be a constant such that $\log\mathsf{len}([\Xi_r]_j)\ge -C_1/\alpha_{rj}$ for every $r$ and $j$ and choose a constant $\delta\in(0,\min\{e^{-C_1},M\}]$ such that $\mathsf{len}([\Xi_r]_j)\delta^{-1/\alpha_{rj}}>1$, $ 1\le r\le R$, $1\le j\le d$.
	On each box $\Xi_r$, consider a Cartesian product of grid points,
	\begin{align*}
\mathcal G_r:=\prod_{j=1}^d \left\{I_{rj}^L+\frac{u_{rj}}{2}, I_{rj}^L+\frac{3u_{rj}}{2} ,I_{rj}^L+\frac{5u_{rj}}{2},\dots,  I_{rj}^L+\mathsf{len}([\Xi_r]_j)-\frac{u_{rj}}{2}\right\},
	\end{align*}
	where $u_{rj}=\mathsf{len}([\Xi_r]_j)/\lceil \mathsf{len}([\Xi_r]_j)\delta^{-1/\alpha_{rj}}\rceil$ is the mesh-size and $I_{rj}^L$ is the left-boundary of $\Xi_r$ in coordinate $j$ (cf. the grid $\tilde{\mathcal G}_r$ used in the proof of  Lemma~\ref{lmm:entropy}).
	Note that
	\begin{align}
	m_r:=|\mathcal G_r|=\prod_{j=1}^d \lceil \mathsf{len}([\Xi_r]_j)  \delta^{-1/\alpha_{rj}}\rceil\ge\mathsf{vol}(\Xi_r)\delta^{-d/\bar\alpha}.
	\label{eqn:mbo2}
	\end{align}
We write the elements of ${\mathcal G}_r$ as $x_{r}^\ell=(x_{r1}^\ell,\dots,x_{rd}^\ell)^\top$, i.e., $x_r^\ell\in {\mathcal G}_r$, $\ell=1,\dots, m_r$, $r=1,\dots, R$.
We define the univariate kernel $\mathcal K(t)= t \mathbbm 1(|t|\le 1/2) + ({\rm sgn}(t)-t)\mathbbm 1 (1/2<|t|\le 1)$, $t\in\mathbb R$, supported on $[-1,1]$. Clearly, $\mathcal K$ is 1-Lipschitz and satisfies $\int \mathcal K(t)dt =0$.

We define the function
\begin{align*}
	\phi_r^\ell(x)=\frac{\delta}{2^{d+1}}\prod_{j=1}^d\mathcal K\left(\frac{x_j- x_{rj}^\ell}{u_{rj}/2}\right), \quad 1\le\ell\le m_r, \quad  1\le r\le R,
\end{align*}
which is supported on $\mathcal X_r^\ell\coloneqq\prod_{j=1}^d[x_j^\ell-u_{rj}/2,x_j^\ell+u_{rj}/2]$ with the center $x_r^\ell$. 
As $\lVert\mathcal K\rVert_\infty = 1/2$, we obtain $\lVert \phi_r^\ell \rVert_\infty \le \delta\lVert\mathcal K\rVert_\infty^d /2^{d+1} \le 1$ for a suitable $C_1>0$.
Using  the  Lipschitz continuity of $\mathcal K$ and the inequality $|\prod_j a_j-\prod_j b_j|\le \sum_j | a_j-b_j | $ for any $a_j,b_j\in[-1,1]$,
we have that for any $x,y$ on the support $\mathcal X_r^\ell$,
\begin{align*}
	|\phi_r^\ell(x)-\phi_r^\ell(y)| \le \frac{\delta}{2}\sum_{j=1}^d \left|\frac{x_j-y_j}{u_{rj}}\right|\le 	\frac{\delta}{2}\sum_{j=1}^d \left|\frac{x_j-y_j}{u_{rj}}\right|^{\alpha_{rj}}\le \sum_{j=1}^d |x_j-y_j|^{\alpha_{rj}},
\end{align*}	
	where we used  the inequalities $x\le x^a$ for any $x\in[0,1]$ and $a\in[0,1]$, and $u_{rj}\ge 1/ (2\delta^{-1/\alpha_{rj}})$ as soon as $\mathsf{len}([\Xi_r]_j)\delta^{-1/\alpha_{rj}}\ge 1/2$ (note that $\lceil x \rceil\le 2x$ for $x\ge 1/2$).
	This shows that $\phi_r^\ell\in\mathcal H_1^{\alpha_r,d}(\mathcal X_r^\ell)$ for every $1\le\ell\le m_r$ and $1\le r\le R$.
For a binary vector $\tilde\omega_r=(\tilde\omega_r^1,\dots,\tilde\omega_r^{m_r})^\top\in\{0,1\}^{m_r}$, define the continuous function $h_{\tilde\omega_r}=\sum_{\ell=1}^{m_r} \tilde\omega_r^\ell \phi_r^\ell$ supported on $\Xi_r$.
As $\int \phi_r^\ell (x) dx_j=0$ for every $j$ and each $\phi_r^\ell$ is a shifted copy of another, we obtain $\int h_{\tilde\omega_r}(x)dx_j=0$ for every $j$ and $h_{\tilde\omega_r}\in\mathcal H_1^{\alpha_r,d}(\Xi_r)$.
Let $m=\sum_{r=1}^R m_r$, which satisfies $m\ge \delta^{-d/\bar\alpha}$ by \eqref{eqn:mbo2}. We write $\omega=(\omega_1,\dots,\omega_m)^\top=(\tilde\omega_1^\top,\dots,\tilde\omega_R^\top)^\top\in\{0,1\}^{m}$ and define $h_{\omega}=\sum_{r=1}^R h_{\tilde\omega_r}$.
Then, as $\lVert h_{\omega} \rVert_\infty = 	\max_{r,\ell}\lVert\phi_r^\ell\rVert_\infty\le \delta\le M$ and each $h_{\tilde\omega_r}$ is zero at all points on the boundary of $\Xi_r$, it is easy to see that
$h_{\omega}\in \overline{\mathcal H}_{1,M}^{A_{\bar\alpha},d}(\mathfrak X_0)\cap {\mathcal C}([0,1]^d)$ and $\int h_{\omega} (x) dx_j=0$. We also have that for any $\omega,\omega'\in\{0,1\}^m$,
\begin{align}
	\lVert h_{\omega} - h_{\omega'}\rVert_2^2 &\ge \left[\sum_{b=1}^m (\omega_b-\omega_b')^2\right] \min_{r,\ell}\int [\phi_r^ {\ell}(x)]^2dx=  \rho(\omega,\omega')\left(\frac{\delta^2\lVert \mathcal K \rVert_2^{2d}}{2^{3d+2}}\right) \min_{r}\prod_{j=1}^d u_{rj},
	\label{eqn:l2bo}
\end{align}	
where $\rho(\omega,\omega')=\sum_{b=1}^m \mathbbm 1(\omega_b\ne \omega_b')$ is the Hamming distance between $\omega$ and $\omega'$. As $m\ge\delta^{-d/\bar\alpha}\ge \delta^{-1}\ge e^{C_1}> 8$, the Gilbert-Varshamov bound (Lemma~2.9 of \citet{tsybakov2008introduction}) says that
there exist $N\ge 2^{m/8}$ binary strings $\omega^{(1)},\dots,\omega^{(N)}\in\{0,1\}^m$ such that $\rho(\omega^{(\ell)},\omega^{(\ell')})\ge m/8$, $0\le\ell<\ell'\le N$, with $\omega^{(0)}=0$. 
As $\min_{r}\prod_{j=1}^d u_{rj}\ge 1/ (2^d\delta^{-d/\bar\alpha})$ and $\lVert\mathcal K\rVert_2^2 = 1/6$, the lower bound in \eqref{eqn:l2bo} gives that for every $0\le\ell<\ell'\le N$,
\begin{align*}
	\lVert h_{\omega^{(\ell)}} - h_{\omega^{(\ell')}}\rVert_2^2 \ge\frac{m}{8}\left(\frac{\delta^2}{6^d 2^{3d+2}}\right) \min_{r}\prod_{j=1}^d u_{rj}\ge \frac{\delta^2}{2^5 96^d}.
\end{align*}	
Letting $\epsilon=\delta/\sqrt{2^5 96^d}$, the previous lower bound gives $\lVert h_{\omega^{(\ell)}} - h_{\omega^{(\ell')}}\rVert_2\ge\epsilon$ while $N\ge2 ^{m/8}\ge \exp(\delta^{-d/\bar\alpha}(\log 2)/8)\ge\exp(1/(2^8 96^{d/2}\epsilon)^{d/\bar\alpha})$. As $\delta$ is a constant, this holds for every $\epsilon<\epsilon_1/96^{d/2}$ for some $\epsilon_1>0$.
\end{myproof}

\begin{mylemma}[Entropy with sparsity]
	For $d>0$, $\lambda>0$, $R>0$, a partition $\mathfrak X_0=\{\Xi_1,\dots,\Xi_R\}$ of $[0,1]^d$, and a smoothness parameter $A_{\bar\alpha}\in \mathcal A_{\bar\alpha}^{R,d}$ for $\bar\alpha\in(0,1]$ such that, there exist $\epsilon_{2}>0$ and $M_2>1$ such that for any $\epsilon<\epsilon_{2}$ and $\epsilon'<\epsilon_{2}^d$,
	\begin{align}
		\label{eqn:packbound3}
		\log N \big(\epsilon,  \overline{\Gamma}_{\lambda,M}^{A_{\bar\alpha},d,p}(\mathfrak X_0), \lVert \cdot\rVert_2\big)&\le \log \binom{p}{d}+\Big(\frac{M_2  \lambda d}{\epsilon}\Big)^{d/\bar\alpha},\\
		\label{eqn:packbound1}
		\log D \big(\epsilon',  \overline{\Gamma}_{\lambda,M}^{A_{\bar\alpha},d,p}(\mathfrak X_0)\cap {\mathcal C}([0,1]^p), \lVert \cdot\rVert_2\big)&\ge \log \binom{p}{d} +\Big(\frac{\lambda}{M_2^d\epsilon'}\Big)^{d/\bar\alpha}.
	\end{align}
	\label{lmm:packbound}
\end{mylemma}

\begin{myproof}
	We only need to verify the assertion for $\lambda=1$ since $D(\epsilon,\lambda\mathcal F,\lVert\cdot\rVert_2)=D(\epsilon/\lambda,\mathcal F,\lVert\cdot\rVert_2)$ and $N(\epsilon,\lambda\mathcal F,\lVert\cdot\rVert_2)=N(\epsilon/\lambda,\mathcal F,\lVert\cdot\rVert_2)$ for any set $\mathcal F$.
	We first verify the upper bound \eqref{eqn:packbound3}.
	For every $\epsilon<\epsilon_0$, Lemma~\ref{lmm:entropy} gives 
	$\log N (\epsilon, \overline{\mathcal H}_{1,M}^{A_{\bar\alpha},d}(\mathfrak X_0), \lVert \cdot\rVert_2)\le (M_0d/\epsilon)^{d/\bar\alpha}$.
	As $\overline{\Gamma}_{1,M}^{A_{\bar\alpha},d,p}(\mathfrak X_0)$ is a union of $\binom{p}{d}$ many $\overline{\mathcal H}_{1,M}^{A_{\bar\alpha},d}(\mathfrak X_0)$, the assertion easily follows.
	
	Next, we verify \eqref{eqn:packbound1}.
	By Lemma~\ref{lmm:funcex}, for every $\epsilon'<\epsilon_1^d$, there are functions $h_0=0$, $h_i\in\overline{\mathcal H}_{1,M}^{A_{\bar\alpha},d}(\mathfrak X_0)\cap {\mathcal C}([0,1]^d)$, $1\le i\le N$ satisfying \eqref{eqn:lmm01} and \eqref{eqn:lmm02}, with $N\ge \exp\{1/(M_1^d\epsilon')^{d/\bar\alpha}\}$. This means that for any $S\subseteq\{1,\dots,p\}$ such that $|S|=d$, we have that $W_S^p h_i\in\overline{\Gamma}_{1,M}^{A_{\bar\alpha},d,p}(\mathfrak X_0)\cap {\mathcal C}([0,1]^p)$ for every such $h_i$ , $0\le i\le N$. Therefore,
	\begin{align*}
		\mathcal W(\epsilon'):=\bigcup_{S\subseteq\{1,\dots,p\}:|S|=d} \{ W_S^p h_i:1\le i\le N\}\subseteq \overline{\Gamma}_{1,M}^{A_{\bar\alpha},d,p}(\mathfrak X_0)\cap {\mathcal C}([0,1]^p).
	\end{align*}
	Now, for any $S\ne S'\subseteq\{1,\dots,p\}$ and $1\le i\le k \le N$, observe that $\lVert W_S^p h_i -  W_{S'}^p h_k \rVert_2=(\lVert h_i \rVert_2^2 + \lVert h_k \rVert_2^2)^{1/2}\ge \epsilon' $ by \eqref{eqn:lmm02}, as $\langle W_S^p h_i,  W_{S'}^p h_k \rangle=0$ owing to \eqref{eqn:lmm01}, where we used $h_0=0$. 
	Also for any $S\subseteq\{1,\dots,p\}$, it is easy to see that $\lVert W_S^p h_i - W_S^p h_k\rVert_2=\lVert h_i - h_k\rVert_2\ge \epsilon'$ by \eqref{eqn:lmm02}. These imply that $\mathcal W(\epsilon')$ is $\epsilon'$-separated, and hence the packing number $ D (\epsilon', \overline{\Gamma}_{1,M}^{A_{\bar\alpha},d,p}(\mathfrak X_0)\cap {\mathcal C} ([0,1]^p), \lVert \cdot\rVert_2)$ is bounded below by the cardinality of $\mathcal W(\epsilon')$, which is $\binom{p}{d} N$. This leads to the assertion.
\end{myproof}

\begin{myprooftitle}[Proof of Theorem~\ref{thm:minimax}]
	Let the right-hand sides of \eqref{eqn:packbound3} and \eqref{eqn:packbound1} be $V^\ast(\epsilon)$ and $V_\ast(\epsilon)$, respectively.
As $L_2(Q)$-norm can be replaced by $L_2$-norm under Assumption \ref{asm:density},
Lemma~\ref{lmm:yangbarron} implies that a sequence $\gamma_n$ is a minimax lower bound if $V_\ast(\gamma_n)= n\bar\gamma_n^2$ and $V^\ast(\bar\gamma_n)= n\bar\gamma_n^2$ for some $\bar\gamma_n$.  

Let $\hat\gamma_n= \sqrt{n^{-1}\log\binom{p}{d}}+((\lambda d)^{d/\bar\alpha}/n)^{\bar\alpha/(2\bar\alpha+d)}$ and $\bar\gamma_n$ be the solution to $V^\ast(\bar\gamma_n)=n\bar\gamma_n^2$.
As $V^\ast(\epsilon)$ is nondecreasing in $\epsilon$, we obtain
\begin{align*}
	V^\ast(M_2\hat\gamma_n)\le V^\ast\Big(M_2((\lambda d)^{d/\bar\alpha}/n)^{\bar\alpha/(2\bar\alpha+d)}\Big) = n\hat\gamma_n^2 \le M_2^2 n\hat\gamma_n^2.
\end{align*}
This shows that
$\bar\gamma_n\le M_2\hat\gamma_n$. Now, define
	$\kappa_n=\max\left\{\sqrt{n^{-1}\log\binom{p}{d}},((\lambda d)^{d/\bar\alpha}/n)^{\bar\alpha/(2\bar\alpha+d)}\right\}$. It follows that $V^\ast(\hat\gamma_n/2)\ge V^\ast(\kappa_n)$ because $\hat\gamma_n/2\le \kappa_n$. 
	If $\sqrt{n^{-1}\log\binom{p}{d}} \le ((\lambda d)^{d/\bar\alpha}/n)^{\bar\alpha/(2\bar\alpha+d)}$,
\begin{align*}
V^\ast(\kappa_n)=V^\ast\Big(((\lambda d)^{d/\bar\alpha}/n)^{\bar\alpha/(2\bar\alpha+d)}\Big)\ge M_2^{d/\bar\alpha}n\kappa_n^2\ge n\hat\gamma_n^2/4,
\end{align*}
as $M_2^{d/\bar\alpha}\ge 1$. 	If $\sqrt{n^{-1}\log\binom{p}{d}} > ((\lambda d)^{d/\bar\alpha}/n)^{\bar\alpha/(2\bar\alpha+d)}$,
\begin{align*}
 V^\ast(\kappa_n)=V^\ast\bigg(\sqrt{n^{-1}\log\binom{p}{d}}\bigg)\ge n\kappa_n^2\ge n\hat\gamma_n^2/4.
\end{align*}
Putting the bounds together, we obtain $\bar\gamma_n\ge \hat\gamma_n/2$. This concludes $\bar\gamma_n\asymp \hat\gamma_n$.

Now, let $\tilde \gamma_n=\sqrt{n^{-1}\log\binom{p}{d}}+M_2^{-d}( \lambda^{d/\bar\alpha}/(d^2n))^{\bar\alpha/(2\bar\alpha+d)}$ and $\gamma_n$ be the solution to $V_\ast(\gamma_n)= n\hat\gamma_n^2$. Then, it is easy to see that
\begin{align*}
	V_\ast(\tilde\gamma_n)\le V_\ast\Big(M_2^{-d}( \lambda^{d/\bar\alpha}/(d^2n))^{\bar\alpha/(2\bar\alpha+d)}\Big) = n\hat\gamma_n^2,
\end{align*}
which implies $\gamma_n\le \tilde\gamma_n$. Let
	$\tilde\kappa_n=\max\left\{\sqrt{n^{-1}\log\binom{p}{d}},M_2^{-d}( \lambda^{d/\bar\alpha}/(d^2n))^{\bar\alpha/(2\bar\alpha+d)}\right\}$ and note that $V_\ast(\tilde\gamma_n/2)\ge V_\ast(\tilde\kappa_n)$. Similar to the above, if $\sqrt{n^{-1}\log\binom{p}{d}}\le M_2^{-d}( \lambda^{d/\bar\alpha}/(d^2n))^{\bar\alpha/(2\bar\alpha+d)}$,
\begin{align*}
	V_\ast(\tilde\kappa_n) = V_\ast\Big(M_2^{-d}( \lambda^{d/\bar\alpha}/(d^2n))^{\bar\alpha/(2\bar\alpha+d)}\Big) = n\hat\gamma_n^2,
\end{align*}
 and if $\sqrt{n^{-1}\log\binom{p}{d}} > M_2^{-d}( \lambda^{d/\bar\alpha}/(d^2n))^{\bar\alpha/(2\bar\alpha+d)}$,
\begin{align*}
	V_\ast(\tilde\kappa_n) =V_\ast\bigg(\sqrt{n^{-1}\log\binom{p}{d}}\bigg)\ge n\tilde \kappa_n^2\ge n\hat\gamma_n^2/4. 
\end{align*}
These give $\gamma_n\ge \tilde \gamma_n/2$, and hence $\gamma_n \asymp \tilde \gamma_n$. Lemma~\ref{lmm:yangbarron} concludes that $	r_n \big(\overline{\Gamma}^{A_{\bar\alpha},d,p}_{\lambda,M}(\mathfrak X_0)\cap {\mathcal C}([0,1]^p)\big)\gtrsim\tilde\gamma_n$. As $M_2^{-d}d^{-2\alpha/(2\alpha+d)}\ge M_2^{-d}d^{-2}$, $\tilde\gamma_n$ is bounded below by the lower bound in Theorem~\ref{thm:minimax} for some $M_d>1$ depending only on $d$.
\end{myprooftitle}

\subsection{Proofs of Theorems~\ref{thm:nonregrd}--\ref{thm:addreg}}
\label{sec:proofthm4-5}
This section provides proofs of Theorems~\ref{thm:nonregrd}--\ref{thm:addreg}. The proofs are largely based on the proof of Theorem~\ref{thm:nonreg}. We often refer to the reader to the proof of Theorem~\ref{thm:nonreg} rather than showing all details.

\bigskip

\begin{myprooftitle}[Proof of Theorem~\ref{thm:nonregrd}]
	Let $p_{f,\sigma^2}$ be the density of model~\eqref{eqn:modelregrd} with $f$ and $\sigma^2$.
	By Lemma~B.1 of \citet{xie2018adaptive}, the Hellinger distance $\rho_{\rm H}$ satisfies
	\begin{align}
		\lVert f_1-f_2\rVert_{2,Q}^2+|\sigma_1^2-\sigma_2^2|^2\lesssim \rho_{\rm H}^2(p_{f_1,\sigma_1^2},p_{f_2,\sigma_2^2})\lesssim \lVert f_1-f_2\rVert_{1,Q}+|\sigma_1^2-\sigma_2^2|^2,
		\label{eqn:hellbound}
	\end{align}
	if $f_1,f_2,\log\sigma_1,\log\sigma_2$ are uniformly bounded (we use variance parameters in place of standard deviations; both are identical up to constants under the boundedness assumption). Hence, it suffices to show the assertion with respect to the Hellinger distance.
	
	By the well-known theory of posterior contraction (e.g., Theorem 2.1 of \citet{ghosal2000convergence}), we need to verify that there exists $\Theta_n\subseteq\mathcal F\times [\overline C_2^{-1},\overline C_2]$ such that for some $\bar c>0$ and a sufficiently large $\bar c'>0$,
	\begin{align}
		\label{eqn:con12}
		\Pi(B_n)&\ge e^{-\bar c n\epsilon_n^2 },\\
		\label{eqn:con22}
		\log N(\epsilon_n, \Theta_n, \rho_{\rm H})&\lesssim n\epsilon_n^2,\\
		\label{eqn:con32}
		\Pi((f,\sigma^2)\notin \Theta_n)&\ll e^{-\bar c' n\epsilon_n^2 },
	\end{align}
	similar to \eqref{eqn:con1}--\eqref{eqn:con3}, where $B_n=\{ f:  K(p_0,p_{f,\sigma^2})\le\epsilon_n^2 , V(p_0,p_{f,\sigma^2})\le \epsilon_n^2\}$.  Using \eqref{eqn:hellbound}, the conditions \eqref{eqn:con22} and \eqref{eqn:con32} can be similarly verified as in the proof of Theorem~\ref{thm:nonreg}; only difference is that we use truncated priors, so \eqref{eqn:con32} is even more easily satisfied.
	For \eqref{eqn:con12}, note that by Lemma~B.2 of \citet{xie2018adaptive},
	\begin{align*}
		\max\left\{ K(p_0,p_{f,\sigma^2}) , V(p_0,p_{f,\sigma^2}) \right\}\lesssim \lVert f-f_0\rVert_{2,Q}^2+|\sigma^2-\sigma_0^2|,
	\end{align*}
	as $\lVert f_0\rVert_\infty$ and $|\log\sigma_0|$ are bounded and the priors are truncated.
	Hence, there exists a constant $C_1>0$ such that
	\begin{align*}
		B_n&\supseteq \{(f,\sigma^2): \lVert f-f_0\rVert_{2,Q} \le C_1 \epsilon_n , | \sigma^2-\sigma_0^2|\le C_1 \epsilon_n^2 \}.
	\end{align*}
	Note that $\lVert f-f_0\rVert_{2,Q}\lesssim \lVert f-f_0\rVert_2$ if the density of $Q$ is bounded. It is easy to see that $\log \Pi(\sigma^2: |\sigma^2-\sigma_0^2|\le C_1\epsilon_n^2)\gtrsim -\log n$, as $|\log\sigma_0^2|$ is bounded.
 Uisng Lemmas~\ref{lmm:priorcon}--\ref{lmm:cond3}, the rest of the proof follows similarly to that of Theorem~\ref{thm:nonreg}.
\end{myprooftitle}

\begin{myprooftitle}[Proof of Theorem~\ref{thm:density}]
	It is well known that the Hellinger distance possesses an exponentially powerful local test with respect to both the type-I and type-II errors (e.g., Section~7 of \citet{ghosal2000convergence} or Lemma~2 of \citet{ghosal2007convergence}). Therefore by the general posterior contraction theory, it suffices to show that there exists $\Theta_n\subseteq \mathcal F$ such that for some $\bar c>0$ and a sufficiently large $\bar c'>0$,
	\begin{align}
		\label{eqn:con1-2}
		\Pi(B_n)&\ge e^{-\bar c n\epsilon_n^2 },\\
		\label{eqn:con2-2}
		\log N(\epsilon_n, \Theta_n, \rho_{\rm H})&\lesssim n\epsilon_n^2,\\
		\label{eqn:con3-2}
		\Pi(f\notin \Theta_n)&\ll e^{-\bar c' n\epsilon_n^2 },
	\end{align}
	where $B_n=\{f: K(p_0,p_f)\le \epsilon_n^2 , V(p_0,p_f)\le \epsilon_n^2\}$. The last condition \eqref{eqn:con3-2} follows directly from the proof of Theorem~\ref{thm:nonreg}, so we only need to verify \eqref{eqn:con1-2} and \eqref{eqn:con2-2}.
	
	By Lemma~3.1 of \citet{van2008rates}, for any measurable $f,g,$
	\begin{align}
		\begin{split}
			K(p_f,p_g)&\lesssim  \lVert f-g \rVert_\infty^2 e^{\lVert f-g \rVert_\infty}(1+\lVert f-g \rVert_\infty),\\
			V(p_f,p_g)&\lesssim  \lVert f-g \rVert_\infty^2 e^{\lVert f-g \rVert_\infty}(1+\lVert f-g \rVert_\infty)^2,\\
			\rho_{\rm H}(p_f,p_g)&\le \lVert f-g \rVert_\infty e^{\lVert f-g \rVert_\infty/2}.
		\end{split}
		\label{eqn:densitybo}
	\end{align}
(The uniform norm is used in \citet{van2008rates} but can be easily replaced by the $L_\infty$-norm.)
	The first two assertions imply that there exists $C_1>0$ such that $B_n\supseteq\{f:\lVert f-f_0\rVert_\infty\le C_1\epsilon_n \}$ if $\epsilon_n\rightarrow 0$. Hence we follow the calculation in the proof of Theorem~\ref{thm:nonreg} to conclude that there exists a constant $\bar c>0$ such that $\Pi(B_n)\ge e^{-\bar c n\epsilon_n^2 }$.
	The last assertion of \eqref{eqn:densitybo} enables us to work with the supremum norm in the calculation of the Hellinger covering number. The entropy calculation in Theorem~\ref{thm:nonreg} also verifies \eqref{eqn:con2-2}, completing the proof.
\end{myprooftitle}

\begin{myprooftitle}[Proof of Theorem~\ref{thm:binary}]
	Denote by $p_f(x,y)$ the density of model \eqref{eqn:binary} and by $p_0(x,y)$ the true density.
	We also write $f_0=H^{-1}(\varphi_0)$.
	From the fact that $|p_f(0|x)-p_0(0|x)|=|p_f(1|x)-p_0(1|x)|=|H(f(x))-H(f_0(x))|$, 
	it follows that	$\lVert p_f-p_0 \rVert_2 =\sqrt{2} \lVert H(f)-H(f_0) \rVert_{2,Q} $.
	The $L_2$-norm is bounded by a multiple of the Hellinger distance as $p_f$ and $p_0$ are uniformly bounded,  (see, for example, Lemma~B.1 of \citet{ghosal2017fundamentals}). Hence, it suffices to show the contraction rate results with respect to the Hellinger distance.
	This means that the assertion can be verified if there exists $\Theta_n\subseteq \mathcal F$ satisfying 
	\eqref{eqn:con1-2}--\eqref{eqn:con3-2} for some $\bar c>0$. 
	By Lemma~2.8 of \citet{ghosal2017fundamentals}, $K(p_0,p_f)\lesssim \lVert f-f_0 \rVert_{2,Q}^2$ and $V(p_0,p_f)\lesssim \lVert f-f_0 \rVert_{2,Q}^2$. We also have that $\rho_{\rm H}(p_f,p_g)\lesssim \lVert f-g \rVert_{2,Q}$ for every measurable $f,g$ by the same lemma.
	Similar to the proof of Theorem~\ref{thm:nonregrd}, the proof is completed by following that of Theorem~\ref{thm:nonreg}.
\end{myprooftitle}

\begin{myprooftitle}[Proof of Theorem~\ref{thm:addreg}]
	It suffices to verify \eqref{eqn:con1}--\eqref{eqn:con3} for the given model. 
	Following the proof of Theorem~\ref{thm:nonreg}, one can easily see that \eqref{eqn:con1} is verified as soon as
	\begin{align}
		\log\Pi(\widehat{\mathcal E})+\log\Pi(f\in\mathcal F_{\widehat{\mathcal E}}:\lVert f-f_0\rVert_n \le C_1\epsilon_n^\ast)\gtrsim -n(\epsilon_n^\ast)^2,
		\label{eqn:wts1}
	\end{align}
	for an approximating ensemble $\widehat{\mathcal E}$.
	Assumption \ref{asm:spnet} says that for each $1\le t\le T_0$, there exists a $\mathcal Z$-tree partition $\widehat{\mathcal T}^t$ such that $\lVert\hat f_{0t}-f_{0t}\rVert_n\lesssim \bar\epsilon_{t,n}$ for some $\hat f_{0t}\in\mathcal F_{\widehat{\mathcal T}^t}$.
	We index $\widehat{\mathcal E}=(\widehat{\mathcal T}^1,\dots,\widehat{\mathcal T}^T)$ with $\widehat{\mathcal T}^t=\{[0,1]^p\}$, $t=T_0+1,\dots, T$. Then,
	\begin{align*}
		\log \Pi(\widehat{\mathcal E})
		=\sum_{t=1}^{T_0}\log \Pi(\widehat{\mathcal T}^t)+(T-T_0) \log (1-\nu)
		\gtrsim  -\sum_{t=1}^{T_0}\widehat K^t \log n -\sum_{t=1}^{T_0} d_t\log p \gtrsim - n(\epsilon_n^\ast)^2,
	\end{align*}
	by Lemma~\ref{lmm:priorcon}. 
	Constructing $\hat f_{0t}$ as in the proof of Theorem~\ref{thm:approx2}, we denote every $\hat f_{0t}$ by $f_{0t,\widehat{\mathcal T}^t,\widehat\beta^t}$, where $\widehat\beta^t$ is the corresponding step-heights.
	Then the approximator of $f_0$ can be expressed as $f_{0,\widehat{\mathcal E},\widehat B}=\sum_{t=1}^{T_0} f_{0t,\widehat{\mathcal T}^t,\widehat\beta^t}$ with the ensemble components $(\widehat{\mathcal E},\widehat B)$, where $\widehat B=({\widehat\beta^{1\top}},\dots,{\widehat\beta^{T_0 \top}},0,\dots,0)^\top\in \mathbb R^{\widehat K_\ast}$ with $\widehat K_\ast=\sum_{t=1}^{T_0} \widehat K^t+T-T_0$. This gives us that
	\begin{align*}
		\lVert f-f_0\rVert_\infty\le \lVert f-f_{0,\widehat{\mathcal E},\widehat B}\rVert_\infty+\sum_{t=1}^{T_0} \lVert f_{0t,\widehat{\mathcal T}^t,\widehat\beta^t}-f_{0t}\rVert_\infty \lesssim \lVert f-f_{0,\widehat{\mathcal E},\widehat B}\rVert_\infty+ \sum_{t=1}^{T_0} \epsilon_{t,n}.
	\end{align*}
	Therefore, using $\sum_{t=1}^{T_0} \epsilon_{t,n} \le \sqrt{T_0} \epsilon_n^\ast $, we obtain that
	\begin{align}
		\Pi(f\in\mathcal F_{\widehat{\mathcal E}}:\lVert f-f_0\rVert_\infty\le C_1\epsilon_n^\ast)&\ge \Pi\left(f\in\mathcal F_{\widehat{\mathcal E}}:\lVert f-f_{0,\widehat{\mathcal E},\widehat B} \rVert_\infty\le C_2\epsilon_n^\ast\right).
		\label{eqn:poloba}
	\end{align}
	For any $B_1=({\beta_1^{1\top}},\dots, {\beta_1^{T_0\top}},\beta_1^{T_0+1},\dots,\beta_1^T)^\top,B_2=({\beta_2^{1\top}},\dots, {\beta_2^{T_0\top}},\beta_2^{T_0+1},\dots,\beta_2^T)^\top\in\mathbb R^{\widehat K_\ast}$, we write $f_{\widehat{\mathcal E},B_1}$, $f_{\widehat{\mathcal E},B_2}\in\mathcal F_{\widehat{\mathcal E}}$ to denote two additive tree functions that lie on the same partition ensemble $\widehat{\mathcal E}$. From \eqref{eqn:treelearner}, it is easy to see that $\lVert f_{\widehat{\mathcal E},B_1}-f_{\widehat{\mathcal E},B_2}\rVert_\infty\le \lVert B_1-B_2\rVert_2\widehat K_\ast^{1/2}$. As $\widehat K_\ast\log n\lesssim\sum_{t=1}^{T_0} \widehat K^t \log n\lesssim n(\epsilon_n^\ast)^2$, one can follow the proof of Theorem~\ref{thm:nonreg} to lower bound the logarithm of \eqref{eqn:poloba} by a constant multiple of $-n(\epsilon_n^\ast)^2$. Combined with the lower bound of $\Pi(\widehat{\mathcal E})$, this verifies \eqref{eqn:wts1}.
	The conditions in \eqref{eqn:con2} and \eqref{eqn:con3} follow directly from the proof of Theorem~\ref{thm:nonreg}, but with the rate $\epsilon_n^\ast$ for the additive regression.
\end{myprooftitle}

\section{Auxiliary Result: Dirichlet Prior Concentration}

The following lemma is a slight modification of Theorem 2.1 of \citet{yang2014minimax}. We provide the complete proof for a self-contained result. Similar results are also available in the literature (e.g., Lemma G.13 of \citet{ghosal2017fundamentals}).

\begin{mylemma}[Concentration of Dirichlet priors]
	Suppose that $\eta\in\mathbb S^p$ has a Dirichlet prior in \eqref{eqn:dir} with $\zeta>0$ and $\xi>1$. For any $\eta^\ast\in\mathbb S^p$ such that $\sum_{j=1}^p \mathbbm 1(\eta_j^\ast\ne 0)=s$ and any $\epsilon\in(0,1)$, there exists a constant $C>0$ such that
	\begin{align} 
		\Pi(\lVert\eta-\eta^\ast\rVert_1\le \epsilon)&\ge\exp\{-C\xi s \log (p/\epsilon)\},\label{eqn:dirineq1} \\
		\Pi\Bigg(\min_{S:|S|= s}\sum_{j\notin S}\eta_j \ge \epsilon\Bigg)&\le \exp\{-C(\xi-1) s \log p-\log \epsilon\}. \label{eqn:dirineq2}
	\end{align}
	\label{lmm:dir}
\end{mylemma}
\begin{myproof}
	We first prove \eqref{eqn:dirineq1}.
	Without loss of generality, we assume that the index set of nonzero entries of $\eta^\ast$ is $\{1,2,\dots,s-1,p\}$, i.e., $\eta_j^\ast=0$, $j=s,s+1,\dots, p-1$.
	By the inequality $|\eta_p-\eta_p^\ast|=|\sum_{j=1}^{p-1}\eta_j-\sum_{j=1}^{p-1}\eta_j^\ast|\le \sum_{j=1}^{p-1}|\eta_j-\eta_j^\ast|$, observe that $\lVert\eta-\eta^\ast\rVert_1\le 2\sum_{j=1}^{p-1} |\eta_j-\eta_j^\ast|= 2\sum_{j=1}^{s-1} |\eta_j-\eta_j^\ast|+ 2\sum_{j=s}^{p-1}\eta_j$. Hence, for $b_0=\epsilon/(4s)$ and $b_1=\epsilon/(4p-4s)$,
	\begin{align*}
		\mathcal S&=\{\eta\in\mathbb S^p:|\eta_j-\eta_j^\ast|\le b_0, j=1,\dots,s-1,\eta_j\in(0,b_1], j=s,\dots,p-1\}\\
		&\subseteq \{\eta\in\mathbb S^p:\lVert\eta-\eta^\ast\rVert_1\le\epsilon\}.
	\end{align*}
	Using this, we obtain
	\begin{align*}
		\Pi(\lVert\eta-\eta^\ast\rVert_1\le \epsilon)&\ge\Pi(\mathcal S)\\
		&=\int_{\mathcal S} \frac{\Gamma(\zeta/p^{\xi-1})}{\Gamma^p(\zeta/p^{\xi})} \prod_{j=1}^{p-1}\eta_j^{\zeta/p^\xi-1}\Bigg(1-\sum_{j=1}^{p-1} \eta_j\Bigg)^{\zeta/p^\xi-1}d\eta_1\dots d\eta_{p-1}\\
		& \ge \frac{\Gamma(\zeta/p^{\xi-1})}{\Gamma^p(\zeta/p^{\xi})} \left\{\prod_{j=1}^{s-1}  \int_{\max\{0,\eta_j^\ast-b_0\}}^{\min\{1,\eta_j^\ast+b_0\}} \eta_j^{\zeta/p^\xi-1} d\eta_j \right\} \left\{\prod_{j=s}^{p-1} \int_0^{b_1} \eta_j^{\zeta/p^\xi-1} d\eta_j \right\},
	\end{align*}
	where we used the fact that $\eta_p\le1$ and $\zeta/p^\xi-1<0$ for large enough $p$. As the Taylor expansion of $\Gamma$ gives that $x\Gamma(x)=1-\gamma_0 x+O(x^2)$ for the Euler-Mascheroni constant $\gamma_0$, we obtain $\Gamma(x)\asymp 1/x$ for every small enough $x$. Therefore, the last display is bounded below by a constant multiple of 
	\begin{align*}
		\frac{(\zeta/p^{\xi})^p}{\zeta/p^{\xi-1}} (2b_0)^{s-1} \left(\frac{p^\xi}{\zeta} b_1^{\zeta/p^\xi}\right)^{p-s} &= \zeta^{s-1} p^{-\xi(s-1)-1} \left(\frac{\epsilon}{2s}\right)^{s-1} \left(\frac{\epsilon}{4p-4s}\right)^{\zeta p^{-(\xi-1)}(1-s/p)} \\
		&\ge \zeta^{s-1} p^{-\xi(s-1)-1} \left(\frac{\epsilon}{2s}\right)^{s-1} \left(\frac{\epsilon}{4p}\right)^\zeta,
	\end{align*}
	where for the inequality we used the fact that $\xi \ge 1$. The logarithm of the rightmost side leads to the desired assertion.
	
	Now, we verify \eqref{eqn:dirineq2}. Consider a Dirichlet process $\text{DP}(\zeta/p^{\xi-1},Q_0)$ with concentration parameter $\zeta/p^{\xi-1}$ and uniform measure $Q_0$ on $[0,1]$. Suppose a random measure $P \sim \text{DP}(\zeta/p^{\xi-1},Q_0)$. Then, for the intervals $\mathcal I_j=[j-1)/p,j/p)$, $j=1,\dots,p$, we have 
	\begin{align*}
		(P(\mathcal I_1),\dots, P(\mathcal I_p)) \sim \text{Dir}(\zeta/p^{\xi},\dots,\zeta/p^{\xi}).
	\end{align*}
	This allows us to define $\eta$ as $\eta=(P(\mathcal I_1),\dots, P(\mathcal I_p))^\top$ using the Dirichlet process above. The stick-breaking representation of a Dirichlet process gives an expression
	$P=\sum_{k=1}^\infty w_k\delta_{z_k}$ for $z_k\sim Q_0$ and
	\begin{align*}
		w_k=v_k\prod_{j=1}^{k-1}(1-v_j), \quad v_k\sim \text{Beta}(1,\zeta/p^{\xi-1}).
	\end{align*}
	For every $k$, let $j_k$ be the index such that $z_k\in \mathcal I_{j_k}$. It follows that
	\begin{align*}
		\max_{S:|S|\le s}\sum_{j\in S}\eta_j \ge \sum_{k=1}^s\eta_{j_k} = \sum_{k=1}^s P(\mathcal I_{j_k}) = \sum_{1\le\ell<\infty:z_\ell\in\cup_{k=1}^s \mathcal I_{j_k}} w_\ell \ge \sum_{k=1}^s w_k,
	\end{align*}
	where the last inequality holds as $z_k\in \mathcal I_{j_k}$, $k=1,\dots, s$. This gives that
	\begin{align*}
		\min_{S:|S|= s}\sum_{j\notin S}\eta_j \le 1-\sum_{k=1}^s w_k =  1-\sum_{k=1}^s v_k\prod_{j=1}^{k-1}(1-v_j)=\prod_{j=1}^s(1-v_j),
	\end{align*}
	where the last equality can be verified by induction. Letting $\bar v_j=1-v_j\sim \text{Beta}(\zeta/p^{\xi-1},1)$, $j=1,\dots,s$, we obtain 
	\begin{align*}
		\Pi\Bigg(\min_{S:|S|= s}\sum_{j\notin S}\eta_j\ge \epsilon\Bigg)\le \Pi\Bigg( \prod_{j=1}^s \bar v_j\ge \epsilon \Bigg)\le 
		\frac{\zeta^s}{\epsilon(\zeta+p^{\xi-1})^s}\le \epsilon^{-1}\zeta^s p^{-s(\xi-1)},
	\end{align*}
	using the Markov inequality. The rightmost side verifies the assertion.
\end{myproof}

\vskip 0.2in
\bibliography{ref}

\end{document}